\documentclass[%
a4paper]{mpi2015-cscpreprint}

\usepackage{graphicx,epstopdf}
\usepackage{subfig}
\usepackage{enumitem} 
\usepackage{tikz}

\usepackage{breqn} 
\usepackage{algorithm}
\usepackage{algorithmic}
\usepackage{lipsum}
\usepackage{amsfonts}

\newtheorem{remark}{Remark}

\usepackage{amssymb}
\usepackage{bbm}
  
\newcommand{\td}{\tilde}
\newcommand{\mbb}{\mathbb}
\newcommand{\mcal}{\mathcal}
\newcommand{\parSp}{\mcal Z}		
\newcommand{\opn}{\operatorname}
\newcommand{\hsp}{\hspace{0.1cm}}
\newcommand{\hspB}{\hspace{0.3cm}}
\newcommand{\pd}{\partial}
\newcommand{\lan}{\left\langle}
\newcommand{\ran}{\right\rangle}

\newcommand{\zRef}{z_{\opn{ref}}}
\newcommand{\ztrain}[1]{z^{(#1)}}

\newcommand{\snapMat}[1]{\mcal S_{#1}}

\DeclareMathOperator*{\argmax}{arg\,max}
\DeclareMathOperator*{\argmin}{arg\,min}

\newcommand{\sol}[1]{u(\cdot,{#1})}
\newcommand{\solFV}[1]{u_N(\cdot,{#1})}
\newcommand{\solTransf}[2]{u(#1,#2)}
\newcommand{\solTransfFV}[2]{u_N(#1,#2)}

\newcommand{\solROM}[2]{u_{#1}(\cdot,#2)}

\newcommand{\spTransf}[2]{\varphi(\cdot,#1,#2)}
\newcommand{\dev}[2]{\Psi(\cdot,#1,#2)}

\newcommand{\coeffpod}[2]{\alpha_{#1}(#2)}
\newcommand{\coeffpodR}[2]{\alpha^{\mcal R}_{#1}(#2)}

\newcommand{\sampleTrain}{m_{\opn{tr}}}

\newcommand{\basisPOD}[1]{v_{#1}}

\newcommand{\nPsi}{n_{\tilde{\Psi}}}

\newcommand{\Eproj}[1]{E_n^{proj}(#1)}

\usetikzlibrary{decorations.pathreplacing}
\usetikzlibrary{decorations.shapes}
\usetikzlibrary{calc}

   \newtheorem{theorem}{Theorem}

   \newtheorem{definition}[theorem]{Definition}

\tikzset{decorate sep/.style 2 args=
{decorate,decoration={shape backgrounds,shape=circle,shape size=#1,shape sep=#2}}}

\usepackage{amsopn}

\begin{document}
\title{Data-Driven Model Order Reduction for Problems with Parameter-Dependent Jump-Discontinuities}
  
\author[$\ast$]{Neeraj Sarna}
\affil[$\ast$]{Max Planck Institute for Dynamics of Complex Technical Systems, 39106 Magdeburg, Germany.
  \email{sarna@mpi-magdeburg.mpg.de}, \orcid{0000-0003-0607-2067}}
  
\author[$\dagger$]{Peter Benner}
\affil[$\dagger$]{Max Planck Institute for Dynamics of Complex Technical Systems, 39106 Magdeburg, Germany, and Faculty for Mathematics, Otto von Guericke University Magdeburg, 39106 Magdeburg, Germany.\authorcr
  \email{benner@mpi-magdeburg.mpg.de}, \orcid{0000-0003-3362-4103}}
  
\shorttitle{Data-driven model order reduction jump-discontinuities}
\shortauthor{Neeraj Sarna, Peter Benner}
\shortdate{}
  


\abstract{We propose a data-driven model order reduction (MOR) technique for parametrized partial differential equations that exhibit parameter-dependent jump-discontinuities. Such problems have poor-approximability in a linear space and therefore, are challenging for standard MOR techniques. We build upon the methodology of approximating the map between the parameter domain and the expansion coefficients of the reduced basis via regression. The online stage queries the regression model for the expansion coefficients and recovers a reduced approximation for the solution. We propose to apply this technique to a transformed solution that results from composing the solution with a spatial transform. Unlike the (untransformed) solution, it is sufficiently regular along the parameter domain and thus, is well-approximable in a low-dimensional linear space. To recover an approximation for the (untransformed) solution, we propose an online efficient regression-based technique that approximates the inverse of the spatial transform. Our method features a decoupled online and offline stage, and benchmark problems involving hyperbolic and parabolic equations demonstrate its effectiveness.}

\novelty{
\begin{enumerate}
\item Combination of image registration and regression to approximate problems with parameter-dependent jump-discontinuities.
\item Image registration provides a transformed solution, which is well approximable in a linear reduced space.
\item Regression approximates the mapping from the parameter domain to the reduced basis coefficients of the transformation solution.
\item A de-transformation step recovers an approximation for the (untransformed) solution.
\item Numerical experiments report significant improvements over a standard linear approximation.
\end{enumerate}
}
\maketitle

\section{Introduction}
We consider parametrized partial differential equations (pPDEs) in a multi-query scenario where we seek a solution at multitudes of different parameter instances. Such scenarios arise in applications related to active control \cite{Benner2014}, design optimization \cite{Amsallem2015,Han2005,Yue2013}, uncertainty quantification \cite{Benner2015uncertainty,Crisovan2019model}, etc. Although a high-fidelity finite-difference/element/volume type solver can approximate the solution at any given parameter, it is prohibitively expensive for multi-query scenarios. We, therefore, resort to a MOR-based surrogate. 

The MOR technique splits the solution procedure into an offline and an online stage. The offline stage bears the cost of the most expensive computations and is decoupled from the online stage, allowing for online efficiency. 
Using the high-fidelity solver, it collects solution snapshots and computes a set of reduced basis vectors $\mcal X_n := \{\basisPOD{i}\}_{i=1,\dots,n}$. The online stage computes an approximation in the span of this basis. Precisely, let $\sol{z}$, defined over a spatial domain $\Omega$, represent a solution (or its high-fidelity approximation) to our pPDE. Furthermore, let $z\in\parSp$ be some parameter-of-interest, and let $\solROM{n}{z}\in \opn{span}(\mcal X_n)$ be an approximation to $\sol{z}$ such that 
\begin{gather}
\solROM{n}{z} := \sum_{i=1}^n\alpha_{u,i}(z)\basisPOD{i}.
\end{gather}
Then, the online stage computes the coefficients $\coeffpod{u,i}{z}$. We collect all these coefficients in a vector $\alpha_u(z)\in\mbb R^n$.

We consider the proper-orthogonal-decomposition (POD) approach to construct $\mcal X_n$. Since the POD is data-driven, the intrusive or non-intrusive nature of our MOR technique hinges on the methodology used to compute the vector $\coeffpod{u}{z}$. An intrusive technique computes $\coeffpod{u}{z}$ by projecting the pPDE onto the approximation space $\opn{span}(\mcal X_n)$. In doing so, it accesses the discrete (high-fidelity) evolution operators---the books \cite{JanStamm,PeterBook} and the review paper \cite{PeterReview} discuss this approach at length. In contrast to the intrusive approach, a non-intrusive approach treats the high-fidelity solver as a black-box. It computes $\coeffpod{u}{z}$ either via regression or via a data-driven approximation of its evolution operator. 

Broadly speaking, non-intrusive techniques are either physics or non-physics informed. The former requires structural information of the non-linearities in the underlying pPDE. This information is then used to infer the evolution operator of $\coeffpod{u}{z}$ \cite{Kramer2020,Benjamin2016,Qian2020,NonlinearDMD}. In contrast, the latter “ignores” the underlying pPDE altogether. It treats the solution $\sol{z}$ like a generic parametrized function that might as well not even correspond to the solution of a pPDE.

The non-physics informed methods either (i) directly approximate the vector $\coeffpod{u}{z}$ using regression; or (ii) first approximate $\coeffpod{u}{z}$'s evolution operator---without using the structural properties of the non-linearities---followed by time-stepping. Authors in \cite{NairNonIntrusive2013,GreedyJan2018}, \cite{Guo2019} and \cite{JanNN2018} undertake the first approach and use radial basis functions (RBFs), Gaussian processes (GP) and neural networks, respectively, for regression. Authors in \cite{RBFXiao2015,Xiao2015} opt for the second approach and reduce the Navier--Stokes equations by approximating the evolution operator using RBFs---\cite{Xiao2016} and \cite{MovingBC2021} provide extensions to problems involving fluid-structure interaction and moving boundaries, respectively. 

\subsection{Regression based non-intrusive MOR}
We consider a non-physics informed technique that approximates the map between the parameter and the POD coefficients using regression. In particular, consider a set of parametrized functions given as 
$
\{g(\cdot,z)\hsp :\hsp z\in\mcal Z\}.
$
The function $g(\cdot,z)$ can correspond to a solution of a pPDE, a transformed solution of a pPDE (discussed below), or any other quantity-of-interest. \Cref{algo: offline block} and \Cref{algo: online block} outline the building blocks of an algorithm that can approximate any function of the above set in a data-driven fashion. The main ingredient of the algorithm is a regression model for the mapping $z\mapsto \coeffpod{g}{z}$, which the online phase queries to recover an approximation for $g(\cdot,z)$. 

We use GPR to perform regression. The main reason being that we also develop a regression-based surrogate for the error introduced by our MOR technique. Since a GPR is probablistic, its confidence region is helpful in devising an accurate error surrogate. However, if one is not interested in error modelling then, any other regression technique would also suffice---RBFs being a noteworthy example.

\begin{algorithm}[ht!]
\caption{Regression based MOR: Offline stage}
\begin{algorithmic}[1] \label{algo: offline block}
\STATE Collect the training samples $\{\ztrain{i}\}_i$ and the snapshots $\{g(\cdot,\ztrain{i})\}_i$.
\STATE Compute the reduced basis $\mcal X_n:=\{v_i\}_{i=1,\dots,n}$.
\STATE Orthogonally project $\{g(\cdot,\ztrain{i})\}_i$ onto $\{v_i\}_{i=1,\dots,n}$ to compute the training data $\{(\ztrain{i},\coeffpod{g}{\ztrain{i}}\}_i$. 
\STATE Using regression, approximate $z\mapsto \coeffpod{g}{z}$ by $z\mapsto \coeffpodR{g}{z}$.
\end{algorithmic}
\end{algorithm}

\begin{algorithm}[ht!]
\caption{Regression based MOR: Online stage}
\begin{algorithmic}[1] \label{algo: online block}
\STATE For any $z\in \mcal Z$, query the regressed mapping $z\mapsto\coeffpodR{g}{z}$.
\STATE Approximate $g(\cdot,z)$ by $g_n(\cdot,z) := \sum_i\alpha^{\mcal R}_{g,i}(z)\basisPOD{i}$. 
\end{algorithmic}
\end{algorithm}

\subsubsection{Shortcomings of the standard approach}We consider pPDEs that exhibit parameter-dependent jump-discontinuities or steep-gradients. Several problems of practical interest exhibit such a behaviour. A standard example being that of non-linear hyperbolic equations, which mostly appear in applications involving fluid flows. Even for a smooth initial and boundary data, such problems can develop spatial discontinuities that move in the parameter domain \cite{Cagniart2019,Welper2017}. A diffusion equation whose diffusion coefficient has a moving discontinuity, also belongs to a similar category \cite{RegisterMOR}.

In the above algorithms, one can choose 
\begin{gather}
g(\cdot,z) = \sol{z}
\end{gather}
and recover a non-intrusive MOR technique \cite{NairNonIntrusive2013,Guo2019,JanNN2018}. However, such an approach is inefficient for the aforementioned problems. Since a solution with parameter-dependent jump-discontinuities (or steep-gradients) has poor approximability in a linear space, only a large set of POD modes can provide a reasonable accuracy. This makes the MOR technique inefficient---see \cite{PeterBook,Welper2017,WaveKolmogorov} for related proofs. Note that in the context of POD, a slow singular value decay is indicative of the poor approximability in POD basis.

It is noteworthy that even a sufficiently large set of POD modes does not guarantee a physically accurate and stable solution. 
Discontinuities in the parameter domain can trigger oscillations in the POD modes. Analogous to the Gibbs phenomenon for a Fourier series expansion, the oscillation frequency increases with the order of the POD modes. Therefore, although increasing the number of POD modes can provide a better approximation in the $L^2$(or $L^1$)-sense, it results in a solution with un-physical high frequency oscillations---results in \cite{Xiao2016,Kevin2015adaptive,KevinAuto} showcase these oscillations. 

\subsubsection{A transformation and de-transformation approach}
For the above reasons, we refrain from directly approximating the solution. Rather, we undertake a two step procedure comprising of a solution transformation followed by de-transformation. In the transformation step, we introduce a spatial transform 
\begin{gather}
\spTransf{\zRef}{z}:\Omega\to\Omega,
\end{gather}
with $\zRef\in\parSp$ being a reference parameter. This spatial transform is such that the transformed solution given by  
\begin{gather}
g(\cdot,z) = \solTransf{\spTransf{\zRef}{z}}{z}, \label{trasf sol}
\end{gather}
at least ideally, has no discontinuities along $\parSp$. This ensures that the transformed solution (with some additional assumptions) is well approximable in a sufficiently low-dimensional linear space---we refer to \cite{Welper2017,WelperAdaptive,SarnaCalib2020} for further details. One may associate a physical relevance to $\varphi$ by interpreting it as a transformation to a Lagrangian coordinate system where the discontinuities do not move in $\parSp$ \cite{MojganiLagrangian,RegisterMOR}. This robs $g(\cdot,z)$ of its \textit{transport-dominant} nature and  makes it well-approximable in a low-dimensional linear space.

The algorithms discussed earlier provide the data-driven approximation 
$g(\cdot,z) \approx g_n(\cdot,z),$ where $g_n(\cdot,z)$ belongs to the POD space constructed using the snapshots of $g(\cdot,z)$. Let us recall that our goal was to approximate the solution $\sol{z}$. Therefore, we need to circle back and recover an approximation for $\sol{z}$ from $g_n(\cdot,z)$. This is what we refer to as de-transformation. De-transformation requires an inversion of $\spTransf{\zRef}{z}$. Since an exact inversion is prohibitively expensive, we propose an efficient GPR-based technique to approximate $\spTransf{\zRef}{z}^{-1}$. Analogues to the approximation for $g(\cdot,z)$, our technique first constructs POD basis for $\spTransf{\zRef}{z}^{-1}$ and then uses GPR to approximate the POD coefficients.

\subsection{Intrusive vs.  non-intrusive approach}
Comments that motivate a non-intrusive approach are in order. We particularly emphasize on the first two points, which, we believe, are exclusive to problems that exhibit parameter-dependent jump-discontinuities. 
\begin{enumerate}
\item  Firstly, an intrusive approach derives a Lagrangian pPDE for $\solTransf{\spTransf{\zRef}{z}}{z}$ and reduces it via a Galerkin projection over the POD modes \cite{RegisterMOR,MojganiLagrangian}. Even for an affine-in-parameter pPDE, the Lagrange equations can be non-affine. This makes the reduced-order model as expensive as the high-fidelity solver. To gain efficiency, one resorts to hyper-reduction, which adds a layer of approximation and complexity to the MOR technique. In contrast, decoupled from the underlying pPDE, the non-intrusive technique treats the affine and non-affine parameter dependence alike---similar comments holds for non-linear pPDEs. Note that rather than transforming to Lagrangian coordinates, one can directly reduce the pPDE in Eulerian coordinates \cite{Nair,MATS}. In Eulerian coordinates, the reduced approximation space is non-linear, and gaining efficiency requires sophisticated non-standard hyper-reduction techniques. 

\item Secondly, the Lagrangian pPDE contains the inverse of the derivatives of $\spTransf{\zRef}{z}$. If ill-conditioned, we speculate, these terms can result in stability issues with the Galerkin projection. 
\item Lastly, the underlying pPDE might be unavailable or the discrete finite-difference/element operators might be inaccessible. The first scenario, for instance, corresponds to solution snapshots collected from experiments, and the latter corresponds to legacy codes or commercial solvers that provide only the solution and not the discrete evolution operators.
\end{enumerate}

\subsection{Relation to previous works}
Particularly in the context of hyperbolic pPDEs with moving discontinuities, non-intrusive MOR techniques are not new in the literature. Authors in \cite{Metric2019} embed the solution manifold in the Wasserstein space and approximate the solution using Wasserstein barycentres. The scheme works well for conservative hyperbolic pPDE. However, extensions to multi-dimensional spatial domains and non-conservative pPDEs with boundary conditions are unavailable, as of yet. Closer to our approach is the transformed-snapshot-interpolation (TSI) proposed in \cite{Welper2017}. TSI is an Eulerian method that approximates $\sol{z}$ using Lagrange polynomial interpolation over the transformed snapshots. The sample parameters must lie on a tensorized grid over $\mcal Z$. In contrast, our method is Lagrangian and allows for a POD based approximation. The sample parameters need not be tensorized and can result from any of the sampling techniques summarized in \cite{RBBook}. This is particularly appealing for high-dimensional parameter domains where tensorized grids result in a large number of parameter samples, making the offline step unaffordable. Furthermore, a GPR is probabilistic and provides a confidence region that quantifies the quality of regression. Such a quantification can either be used to increase the size of the training set \cite{GPR}, to quantify the extrapolation capabilities \cite{MovingBC2021}, or to develop an error model \cite{Guo2019}. A Lagrange polynomial interpolation is deterministic and does not offer such flexibility. 

The novelty of our work is in combining the regression based MOR techniques developed earlier in \cite{JanNN2018,Guo2019} with two additional steps: transformation and de-transformation. Thereby, we extend the validity of these techniques to pPDEs that exhibit parameter-dependent jump-discontinuities. We emphasize that these additional steps do not interfere with a pre-existing numerical implementation of the earlier mentioned two algorithms---recall that these algorithms were the building blocks of the MOR techniques developed in \cite{JanNN2018,Guo2019}. Therefore, having implemented the (de-)transformation step, a pre-existing numerical implementation can be smoothly extended to accommodate parameter-dependent jumps. Furthermore, as we clarify later, both of these steps can be easily deactivated for problems that do no exhibit parameter-dependent jumps or steep gradients and thus, the extra cost associated with (de-)transformation can be avoided for such problems. These problems can be identified by first passing the solution snapshots through a shock detector \cite{Limiter1984}. Our recommendation is to deactivate the above two steps in case the shock detector returns an empty set. 

Our transformation step is inspired by the registration-based MOR technique developed in \cite{RegisterMOR,Taddei2020ST}. Indeed, to compute the transform $\varphi$, we use the same image registration technique as that developed in \cite{RegisterMOR}. Nonetheless, there are some key differences that we outline as follows. Firstly, our technique is non-intrusive as opposed to the intrusive technique developed in \cite{RegisterMOR,Taddei2020ST}. Secondly, and most importantly, we also perform the de-transformation step to recover an approximation for the (untransformed) solution. To the best of our knowledge, authors in \cite{RegisterMOR,Taddei2020ST}, cater to developing an accurate approximation for the transformed solution and do not perform any de-transformation. The $L^2$-error in approximating the (untransformed) solution is computed using the Jacobians of $\varphi$, which does not require an explicit computation of $\varphi^{-1}$. In our framework, we explicitly approximate $\varphi^{-1}$ and recover an (explicit) approximation for the solution. We acknowledge that the possibility of de-transformation, and the associated difficulties, were discussed in Remark-3.3 of \cite{RegisterMOR}. Our de-transformation step is a possible solution to the problem posed in this remark. 

\subsection{Organization}
Rest of the article is organized as follows. \Cref{sec: GPR} provides a brief summary of the GPR employed. \Cref{sec: approx transf sol} outlines a data-driven technique to approximate the transformed solution described above. \Cref{sec: approx u} outlines the de-transformation step and presents an efficient computation of $\spTransf{\zRef}{z}^{-1}$. \Cref{sec: summary} presents a summary of our technique. \Cref{sec: num results} presents numerical results, and \Cref{sec: conclusion} closes the article with a conclusion.

\section{Gaussian process regression (GPR)}\label{sec: GPR}
The forthcoming sections will extensively use GPR and for completeness, we summarize it here. We refer the reader to the book \cite{GPR} for an exhaustive discussion. GPR solves the following regression problem: Given the training data
$\{(x_i,h(x_i))\}_{i=1,\dots,\sampleTrain},$ 
where $x_i\in\mbb R$ and $h:\mbb R\to\mbb R$, approximate $h(x^*)$, where $x^*\not\in \{x_i\}_i$. For convenience, we collect the training points in the set 
\begin{gather}
\mcal D_x:=\{x_i\}_i,
\end{gather}
and denote the approximation via
\begin{gather}
h(x^*)\approx\mcal R_\lambda(\theta,\mcal D_x,x^*),
\end{gather}
where $\mcal R_\lambda(\theta,\mcal D_x,x^*)$ is the GPR model, $\theta$ is a hyper-parameter, and $\lambda$ is a user-defined parameter---below, we clarify the definition of $\theta$ and $\lambda$. 

The GPR model follows from a two-step online-offline decomposition based procedure:
\begin{enumerate}
\item Offline, the training phase  computes $\theta$ using the training data;
\item Online, the prediction phase assigns a value to the GPR model at a given $x^*\in\mbb R$.
\end{enumerate}
We start with the details of the training phase.
 \subsection{Training a GPR}\label{sec: GPR train}
Training a regression model corresponds to computing its hyper-parameters that we collect in the vector $\theta\in\mbb R^r.$
The value of $r$ and the type of $\theta$ is regression technique-dependent. To present a GPR's hyper-parameters, we first define a few objects.  A GPR models $h(x)$ as a Gaussian process (GP) meaning that for any $m\in\mbb N$, the vector $(h(x_1),\dots,h(x_{m}))^T$ is normally distributed. A precise definition of a GP is as follows.
\begin{definition}[Gaussian process (GP)]\label{def: GP}
For some $m\in\mbb N$, let $\mcal D = \{x_i\}_i$ be a set of $m$ points in $\mbb R$. Let $\vartheta :\mbb R\to\mbb R$ be a mean function, and let $\zeta:\mbb R\times\mbb R\to\mbb R$ be a positive definite kernel function. Let $h:\mbb R\to\mbb R$ be a random function, and let $h_{\mcal D}\in \mbb R^m$ be a random vector such that 
\begin{gather}
\left(h_{\mcal D}\right)_i := h(x_i),\hspB\forall i\in\{1,\dots,m\}.
\end{gather}
Then, $h(x)$ is a Gaussian process if $h_{\mcal D}$ follows a multivariate normal distribution with the mean vector $\vartheta_{\mcal D}\in\mbb R^{m}$ and the covariance matrix $\zeta_{\mcal D\mcal D}\in\mbb R^{m\times m}$ given as 
\begin{gather}
\vartheta_{\mcal D} := (\vartheta(x_1),\dots,\vartheta (x_m))^T,\hspB \left(\zeta_{\mcal D\mcal D}\right)_{ij} := \zeta(x_i,x_j),\hspB \forall i,j\in\{1,\dots,m\}.
\end{gather}
\end{definition}
In practise, to well-define GPR, we restrict to a parametrized sub-class of mean and kernel functions. For simplicity, we do not introduce any new notations and denote a mean function and a kernel function of this sub-class via $\vartheta(\cdot,\beta)$ and $\zeta(\cdot,\cdot,\kappa)$, respectively, where $\beta$ and $\kappa$ are the parameters.
A first-order polynomial in $\mbb R$ is our mean function $\vartheta(\cdot,\beta)$, and the kernel function $\zeta(\cdot,\cdot,\kappa)$ is the automatic relevance determination (ARD) squared exponential (SE) kernel. Explicit forms read
\begin{gather}
\vartheta(x,\beta) = \beta^T \Phi(x),\hspB \zeta(x,y,\kappa) = \kappa^2_1 \exp\left(-\frac{1}{2}\frac{(x-y)^2}{\kappa_2^2}\right). \label{mu and kernel}
\end{gather}
The above choice of the kernel function works well for a function $h(x)$ that is smooth \cite{GPRReview}---for our applications, the smoothness property indeed holds true.
The function $\Phi(x) := (1,x)^T$ maps the sample space to the feature space. The vectors $\beta\in\mbb R^2$ and $\kappa := (\kappa_1,\kappa_2)^T$ together form the hyper-parameter given as
\begin{gather}
\theta := (\beta^T,\kappa^T)^T.
\end{gather}

To estimate $\theta$, we consider the Bayesian approach of maximum likelihood estimation (MLE).
We compute $\theta$ such that the joint normal distribution corresponding to the random vector  $h_{\mcal D_x}$ has the maximum possible log-likelihood. This approach, along with the fact that $h_{\mcal D_x}\sim \mcal N(\vartheta_{\mcal D_x},\zeta_{\mcal D_x\mcal D_x})$, leads to the maximization problem
\begin{equation}
\begin{aligned}
\theta = \argmax_{\theta^*\in\mbb R^4} \left(\right.&-\frac{1}{2}(h_{\mcal D_x}-(\beta^*)^T \Phi_{\mcal D_x})^T\zeta_{\mcal D_x\mcal D_x}(\kappa^*)^{-1}(h_{\mcal D_x}-(\beta^*)^T  \Phi_{\mcal D_x})\\
&\left.-\frac{1}{2}\log|\zeta_{\mcal D_x\mcal D_x}(\kappa^*)|\right).
\end{aligned}
\end{equation}
We train the GPR using the \texttt{fitgrp} routine in MATLAB.

\subsection{Prediction with GPR} The GPR model $\mcal R_\lambda(\theta,\mcal D_x,x^*)$ follows from the conditional probability $h(x^*)\rvert h_{\mcal D_x}$ i.e., the probability of $h(x^*)$ given the observed data $h_{\mcal D_x}$. Since $h(x)$ is a GP, the conditional probability $h(x^*)\rvert h_{\mcal D_x}$ is normally distributed with the mean-value $\vartheta^*(\theta,\mcal D_x,x^*)\in\mbb R$ and the covariance $\zeta^*(\theta,\mcal D_x,x^*)\in \mbb R$ given as 
\begin{equation}
\begin{gathered}
\vartheta^*(\theta,\mcal D_x,x^*) := \vartheta(x^*,\beta) + \zeta_{x^*\mcal D_x}(\kappa) \zeta_{\mcal D_x\mcal D_x}(\kappa)^{-1}(h_{\mcal D_x}-\vartheta_{\mcal D_x}(\beta)),\\
 \zeta^*(\theta,\mcal D_x,x^*) := \zeta_{\mcal D_x x^*}(\kappa) - \zeta_{x^*\mcal D_x}(\kappa)\zeta_{\mcal D_x\mcal D_x}(\kappa)^{-1}\zeta_{\mcal D_x x^*}(\kappa).
\end{gathered}
\end{equation}
Recall that $\vartheta$ is the mean function given in \eqref{mu and kernel}, and $\zeta_{\mcal D_x\mcal D_x}$ is the covariance matrix corresponding to the kernel function. We set 
\begin{gather}
\mcal R_{\lambda}(\theta,\mcal D_x,x^*) = \vartheta^*(\theta,\mcal D_x,x^*) + \lambda \sqrt{\zeta^*(\theta,\mcal D_x,x^*)}. \label{def R}
\end{gather}
Thus, $\mcal R_{\lambda}(\theta,\mcal D_x,x^*)$ is $\lambda$ standard-deviations (given by $\sqrt{\zeta^*(x^*,\mcal D_x)}$) away from the mean. For $\lambda = 0$, we recover the so-called mean-value prediction. We use the \texttt{predict} function from MATLAB to compute $\mcal R_{\lambda}(\theta,\mcal D_x,x^*)$.

\section{Approximation of the transformed solution}\label{sec: approx transf sol}
\subsection{Problem description}
Consider a general pPDE of the form  
\begin{gather}
\mcal L(u(x,z),z) = 0\hspB\forall (x,z)\in\Omega\times\parSp,
\end{gather}
where $\mcal L(\cdot,z)$ is some spatio-temporal differential operator, and $(x,z)\mapsto u(x,z) \in \mbb R$ is the solution. For simplicity of notation, we consider a scalar-valued solution. For a vector-valued $\sol{z}$, the same technique applies to each of its components. We assume that for all $z\in\parSp$, $\sol{z}\in L^2(\Omega)$. The set $\parSp \subset \mbb R^p$
 is the parameter-domain, which can (and for the later test cases will) include the temporal domain. We assume that $\parSp$ is (or could be mapped via a bijection to) a hyper-cube. It is noteworthy that $\mcal L(\cdot,z)$ can be non-linear or $\mcal L(u(x,z),\cdot)$ could be non-affine. Our technique treats linear, non-linear, affine and non-affine problems alike.

We restrict to a square spatial domain i.e., $\Omega = (0,1)^2$. The only (major) complexity introduced by a general curved domain is in the computation of the spatial transform $\varphi$ given in \eqref{trasf sol}. Other than that, the entire technique remains the same---we refer to \cite{CurvedDomains} for the computation of $\varphi$ on curved domains. As for now, we refrain from introducing this additional complexity and study the performance of our method on a unit square.  

We denote a high-fidelity approximation of the above pPDE via
\begin{gather}
\sol{z} \approx \solFV{z}\in\mcal X_N\subset L^2(\Omega),
\end{gather}
where $\mcal X_N$ is the high-fidelity approximation space with $\opn{dim}(\mcal X_N) = N$. Usually, $\mcal X_N$ is a finite-element/volume type space. We define $\mcal X_N$ over a shape-regular discretization of $\Omega$ defined as 
\begin{gather}
\Omega = \bigcup_{i=1}^N\mcal I_i,\label{space grid}
\end{gather}
where $\mcal I_i$ represents the $i$-th spatial cell. Note that for the sake of notational simplicity and consistency with our numerical experiments, the number of grid cells equals the dimensionality of $\mcal X_N$; in general, these two numbers can also be different. With $\Pi:L^2(\Omega)\to\mcal X_N$ we denote the orthogonal projection operator. 

We consider problems where, at least for some $x\in\Omega$, the function $u(x,\cdot)$ has jump-discontinuities or steep-gradients.
As explained earlier, for such problems, we first apply the MOR technique (outlined in \Cref{algo: offline block} and \Cref{algo: online block}) to a transformed solution given as
\begin{gather}
g(\cdot,z) = \Pi \solTransfFV{\spTransf{\zRef}{z}}{z}, \label{def g}
\end{gather}
where $\spTransf{\zRef}{z}:\Omega\to\Omega$ is the spatial transform, and $\zRef\in\parSp$ is the reference parameter, the choice of which will be discussed later. We recall that the spatial transform is such that the transformed solution (at least ideally) is not discontinuous along the parameter domain. This allows (along with some additional regularity assumption) for an accurate approximation in a sufficiently low-dimensional linear reduced space---further details can be found in \cite{Welper2017,SarnaCalib2020}. The following discussion elaborates on the different steps involved in a data-driven approximation of the transformed solution. 
\subsection{Snapshots of the transformed solution}\label{sec: snap transformed}
 We collect $\sampleTrain\in\mbb N$ different parameter samples in a set denoted by
\begin{gather}
\mcal D_z : =\{\ztrain{i}\}_{i=1,\dots,\sampleTrain}.\label{def D}
\end{gather}
These samples can either be chosen uniformly, randomly, or with the Lattice hyper-cube sampling technique given in \cite{SamplingLattice}. At all of these parameter samples, we need snapshots of $g(\cdot,z)$. This entails computing snapshots of the high-fidelity solution $\{\solFV{\ztrain{i}}\}_i$  and the spatial transform $\{\spTransf{\zRef}{\ztrain{i}}\}_i$. As stated earlier, the former we compute in $\mcal X_N$. For the latter, we use an optimization-based image registration technique summarized below. We refer to the review paper \cite{ReviewDeformReg} and the article \cite{RegisterMOR} for an image analysis and a MOR perspective, respectively, on the registration technique.

\subsubsection{Snapshots of $\varphi$}\label{sec: varphi} We express $\spTransf{\zRef}{z}$ as 
\begin{gather}
\spTransf{\zRef}{z} = \opn{Id} + \dev{\zRef}{z}.
\end{gather}
 The function $\Psi(\cdot,\zRef,z)$ is referred to as the displacement field---it displaces a point $x$ by a distance of $\Psi(x,\zRef,z)$. Furthermore, the identity mapping is denoted by $\opn{Id}$. \Cref{remark: diffeo phi} below motivates the above splitting. 

We seek a $\dev{\zRef}{z}$ that lies in a span of polynomials $\mcal P_M:=\mbb P_M\times \mbb P_M$ that reads
\begin{gather}
\mbb P_M:=\opn{span}\{L_{ij}\Upsilon\}_{i,j=1,\dots,M}. \label{def PM}
\end{gather}
The function $L_{ij}(x)$ is a product of the Legendre polynomials $l_i(x_1)$ and $l_j(x_2)$, where $i$ and $j$ denote the degrees of the respective Legendre polynomials. Thus, $\mcal P_M$ a $2M^2$-dimensional space. Furthermore, the function $\Upsilon(x):=\prod_{k=1}^2 x_k(1-x_k)$ ensures that for all $x\in\pd\Omega$, we have the boundary conditions $\Psi(\pd\Omega,\zRef,z) = 0$. Note that following \cite{Nair,Welper2017}, we have imposed a stricter set of boundary conditions than in \cite{RegisterMOR}. At least for the test cases we considered, these boundary conditions provide reasonable results without any additional constraints on the Jacobian of $\varphi$. We will later study these Jacobians empirically.

The expansion coefficients for $\dev{\zRef}{z}$ result from the minimization problem
\begin{equation}
\begin{aligned}
\dev{\zRef}{z}=\argmin_{\Psi^*\in \opn{span}(\mcal P_M)}\mcal F(\Psi^*,\zRef,z).  \label{L2 opt}
\end{aligned}
\end{equation}
where
\begin{gather}
\mcal F(\Psi^*,\zRef,z) := \mcal M(\Psi^*,z,\zRef)^2 + \epsilon\times \mcal R(\Psi^*)^2,\label{def F}
\end{gather}
with $\epsilon \geq 0$ being a penalty parameter. We minimize a summation of two objects, the so-called matching criterion $\mcal M$ and the regularizer $\mcal R$. The matching criterion (a term we borrow from image analysis \cite{ReviewDeformReg}) should be chosen such that the transformed solution $g(\cdot,z)$ is sufficiently regular along the parameter domain, making it well approximable in a linear reduced space. An appropriate choice of the matching criterion requires some information of the underlying physical process. For instance, in case of hyperbolic equations, $\mcal M$ can be chosen as the $L^2(\Omega)$ distance between the transformed and the reference snapshots \cite{Welper2017,Nair,RegisterMOR,Taddei2020ST}. \Cref{sec: num results} further elaborates on the choice of $\mcal M$. Note that the term matching criterion is justified for $\mcal M$ in the sense that minimizing $\mcal M(\Psi^*,\zRef,z)$ \textit{matches} the discontinuities between snapshots, which induces regularity in the transformed solution $g(\cdot,z)$. 

The regularizer $\mcal R(\Psi^*)$ penalizes the spatial regularity of $\Psi^*$. Thereby, preventing spurious oscillations and promoting a diffeomorphic $\spTransf{\zRef}{z}$. Several previous works deem the diffeomorphism property desirable \cite{Welper2017,RegisterMOR,MATS}---\Cref{remark: diffeo phi} below provides further elaboration. Furthermore, empirically, one observes that the optimization routine used to compute the above problem provides a more stable and better solution with spatial regularization \cite{Klein2007}. Here, we make the standard choice \cite{ReviewDeformReg}
\begin{gather}
\mcal R(\Psi^*) := \|\Delta \Psi^*\|_{L^2(\Omega;\mbb R^2)},
\end{gather}
where $\Delta $ represents the Laplace operator.

Despite the regularization term, computing the above minimization problem is a challenging task, with the (probable) non-convexity of the objective functional being a major concern. Owing to the limited scope of this article, we refrain from devising a specialized optimization toolbox for the above problem. Rather, following the works in \cite{Nair,RegisterMOR}, we resort to the standard interior point algorithm implemented in the \texttt{fmincon} routine of MATLAB. At least for the test cases presented later, this routine provides reasonable results with all the parameter values set to their default. 

Following comments cater to the several practical considerations one encounters while solving the above problem.
\begin{enumerate}
\item \textbf{Initial guess}: An initial guess of $\Psi^*=0$ can be far-off the desired solution, resulting in inaccuracies---due to the possible non-convexity of the objective functional, the optimization algorithm can get stuck in an undesirable local-minimum. Therefore, we first solve a few sub (optimization) problems and use their solution as an initial guess. Section-3.1.2 of \cite{RegisterMOR} presents the details, which are not repeated here for brevity.

\item \textbf{Choice of $\zRef$:} The reference parameter $\zRef$ is largely determined by the \textit{discontinuity-topology} i.e., the number and the relative orientation of the discontinuities \cite{Welper2017,WelperAdaptive}. 
With a parameter invariant discontinuity-topology, in theory, any $\zRef\in\parSp$ suffices. In case the topology changes along $\parSp$, one can either: (i) partition $\parSp$ such that the topology is preserved on each of the subsets \cite{SarnaCalib2020}; or (ii) choose multiple reference parameters and optimize over both $\dev{\zRef}{z}$ and the span of reference snapshots \cite{Taddei2020ST}. Following \cite{Nair,Welper2017}, for now, we restrict to the examples where the topology is parameter invariant, and choose $\zRef$ as the center of $\parSp$. Studies in \cite{RegisterMOR} indicate that the result of the registration technique might change with $\zRef$. However, trying to optimize over $\zRef$ introduces additional complexity to the numerical scheme; therefore, for now, we fix a value for $\zRef$ and accept the results our choice provides.

\item \textbf{Choice of $\epsilon$: } Intuitively, it seems reasonable to choose a non-zero $\epsilon$ smaller than one. Otherwise, we will heavily penalize the regularity of $\dev{\zRef}{z}$ at the expense of an accurate solution transformation. Empirically, we observe that with minor differences, all the different $\epsilon$ in the set $\{10^{-4},10^{-3},10^{-2}\}$ provide reasonable results. In our numerical experiments we set $\epsilon = 10^{-2}$. 

\item \textbf{Choice of $M$: } We choose $M$ iteratively. Corresponding to $\dev{\zRef}{z}\in \opn{span}(\mcal P_M)$, consider the parameter-average of the $\mcal M$-mismatch between the transformed and the reference solution defined as 
\begin{gather}
\Xi_M := \frac{1}{\sampleTrain}\sum_{z\in\mcal D_z}\mcal M(\dev{\zRef}{z},z,\zRef).
\end{gather} Starting with an initial guess of $M=1$, we continue to increase $M$ till, for some user-defined $\texttt{TOL}_M$, we satisfy
\begin{gather}
\frac{|\Xi_M - \Xi_{M-1}|}{\Xi_{M-1}} \leq \texttt{TOL}_M. \label{def TOL M}
\end{gather}
We set $\texttt{TOL}_M = 10^{-3}$. In our experience, decreasing $\texttt{TOL}_M$ further offered minuscule improvements at an additional computational cost. 
\end{enumerate}

\begin{remark}[Diffeomorphic $\varphi$] \label{remark: diffeo phi}
Instead of directly approximating $\spTransf{\zRef}{z}$, we approximate the displacement field $\Psi(\cdot,\zRef,z)$ in the polynomial space $\mcal P_M$. The reason being that in case $\Psi(\cdot,\zRef,z)$ is small as compared to $\opn{Id}$, we can expect a diffeomorphic $\spTransf{\zRef}{z}$, which, as detailed in \cite{Welper2017,RegisterMOR,SarnaCalib2020}, is a desirable property. In the image registration literature, such a displacement is referred to as a \textit{small-displacement.}

At least the test cases considered later exhibit small-displacements. Obviously, in general, a small-displacement is not guaranteed and therefore, the above technique does not guarantee a diffeomorphic $\varphi(\cdot,\zRef,z)$. This is a limitation of the technique that might be resolved by instead approximating the velocity field induced by $\varphi$---see \cite{LDDM}. 
\end{remark}

\begin{remark}[Recovery of the standard approach]
Note that by choosing 
\begin{gather}
\mcal M(\Psi^*,z,\zRef), \mcal R(\Psi^*) = 0,
\end{gather}
we find that $\spTransf{\zRef}{z}=\opn{Id}$, for all $z\in\parSp$. Consequently, we recover the regression based MOR technique developed in \cite{Guo2019}. This justifies our earlier claim that the transformation and the de-transformation step can be easily deactivated for problems that do not exhibit parametric jump-discontinuities. Such problems can be identified by first passing the snapshots through a shock (or discontinuity or steep gradient) detector---the work in \cite{MRAdetect2014} presents one of the many discontinuity detectors. Our recommendation is to make the above choice for $\mcal M$ and $\mcal R$ in case the shock detector returns an empty set. 
\end{remark}
\subsection{Computing the POD basis}\label{sec: POD}
Using the snapshots $\{g(\cdot,\ztrain{i})\}_i$, we compute the POD basis 
\begin{gather}\mcal X_n = \{\basisPOD{i}\}_{i=1,\dots,n}.
\end{gather}
The computation relies on the singular value decomposition (SVD) of a snapshot matrix $\snapMat{G}\in\mbb R^{N\times \sampleTrain}$ given as
\begin{gather}
\snapMat{G} := \left(G(\ztrain{1}),\dots,G(\ztrain{\sampleTrain})\right), \label{def S}
\end{gather}
where the vector $G(z)\in\mbb R^N$ contains the degrees-of-freedom (Dofs) of the transformed solution $g(\cdot,z)$. SVD of the snapshot matrix provides 
$$
\snapMat{G} = \mcal U\Sigma \mcal V^T,
$$
where $\mcal U\in\mbb R^{N\times N}$ and $\mcal V\in\mbb R^{\sampleTrain\times\sampleTrain}$ are orthogonal matrices containing the left and the right singular vectors of $\snapMat{G}$, respectively. Furthermore, $\Sigma\in\mbb R^{N\times\sampleTrain}$ contains, at its diagonals, the singular values of $\snapMat{G}$ that we assume are arranged in a descending order.

Let the vector $V_i\in\mbb R^N$ contain the Dofs of $v_i\in \mcal X_N$. We set $V_i$ to be the $i$-th column of the matrix $\mcal U$. Equivalently,
$$
\left(V_i\right)_j = \mcal U_{ji},\hspB\forall i\in \{1,\dots,n\},j\in\{1,\dots,N\}.
$$
Using the Schmidt-Eckart-Young theorem (see \cite{OptimalSVD}), we can quantify the best-approximation error of approximating the snapshots in the POD basis.
Collecting all the POD modes in a matrix
\begin{gather}
X_n := \left(V_1,\dots,V_n\right), \label{POD mat}
\end{gather}
a bound for the projection error reads
\begin{gather}
E_{n}^{proj}(\snapMat{G}):= \frac{1}{\|\snapMat{G}\|_F}\|\snapMat{G}-X_n X_n^T\snapMat{G}\|_F = \frac{1}{\|\snapMat{G}\|_F}\sqrt{\sum_{i=n+1}^{\min(N,\sampleTrain)}\sigma_i^2}, \label{def Eproj}
\end{gather}
where $\|\cdot\|_F$ represents the Frobenius norm of a matrix, and $\sigma_i$ is the $i$-th singular value of $\snapMat{G}$. Later (in \Cref{sec: num results}), using numerical experiments, we will study how the POD projection error $E_{n}^{proj}(\snapMat{G})$ decays with $n$.

\subsection{GPR for the POD coefficients}\label{sec: GPR for alpha} Let $\coeffpod{g}{z}\in\mbb R^n$ denote a vector containing the POD coefficients of $g(\cdot,z)$. The crux of our technique is that for any given $z\in \parSp$, we approximate the vector $\coeffpod{g}{z}$ using GPR via $\alpha(z)\approx \alpha^{\mcal R}(z)$ and recover the approximation 
\begin{gather}
g(\cdot,z)\approx g_n(\cdot,z) :=\sum_{i=1}^n\alpha^{\mcal R}_{g,i}(z) v_i\in \opn{span}(\mcal X_n).\label{POD transf sol}
\end{gather}
The discussion below outlines a technique to compute $\coeffpodR{g}{z}$.

 As explained in \Cref{sec: GPR}, GPR relies on an offline training and an online prediction step. The former computes the hyper-parameters of a GPR, whereas the latter, for any $z\in\parSp$, predicts a value for $\coeffpod{g}{z}$. Following are the details of these two steps.
\begin{enumerate}
\item \textbf{Training step:} We orthogonally project each of the snapshots in $\{g(\cdot,\ztrain{i})\}_i$ onto the corresponding POD basis vectors $\{v_i\}_i$ and compute the POD coefficients given as 
\begin{gather}
\alpha_{g,i}(\ztrain{j}) = \frac{\lan v_i,g(\cdot,\ztrain{j})\ran_{L^2(\Omega)}}{\|v_i\|^2_{L^2(\Omega)}},\hspB \forall i\in\{1,\dots,n\},\hsp j\in\{1,\dots,\sampleTrain\}.
\end{gather}
Note that, by definition, the basis vector $\{v_i\}_i$ are $L^2$-orthogonal. The above computation provides the training data $\{(\ztrain{j},\coeffpod{g}{\ztrain{j}}\}_j,$ which we use to train a GPR. For each component of $\coeffpod{g}{z}$, we train a separate uncorrelated GPR---\Cref{remark: vector value GPR} below elaborates on this further. As \Cref{sec: GPR train} explains, training a GPR for the $i$-th component of $\alpha_g(z)$ entails computing the hyper-parameter $\theta_{\alpha,i}\in\mbb R^{2\times(p+1)}$ with $p$ being the dimension of $\parSp$.
\item \textbf{Prediction step:} With the hyper-parameters $\{\theta_{\alpha,i}\}_i$ at hand, we consider the mean-value prediction for $\coeffpod{g}{z}$ that reads
\begin{gather}
\alpha_{g,i}(z)\approx\alpha^{\mcal R}_{g,i}(z):= \mcal R_0(\theta_{\alpha,i},\mcal D_z,z),\hspB\forall i\in \{1,\dots,n\}.\label{def alphaR}
\end{gather}
The above approximation finally provides the POD approximation $g_n(\cdot,z)$ given in \eqref{POD transf sol}. Recall that $\mcal D_z$ is a set containing the parameter samples and is defined in \eqref{def D}. Furthermore, $\mcal R$ denotes the GPR model given in \eqref{def R}.
\end{enumerate}

\begin{remark}[GPR for vector-valued functions]\label{remark: vector value GPR}
While training the GPR, we assume that $\coeffpod{g}{z}$ has uncorrelated components \cite{Guo2019,MovingBC2021}. This certainly introduces some inaccuracies because for most pPDEs, $\coeffpod{g}{z}$ has coupled correlated components. As of yet, particularly for non-linear pPDEs, it is unclear how one can account for this coupling using the cross-correlation technique proposed in \cite{GPRVectorValued}.
\end{remark}

\begin{remark}[Extrapolation capabilities]\label{remark: extrapolation}
Regression is known to be inaccurate with extrapolation \cite{Guo2019,GPR,JanNN2018}. Our MOR technique is no exception to this limitation. Heuristics suggest (see \cite{MovingBC2021}) that the variance of the GP can be used to quantify the extrapolation capabilities, but the success of such a technique---particularly for hyperbolic pPDEs considered in our numerical experiments---is unclear. We include the corners of the parameter domain $\parSp$ in the parameter samples and avoid extrapolation altogether. 

We emphasize that our non-intrusive technique shares this limitation with the intrusive techniques proposed in \cite{Nair,RegisterMOR}. In entirety, these intrusive techniques are only partially intrusive---regression is used to compute the spatial transform $\varphi$, which introduces inaccuracies during extrapolation. 
\end{remark}

\section{Approximation of the untransformed solution}\label{sec: approx u}
We want to recover an approximation for $\solFV{z}$ from the POD approximation of the transformed snapshot $g_n(\cdot,z)$ given in \eqref{POD transf sol}. We emphasize that this recovery needs to be performed online. At first glance, the approximation $\solFV{z}\approx \Pi g_n(\varphi(\cdot,\zRef,z)^{-1},z)$ seems reasonable, with $\varphi(\cdot,\zRef,z)^{-1}$ computed using the non-linear least-squares problem
\begin{gather}
\varphi(x,\zRef,z)^{-1} := \argmin_{y\in\Omega}\|\varphi(y,\zRef,z)-x\|^2_{l^2}.\label{def inv phi}
\end{gather}
However, a solution to the above problem comes at a high cost \cite{Noblet2008}. We need  $\varphi(\cdot,\zRef,z)^{-1}$ inside all the $N$ spatial grid cells i.e., we need to solve the above problem at least $\mcal O(N)$ times, with each solution requiring a few iterations. This procedure can easily dominate the cost of our MOR technique and can make it more expensive than a high-fidelity solver. 

For the above reason, we refrain from solving the least-squares problem online and instead approximate its solution in the POD basis. We follow the same line of procedures as that used to approximate the transformed snapshot $g(\cdot,z)$ in \Cref{sec: approx transf sol}. In the offline phase, we solve the above problem and collect snapshots of $\varphi^{-1}$, compute the POD basis, collect training data for the POD coefficients and train a GPR. In the online phase, we query the GPR for the POD coefficients and recover an approximation for $\varphi^{-1}$. The details are as follows:

\subsection{Snapshots of $\varphi^{-1}$} We solve the above least-squares problem and collect the snapshots $\{\spTransf{\zRef}{\ztrain{i}}^{-1}\}_{i=1,\dots,\sampleTrain}$; recall that $\{\ztrain{i}\}_i$ are the parameter samples defined in \eqref{def D}. Similar to the solution $\sol{z}$, we compute these snapshots in a finite-dimensional high-fidelity approximation space. 

As \Cref{remark: diffeo phi} states, we expect $\spTransf{\zRef}{z}$ to be a diffeomorphism. Therefore, it is reasonable to approximate $\spTransf{\zRef}{z}^{-1}$ in a finite-element space of continuous functions. We choose this space to be the span of continuous piecewise linear functions (or the so-called hat-functions) defined over the triangulation of $\Omega$ given in \eqref{space grid}.
We denote this finite-element space by $\td{\mcal X}_{\td N}$, which is such that $\opn{dim}(\td{\mcal X}_{\td N}) = \td N$. Furthermore, the value of $\td N$ equals the number of vertices in the spatial mesh. 

 We label the high-fidelity approximation of $\spTransf{\zRef}{z}^{-1}$ by
\begin{gather}
\spTransf{\zRef}{z}^{-1} \approx \td\varphi_{N}(\cdot,\zRef,z) \in \td{\mcal X}_{\td N},
\end{gather}
and compute it by projecting $\spTransf{\zRef}{z}^{-1}$ onto $\td{\mcal X}_{\td N}$ i.e., we set
$$
\td\varphi_N(\cdot,\zRef,z) := \td{\Pi}\left(\spTransf{\zRef}{z}^{-1}\right),
$$
where $\td \Pi:L^2(\Omega)\to\td{\mcal X}_{\td N}$ is the orthogonal projection operator. We consider a quadrature routine to compute the projection.  The value of $\varphi(\cdot,\zRef,z)^{-1}$ at the quadrature points results from solving the least-squares problem given in \eqref{def inv phi}; to this end, we use the \texttt{lsqnonlin} routine from MATLAB. For the experiments reported later, we use a tensorized set of $3\times 3$ Gauss-Legendre quadrature points in each spatial cell.

\subsection{Computing the POD basis}
Using the snapshots $\{\td\varphi_N(\cdot,\zRef,\ztrain{i})\}_{i=1,\dots,\sampleTrain}$ computed above, for any $z\in\mcal Z$, we seek an efficient approximation of $\td\varphi_N(\cdot,\zRef,z)$. As discussed in \Cref{remark: diffeo phi}, we prefer to split $\td\varphi_N(\cdot,\zRef,z)$ as 
\begin{gather}
\td\varphi_N(\cdot,\zRef,z) = \opn{Id} + \td{\Psi}_N(\cdot,\zRef,z),
\end{gather}
and approximate the displacement field $\td{\Psi}_N(\cdot,\zRef,z)$. In contrast to the solution $\sol{z}$, which has a jump-discontinuity along $\parSp$, the displacement field is sufficiently regular along $\parSp$---\Cref{remark: approx phi POD} below provides further elaboration. Therefore, it is reasonable to expect that it is well-approximable in a sufficiently low-dimensional linear reduced space; numerical experiments will further justify our expectations. We construct this reduced space using POD. 

Before discussing further details, let us recall that $\td\Psi_N(x,\zRef,z)$ is a two-dimensional vector i.e., 
$$
\td\Psi_N(\cdot,\zRef,z) = \left(\td\Psi^{(1)}_N(\cdot,\zRef,z),\td\Psi^{(2)}_N(\cdot,\zRef,z)\right)^T.
$$
In the following, for brevity, we present a technique to approximate $\td\Psi^{(1)}_N(\cdot,\zRef,z)$. The same procedure applies to $\td\Psi^{(2)}_N(\cdot,\zRef,z)$.
Similar to the snapshot matrix $\snapMat{G}$ defined earlier in \eqref{def S}, using snapshots of $\td\Psi^{(1)}_N(\cdot,\zRef,z)$, we define a snapshot matrix $\snapMat{\td\Psi^{(1)}}$. Then, applying the SVD-based technique outlined earlier in \Cref{sec: POD}, we compute the POD basis for the displacement field given as 
$$\mcal X_{\nPsi}:=\{v_i^{\td\Psi^{(1)}}\}_{i=1,\dots,\nPsi},$$ 
where the parameter $\nPsi\in\mbb N$ denotes the dimensionality of $\mcal X_{\nPsi}$. Note that we use the same $\nPsi$ to approximate both the components of $\td\Psi_N$.
\subsection{GPR for the POD coefficients} Let $\alpha_{\td\Psi^{(1)}}(z)\in\mbb R^\nPsi$ denote a vector containing the POD coefficients of the displacement field $\td{\Psi}^{(1)}_N(\cdot,\zRef,z)$. For any given $z\in\parSp$, we approximate $\alpha_{\td\Psi^{(1)}}(z)$ using GPR. The procedure remains exactly the same as that used to approximate the POD coefficients of $g(\cdot,z)$ in \Cref{sec: GPR for alpha}. Following is a brief recall.

In the offline phase, by projecting each of the snapshots in $\{\td\Psi^{(1)}_N(\cdot,\zRef,\ztrain{i})\}_{i}$ onto the corresponding POD basis vector $\{v_i^{\td\Psi^{(1)}}\}_i$, we collect training data from the mapping $z\mapsto\alpha_{\td\Psi^{(1)}}(z)$. With this training data, for the $i$-th component of $\coeffpod{\td\Psi^{(1)}}{z}$, we compute the hyper-parameter $\theta_{\td\Psi^{(1)},i}$ of the GPR. In the online phase, for any $z\in\parSp$, we approximate the POD coefficient $\alpha_{\td\Psi^{(1)}}(z)$ using the mean-value estimate given as 
\begin{gather}
\alpha_{\td\Psi^{(1)},i}(z)\approx \alpha^{\mcal R}_{\td\Psi^{(1)},i}(z):=\mcal R_0(\theta_{\td\Psi^{(1)},i},\mcal D_z,z),\hspB\forall i \in \{1,\dots,\nPsi\}, \label{def alphaPsiR}
\end{gather}
and recover an approximation for the first-component of the displacement field
$$
\td\Psi^{(1)}_N(\cdot,\zRef,z)\approx \td\Psi^{(1)}_{\nPsi}(\cdot,\zRef,z) := \sum_{i=1}^{\nPsi}\alpha^{\mcal R}_{\td\Psi^{(1)},i}v^{\td\Psi^{(1)}}_i.
$$
Similarly, we can approximate the second-component $\td\Psi^{(2)}_N(\cdot,\zRef,z)$ and recover
\begin{gather}
\td\Psi_N(\cdot,\zRef,z)\approx \td\Psi_{\nPsi}(\cdot,\zRef,z) := \left(\td\Psi^{(1)}_{\nPsi}(\cdot,\zRef,z),\td\Psi^{(2)}_{\nPsi}(\cdot,\zRef,z)\right)^T.
\end{gather}
 This finally provides the following approximation for the spatial transform
\begin{gather}
\td\varphi_N(\cdot,\zRef,z)\approx \td\varphi_\nPsi(\cdot,\zRef,z) := \opn{Id} + \td\Psi_\nPsi(\cdot,\zRef,z), \label{approx phiInv}
\end{gather}
that we use to approximate $\solFV{z}$ via 
\begin{gather}
\solFV{z}\approx u_{n,\nPsi}(\cdot,z):= \Pi g_n(\td\varphi_\nPsi(\cdot,\zRef,z),z). \label{approx u}
\end{gather}
Recall that $\Pi$ is an orthogonal projection operator from $L^2(\Omega)$ onto the high-fidelity space $\mcal X_N$. Furthermore, $g_n(\cdot,z)$ is the reduced approximation for the transformed solution $g(\cdot,z)$ and is given in \eqref{def g}.

\begin{remark}[Both $n$ and $\nPsi$ determine the approximation quality]
Note that two parameters---$n$ and $\nPsi$---control the accuracy of our approximation. Increasing both $n$ and $\nPsi$ increases the number of POD modes used to approximate the transformed solution $g(\cdot,z)$ and the displacement field $\td\Psi_N(\cdot,\zRef,z)$, respectively. As one might expect, keeping one parameter fixed and increasing the other offers diminishing returns in terms of the approximation accuracy. Numerical experiments will further study this behaviour.
\end{remark}

\begin{remark}[Approximability of $\td{\varphi}_N$ in the POD basis]\label{remark: approx phi POD}
We assume that $\td{\varphi}_N(\cdot,\zRef,z)$ is well-approximable in a POD space, which is linear by construction. The observation that for all $x\in\Omega$, $\td{\varphi}_N(x,\zRef,\cdot)$ is sufficiently regular along $\parSp$ motivates this assumption. As stated earlier---and further explained in \cite{Welper2017,SarnaCalib2020}---$\td{\varphi}_N(\cdot,\zRef,z)$ matches discontinuities between the reference solution $u_N(\cdot,\zRef)$ and  some other solution $u_N(\cdot,z)$. For several problems involving parameter-invariant discontinuity-topology (see \Cref{sec: approx transf sol}), the locations of these spatial discontinuities vary smoothly in $\parSp$, resulting in a $\td{\varphi}_N$ that is smooth along $\parSp$. Our numerical experiments will provide further elaboration.
\end{remark}

\subsection{Error surrogate}
To certify the quality of our approximation, following the works in \cite{ROMES,Guo2019}, we develop a surrogate for the relative $L^1$-error $E(z,n,\nPsi)$ defined as
\begin{gather}
E(z,n,\nPsi) := \frac{1}{\|\solFV{z}\|_{L^1(\Omega)}}\|\solFV{z} - u_{n,\nPsi}(\cdot,z)\|_{L^1(\Omega)},\label{error L2}
\end{gather}
where $\solFV{z}$ and $u_{n,\nPsi}(\cdot,z)$ is the high-fidelity and the reduced approximation, respectively. We represent the surrogate via
$$
E(z,n,\nPsi)\approx E^{\mcal R}(z,n,\nPsi),
$$
and compute it using a standard GPR-based technique. Offline, at the training points $\{\ztrain{i}\}_{i=1,\dots,\sampleTrain}$, we collect the training data $\{(\ztrain{i},E(\ztrain{i},n,\nPsi))\}_i$. Using this training data, we train a GPR i.e., we compute the hyper-parameters $\theta_E\in\mbb R^{2\times (p+1)}$. Online, for any $z\in\parSp$, we compute the surrogate $E^{\mcal R}(z,n,\nPsi)$ via
\begin{gather}
E^{\mcal R}(z,n,\nPsi):= \mcal R_2(\theta_E,\mcal D_z,z), \label{def ER}
\end{gather}
where $\mcal R_2$ is the GPR model given in \eqref{def R}.
We quantify the accuracy of our error surrogate by an efficiency index $\eta(z,n,\nPsi)$ defined as 
\begin{gather}
\eta(z,n,\nPsi):= \frac{E^{\mcal R}(z,n,\nPsi)}{E(z,n,\nPsi)}.  \label{eff index}
\end{gather}
Later, we study $\eta(z,n,\nPsi)$ using numerical experiments.

 Note that $\mcal R_2(\theta_E,\mcal D_z,z)$ is two standard deviations away from the mean-value and therefore, is pessimistic as compared to a mean-value estimate. Since we use the same sample parameters to compute the POD basis and then train the GPR, we cannot expect much accuracy from a mean-value prediction of $E(z,n,\nPsi)$. Therefore, we instead consider a pessimistic estimate.

\section{Summary}\label{sec: summary}
\Cref{summary: offline} and \Cref{summary: online} summarize the offline and the online stages of the algorithm, respectively.
\begin{algorithm}[ht!]
\caption{Summary: Offline phase}
\begin{algorithmic}[1] \label{summary: offline}
\STATE \textbf{Input: } $n$, $\nPsi$, $\sampleTrain$
\STATE \textbf{Output: }$\theta_\alpha$, $\theta_{\td \Psi}$, $\theta_E$, $\mcal D_z$, $\{v_i\}_{i=1,\dots,n}$, $\{v^{\td \Psi}_i\}_{i=1,\dots,\nPsi}$
\STATE Collect the parameter samples $\mcal D_z:=\{\ztrain{i}\}_{i=1,\dots,\sampleTrain}$ from $\parSp$.
\STATE Compute the transformed snapshots $\{g(\cdot,\ztrain{i})\}_{i=1,\dots,\sampleTrain}$.
\STATE Compute the POD basis $\{v_i\}_{i=1,\dots,n}$.
\STATE Compute the snapshots $\{\td\Psi_N(\cdot,\zRef,\ztrain{i})\}_{i=1,\dots,\sampleTrain}$.
\STATE Compute the POD basis $\{v^{\td\Psi}_i\}_{i=1,\dots,\nPsi}$.
\STATE At the parameter samples in $\mcal D_z$, collect the training data from the mappings $z\mapsto \coeffpod{g}{z}$, $z\mapsto\alpha_{\td\Psi}(z)$ and $z\mapsto E(z,n,\nPsi)$ and use it to compute the hyper-parameters $\theta_\alpha$, $\theta_{\td\Psi}$ and $\theta_E$, respectively, for the GPR. 
\end{algorithmic}
\end{algorithm}

\begin{algorithm}[ht!]
\caption{Summary: Online phase}
\begin{algorithmic}[1] \label{summary: online}
\STATE \textbf{Input: }$\theta_{\alpha}$, $\theta_{\td \Psi}$, $\theta_E$, $\mcal D_z$, $\{v_i\}_{i=1,\dots,n}$, $\{v^{\td \Psi}_i\}_{i=1,\dots,\nPsi}$, $z\in\parSp$\\
\STATE \textbf{Output: }$u_{n,\nPsi}(\cdot,z)$, $E^{\mcal R}(z,n,\nPsi)$
\STATE Compute $\alpha_g^{\mcal R}(z)$ using \eqref{def alphaR} and recover $g_n = \sum_j\alpha_{g,j}^{\mcal R}(z) v_j$.
\STATE Compute $\alpha^{\mcal R}_{\td\Psi}(z)$ using \eqref{def alphaPsiR} and recover $\td\Psi_{\nPsi}(\cdot,\zRef,z) := \sum_{i=1}^{\nPsi}\alpha^{\mcal R}_{\td\Psi,i}v^{\td\Psi}_i$.
\STATE Compute $u_{n,\nPsi} = \Pi g_n(\opn{Id} + \td\Psi_{\nPsi}(\cdot,\zRef,z),z)$.
\STATE Compute $E^{\mcal R}(z,n,\nPsi)$ using \eqref{def ER}. 
\end{algorithmic}
\end{algorithm}

\subsection{Computational costs: online phase}
We compute the costs of the different steps outlined in \Cref{summary: online}. 
\begin{enumerate}
\item \texttt{line-3:} Querying the GPR model for the $n$ components of $\coeffpodR{g}{z}$ requires $\mcal O(\sampleTrain n)$ operations \cite{GPR}. Computing $g_n(\cdot,z)$ requires $\mcal O(nN)$ operations.
\item \texttt{line-4:} Querying the GPR model for the $\nPsi$ components of $\alpha^{\mcal R}_{\td\Psi}(z)$ requires $\mcal O(\sampleTrain \nPsi)$ operations. Computing $\td\Psi_{\nPsi}(\cdot,\zRef,z)$ requires $\mcal O(\nPsi N)$ operations.
\item \texttt{line-5:} Computing $u_{n,\nPsi}(\cdot,z)$ on a structured grid requires $\mcal O(N)$ operations.
\item \texttt{line-6:} Computing the scalar $E^{\mcal R}(z,n,\nPsi)$ requires $\mcal O(\sampleTrain)$ operations.
\end{enumerate}
The total cost sums up to
\begin{equation}
\begin{aligned}
\mcal C^{MOR}_{tot} =&\underbrace{\mcal O(\sampleTrain n) + \mcal O(n N)}_{\text{POD approximation of $g(\cdot,z)$}}+\underbrace{\mcal O(\sampleTrain \nPsi) + \mcal O(\nPsi N)}_{\text{POD approximation for $\td\Psi_N(\cdot,\zRef,z)$}}\\
& +\underbrace{ \mcal O(N)}_{\text{Computation of $u_{n,\nPsi}(\cdot,z)$}} + \underbrace{\mcal O(\sampleTrain)}_{\text{Computation of $E^{\mcal R}(z,n,\nPsi)$}}.  \label{cost MOR}
\end{aligned}
\end{equation}

We study how $\mcal C^{MOR}_{tot}$ scales with the parameter-domain dimension $p$. Assume that to collect the parameter samples from $\parSp$, we take a tensor-product of $m\in\mbb N$ uniformly placed parameter samples along each parameter dimension. Then, the number of parameter samples scales as 
$\sampleTrain = \mcal O(m^{p}),$ which provides
\begin{equation}
\begin{aligned}
\mcal C^{MOR}_{tot} =&\mcal O(m^p n) +\mcal O(m^p \nPsi) +  \mcal O(m^p)\\
& +  \mcal O(n N) + \mcal O(\nPsi N)+\mcal O(N). 
\end{aligned}
\end{equation}
Observe that for extremely large values of $p$, the cost of querying a GPR---which introduces the $m^p$ dependence above---might outweigh all the other costs. This might result in a MOR technique that is more expensive than a high-fidelity solver. Thus the practical applicability of our method is limited to moderate parameter domain dimensions. We share this limitation with the intrusive technique developed in \cite{RegisterMOR}. The reason being that this technique approximates $\varphi$ using RBFs (see \Cref{remark: extrapolation}) which, similar to a GPR, are expensive to query for extremely large values of $p$.

\section{Numerical Results}\label{sec: num results}
We abbreviate our MOR technique as GPR-TS-MOR. The abbreviation TS stands for transformed snapshots. The numerical experiments compare it to the S-PROJ technique. In S-PROJ, we collect solution snapshots (S)---do not perform any snapshot transformation---compute the POD basis and approximate the solution in the POD basis. We orthogonally project the solution onto the POD basis. The high-fidelity solver is abbreviated as HF.

\subsection{Description of the test cases} We consider the following test cases involving hyperbolic (test case-1 and 2) and parabolic (test case-3) equations.  
\begin{enumerate}
\item \textbf{Test-1 (1D wave equation):} We consider the 1D (in space) wave equation (rewritten as a first order system)
    \begin{gather}
        \pd_t u(x,t,\mu) + A \pd_x u(x,t,\mu) = 0,\hspB\forall (x,t,\mu)\in\Omega\times [0,T]\times\mcal P, \label{wave equation}
    \end{gather}
    where $u = (u_1,u_2)^T$ is the vector-valued solution, and the matrix $A$ reads
    \begin{gather}
    A := \left(\begin{array}{c c}
    0 & 1 \\ 
    1 & 0
    \end{array}\right). 
    \end{gather}
    We choose $\Omega = (-0.3,3)$, and $T = 0.8$. As the initial data, for all $x\in\Omega$, we consider the linear superposition
    \begin{equation}
    \begin{gathered}
   u_1(x,t=0,\mu) =  \frac{1}{\sqrt{2}}\left(w_1(x,\mu) + w_2(x,\mu)\right),\\ u_2(x,t=0,\mu) = \frac{1}{\sqrt{2}}\left(-w_1(x,\mu) + w_2(x,\mu)\right),
    \end{gathered}
    \end{equation}
    where $w_1(\cdot,\mu)$ and $w_2(\cdot,\mu)$ are two sin-function bumps given as
    \begin{equation}
    \begin{aligned}
    w_1(x,\mu) := &\mu\times(\sin(2\pi (x + 0.2)) + 1)\mathbbm{1}_{[\delta_1-0.5,\delta_1]}(x),\\
     w_2(x,\mu) := &\mu\times(\sin(2\pi (x-2.3)) + 1)\mathbbm{1}_{[\delta_2-0.5,\delta_2]}(x). \label{def: sin bumps}
    \end{aligned}
    \end{equation}
    We set $\delta_1:=0.3$, and $\delta_2 := 2.8$. Furthermore, we set $\mcal P:=[0.5,2]$. Thus, the parameter domain is two-dimensional and reads $\mcal Z = [0,T]\times \mcal P$. Along the boundary $\pd\Omega\times \parSp$, we prescribe $u = 0$. 
    \item \textbf{Test-2 (2D Burgers' equation):} We consider the 2D (in space) Burgers' equation given as 
    \begin{gather}
    \pd_t u(x,t) + \left(\frac{1}{2},\frac{1}{2}\right)^T\cdot\nabla u(x,t)^2=0,\hsp \forall (x,t)\in\Omega\times [0,T].
    \end{gather}
    The initial data is a characteristic function over a square and reads
    \begin{gather}
    u(x,t = 0) = \begin{cases}
    1,\hspB &x\in [0,0.5]^2\\
    0,\hspB &\text{else}
\end{cases}. \label{ic burgers}
    \end{gather}
    We set $\Omega = (-0.1,1.5)^2$ and $T=2$. Along the boundary $\pd\Omega\times \parSp$, we prescribe $u = 0$. Time is the sole parameter for this problem, i.e., $\parSp = [0,T]$. 
    \item \textbf{Test-3 (2D Heat conduction):} We consider the 2D (in space) heat conduction problem from \cite{RegisterMOR} given as
    \begin{gather}
    -\nabla\cdot(\beta(x,z)\nabla u(x,z)) = 1,\hsp\forall (x,z)\in\Omega\times\parSp,
    \end{gather}
    where $\Omega = [0,1]^2$, and $\parSp = [-0.05,0.05]$. The heat conductivity coefficient $\beta(x,z)\in\mbb R$ is discontinuous in $\Omega \times \parSp$ and reads 
    \begin{gather}
    \beta(x,z) := 0.1 + 0.9\times \mathbbm{1}_{\Omega_z},\hspB \Omega_{z}:=\left\{x\in\Omega\hsp :\hsp \|x-\bar x_z\|_{\infty}\leq \frac{1}{4}\right\}.
    \end{gather}
 The center $\bar x_z $ equals $[1/2 + z,1/2 + z]$.  Along the boundary $\pd\Omega\times\parSp$, we prescribe $u = 0$.
\end{enumerate}

\begin{remark}[Software and hardware details]
All the simulations are run using MATLAB, in serial, and on a computer with two Intel Xeon Silver 4110 processors, 16 cores each and $92$GB of RAM. 
\end{remark}

\begin{remark}[HF solver]
For test case-1 and 2, our HF solver is a second-order finite-volume scheme with a Van-Leer flux-limiter combined with a second-order explicit Runge-Kutta time-stepping scheme. We use the local-Lax-Friedrich numerical flux, and set the CFL number to $0.5$. For test case-3, our HF solver is a continuous Galerkin $\mbb P^1$ finite-element solver. 
\end{remark}
\subsection{Choice of the matching criterion $\mcal M$}
Following the empirical success reported in \cite{Nair,Welper2017,RegisterMOR,Taddei2020ST}, for test cases 1 and 2, which involve hyperbolic equations, we choose $\mcal M$ (appearing in \eqref{L2 opt}) as the $L^2$-distance between the transformed and the reference snapshot i.e.,
\begin{gather}
\mcal M(\Psi^*,z,\zRef)= \|u_N(\opn{Id} + \Psi^*,z)-\solFV{\zRef}\|^2_{L^2(\Omega)}.
\end{gather}
For test case-3, we assume that one has access to both the solution snapshots and a parametrized description of the boundary $\pd\Omega_z$. Usually this description is already available during the spatial mesh generation process. The jump in the heat conductivity coefficient along $\pd\Omega_z$ can be viewed as a change in the material properties, which is known a-priori, and might be used, for instance, to refine the mesh along the boundary $\pd\Omega_z$. With an access to the boundary $\pd\Omega_z$, we set the matching criterion to \cite{RegisterMOR}
\begin{gather}
\mcal M(\Psi^*,z,\zRef) = \sum_{i=1}^{N_p}\|\hat x^{(i)}_{z}-(\hat x^{(i)}_{\zRef} + \Psi^*(\hat x^{(i)}_{\zRef}))\|^2_{2},
\end{gather}
where $\{\hat x_z^{(i)}\}_{i=1,\dots,N_p}$ represent a set of uniformly placed points placed along $\pd\Omega_{z}$. We set $N_p$ to $400$. Observe that the minimization of the above matching criterion pushes the boundary $\varphi(\pd\Omega_{\zRef},\zRef,z)$ to be aligned with that of $\pd\Omega_z$. Consequently, the transformed solution $g_N(\cdot,z)$ does not exhibit steep gradients along $\parSp$.

Since the choice of $\mcal M$ is application dependent, our technique does not run in an entire black-box fashion and requires at least some information about the underlying physical process. Nevertheless, this information is often available because a user is usually aware of the physical process being simulated---even though one might not have access to the underlying discrete evolution operators.

\subsection{Average error}
We quantify the error over the entire parameter domain via the average relative $L^1$-error $E_{a}(n,\nPsi)$ defined as 
\begin{gather}
E_{a}(n,\nPsi):= \frac{1}{\#\mcal D^{tst}_z}\sum_{z\in \mcal D^{tst}_z} E(z,n,\nPsi). \label{error avg}
\end{gather}
The error $E(z,n,\nPsi)$ is as defined in \eqref{error L2}.
The test samples $\mcal D^{tst}_z$ consist of $200$ uniformly and idependently sampled parameters from $\parSp$. Furthermore, $\#\mcal D^{tst}_z$ represents the size of $\mcal D^{tst}_z$.
 Replacing $E(z,n,\nPsi)$ by the surrogate $E^{\mcal R}(z,n,\nPsi)$ defined in \eqref{def ER}, we recover the approximation $$E_{a}(n,\nPsi)\approx E^{\mcal R}_{a}(n,\nPsi).$$ 

\subsection{Test-1} We discretize $\Omega$ with $N=10^3$ grid-cells, resulting in a grid size of $\Delta x = 3.3\times 10^{-3}$. The training data $\mcal D_z$ is a set of uniformly placed $40\times 20$ points inside $\parSp$. We approximate the displacement field $\dev{\zRef}{z}$ in the polynomial space $\mcal P_{M=4}$, where the value of $M$ results from the procedure outlined in \Cref{sec: varphi}.  Recall that the solution here is vector-valued with $\sol{z} = (u_1(\cdot,z),u_2(\cdot,z))^T$. As stated earlier, we apply the technique developed in the previous sections to each of the components of $\sol{z}$. For brevity, in the following, we present the results for $u_1$. For $u_2$, the results are similar.

\subsubsection{Study of the POD projection error} Consider the POD projection error $E^{proj}_n(\mcal S)$ defined in \eqref{def Eproj}. We compare this error for two different snapshot matrices, $\mcal S = \snapMat{G}$ and $\mcal S = \snapMat{U}.$ The former is defined in \eqref{def S} and contains transformed snapshots. The latter contains the snapshots of the (untransformed) solution and reads
\begin{gather}
\snapMat{U}:=\left(U(\ztrain{1}),\dots,U(\ztrain{\sampleTrain})\right), \label{def SU}
\end{gather}
where $U(z)\in\mbb R^N$ contains the Dofs of $u_N^{(1)}(\cdot,z)$ with $u_N^{(1)}(\cdot,z)$ being a HF approximation to $u_1(\cdot,z)$. Note that GPR-TS-MOR and S-PROJ use the POD modes of $\snapMat{G}$ and $\snapMat{U}$, respectively, to approximate the transformed and the untransformed solution, respectively.

\begin{figure}[ht!]
\centering
\subfloat[]{\includegraphics[width=2.5in]{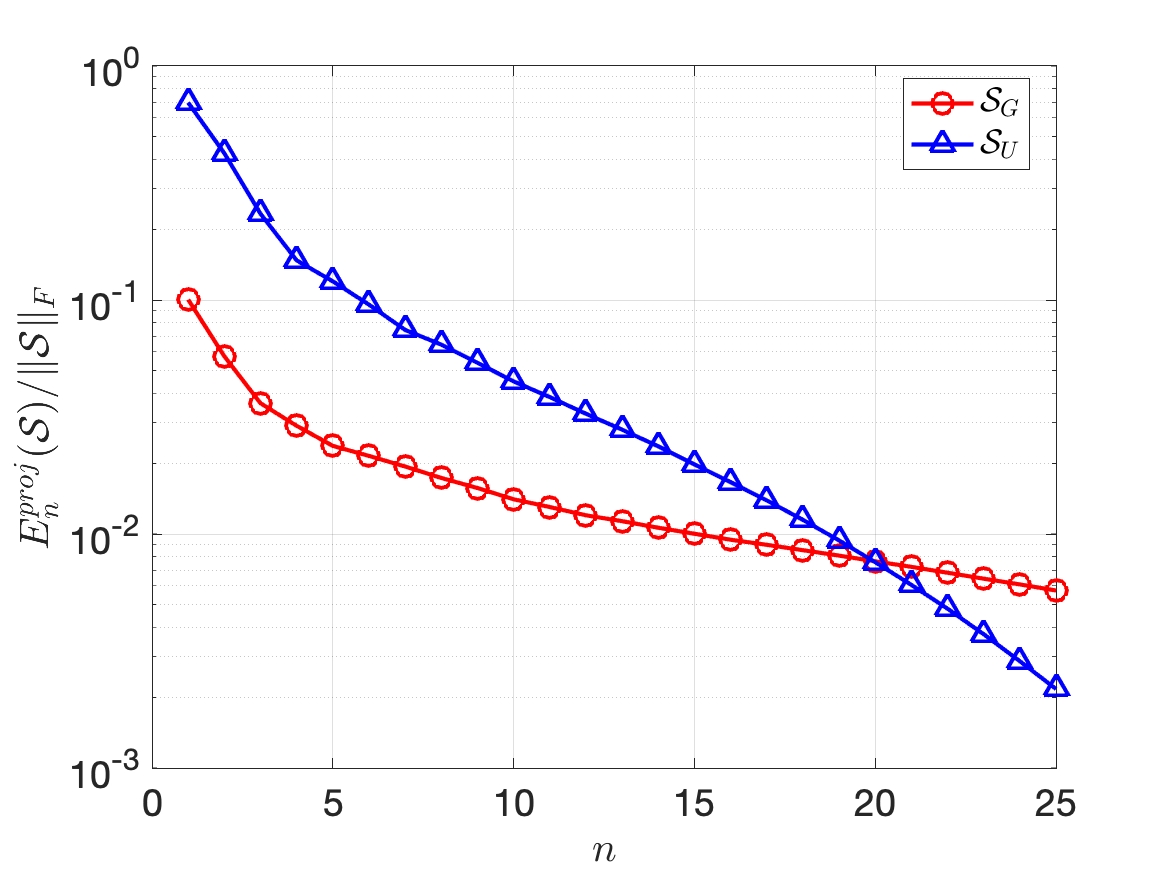}\label{fig: test-1 sigma}} 
\hfill
\subfloat[]{\includegraphics[width=2.5in]{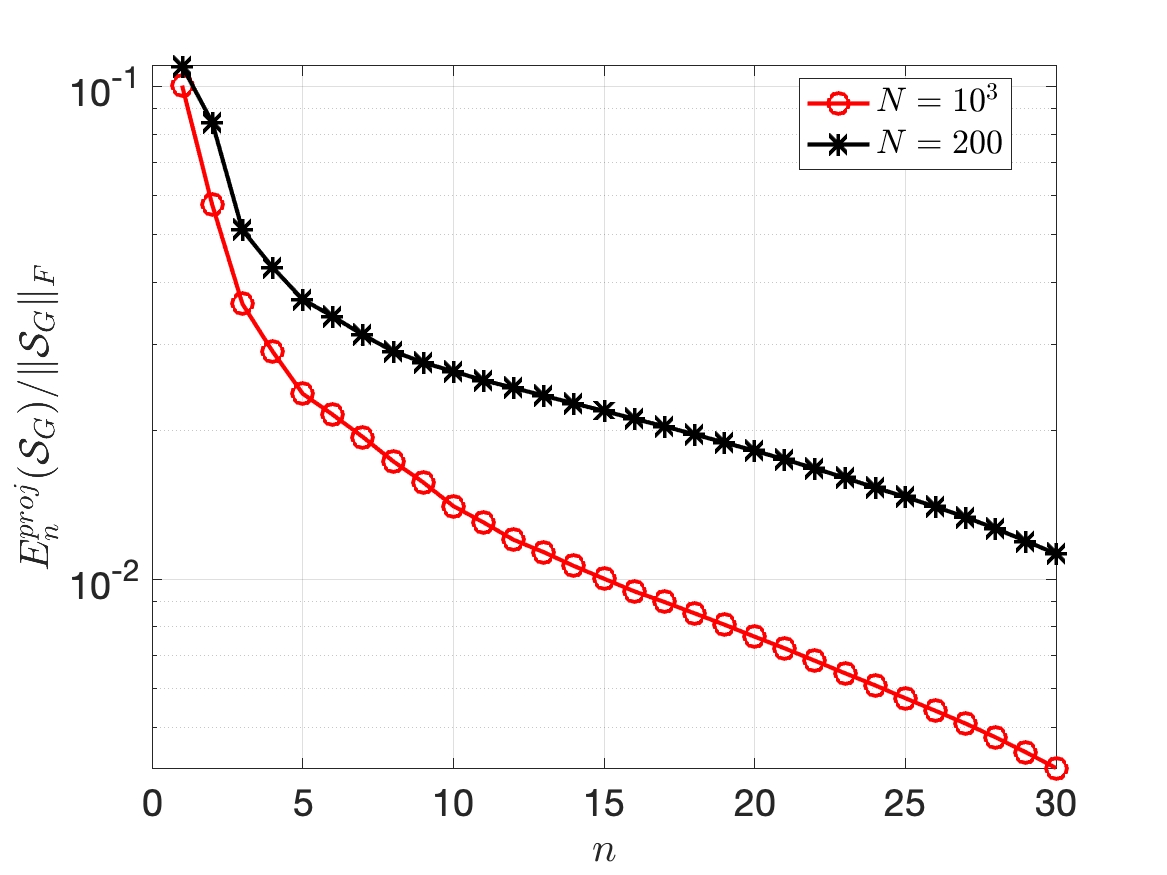}\label{fig: test-1 sigma diffN}}
\caption{\textit{Results for test-1. (a) Compares the POD projection error for the transformed and the untransformed snapshot matrices. (b) Compares the POD projection error for the transformed snapshots for two different grid sizes.}}	
\end{figure}

\Cref{fig: test-1 sigma} presents the results. A few observations are in order. Firstly, for all $n\leq 23$, the relative error is smaller for $\snapMat{G}$. Thus, at least for these smaller values of $n$, snapshot transformation improves the approximability of a snapshot matrix in its POD modes. Secondly, the difference between the two relative errors is dramatic for $n\leq 10$, with $E^{proj}_n(\snapMat{G})$ being at least four time smaller than $E^{proj}_n(\snapMat{U})$. The difference is the most pronounced for $n=1$, for which we find
\begin{gather}
E^{proj}_{n=1}(\snapMat{G}) \approx 0.11,\hspB E^{proj}_{n=1}(\snapMat{U})\approx 0.70. 
\end{gather}
We emphasize that $n=1$ is just $0.1\%$ of the high-fidelity space dimension $N$. Remarkably, for such a small value of $n$, snapshot transformation facilitates a relative error of just $11\%$. Lastly, around $n> 10$, the decay in the relative error $E^{proj}_n(\snapMat{G})$ slows down. The error (almost) starts to stagnate and eventually, at $n\approx 23$, it overshoots $E^{proj}_n(\snapMat{U})$---\cite{RegisterMOR} reports a similar behaviour. 

A plausible reason for this stagnation (as discussed in \cite{SarnaCalib2020}) is that the transformed snapshots have slightly misaligned discontinuities---we further showcase the misalignment below. Beyond a certain $n$, the error from this misalignment dominates the total error, causing (almost) an error stagnation. Usually, the misalignment is $\mcal O(1/N)$---therefore, the point of stagnation decreases upon increasing $N$. \Cref{fig: test-1 sigma diffN} depicts this behaviour. Furthermore, the error at the stagnation point is $\mcal O(1/\sqrt{N})$, which is the same order of accuracy as the HF approximation \cite{SarnaCalib2020}. The implication being that the error introduced by the misalignment is of little practical interest. 

\subsubsection{Study of the transformed solution} 
\begin{figure}[ht!]
\centering
\subfloat[$u_N^{(1)}(x,t,\mu=1)$]{\includegraphics[width=2.5in]{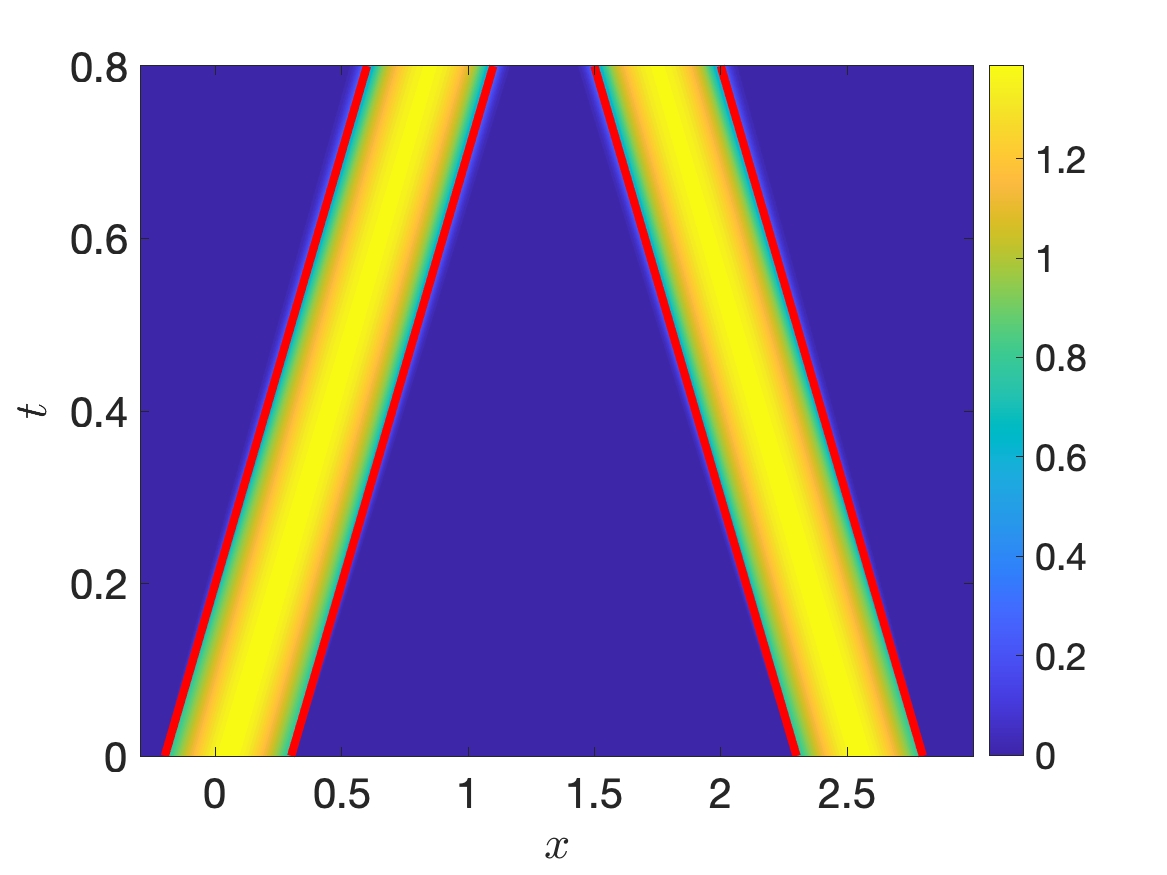}} 
\hfill
\subfloat[$g^{(1)}(x,t,\mu=1)$]{\includegraphics[width=2.5in]{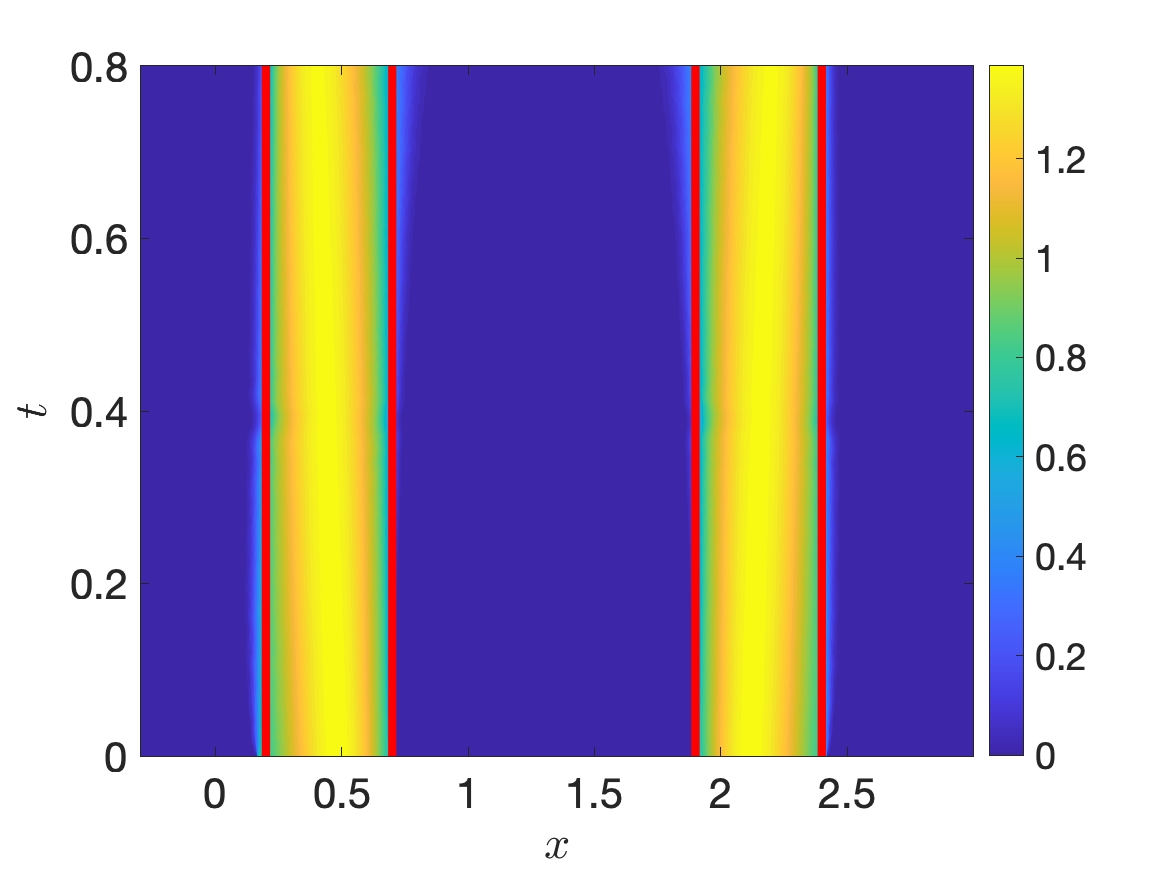}}
\caption{\textit{Results for test-1. Time-trajectory of the (a) untransformed and (b) transformed solution. Location of spatial discontinuities shown in red.}\label{fig: test-1 sol comp}}	
\end{figure}

Let us compare the transformed solution to the untransformed one. We choose $\mu=1$ and study the time-evolution of the transformed $g^{(1)}(\cdot,t,\mu)$ and the untransformed $u^{(1)}_N(\cdot,t,\mu)$ solution. \Cref{fig: test-1 sol comp} presents the results. The spatial discontinuities in $u^{(1)}_N(\cdot,t,\mu)$ (shown in red) move along the time-domain, resulting in poor-approximability in a linear reduced space. However, a composition with a spatial transform almost halts the temporal movement of the discontinuities. This results in a transformed solution that is well-approximable in a linear reduced space. 

It is noteworthy that although the movement of discontinuities in the transformed solution is small, it is not exactly zero i.e., there is some misalignment between the discontinuities. One reason being the $L^2$ objective functional in \eqref{def F}, which, due to its non-convexity, does not guarantee a perfect alignment of spatial discontinuities. To further improve the alignment, one may consider the so-called geometric registration techniques discussed in \cite{ReviewDeformReg}. We plan to pursue such an approach in the future. 

\subsubsection{Study of the average error} 
\begin{figure}[ht!]
\centering
\subfloat[]{\includegraphics[width=2.5in]{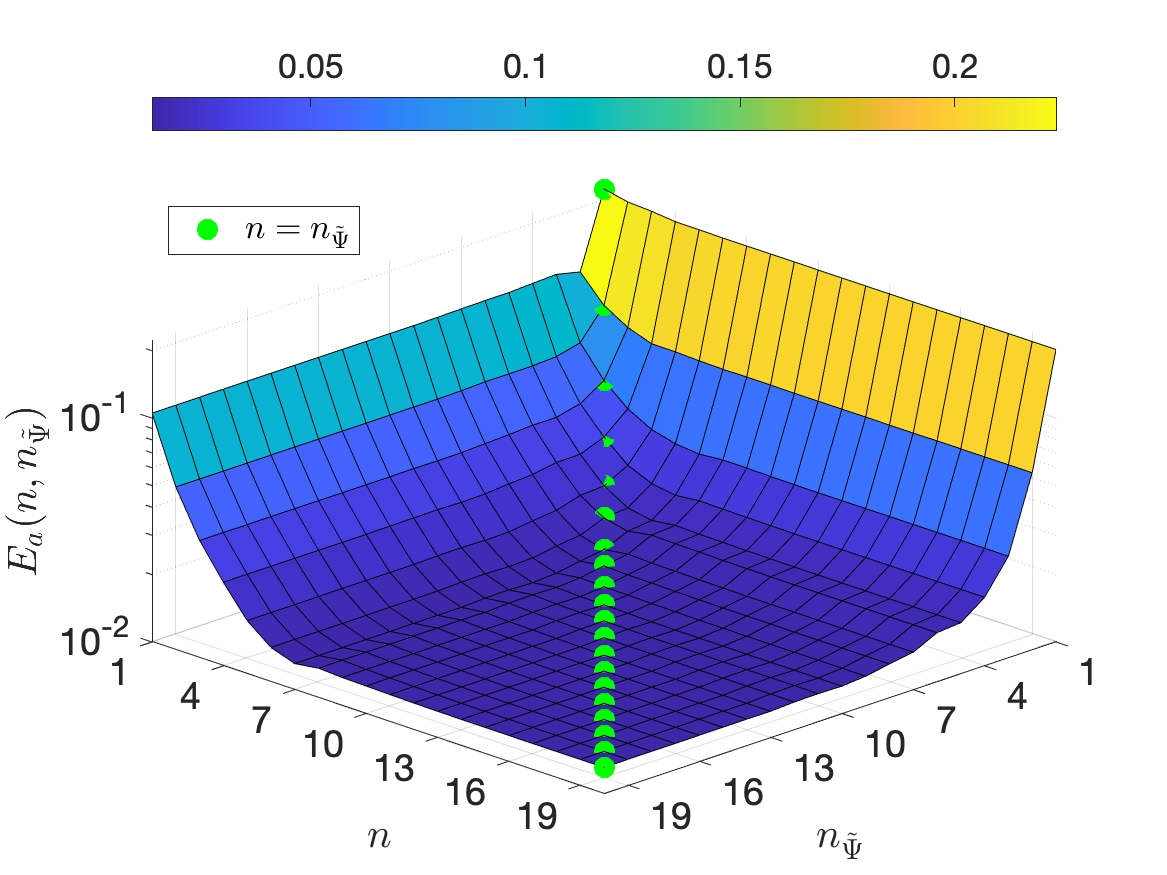}\label{fig: test-1 Eavg}} 
\hfill
\subfloat[]{\includegraphics[width=2.5in]{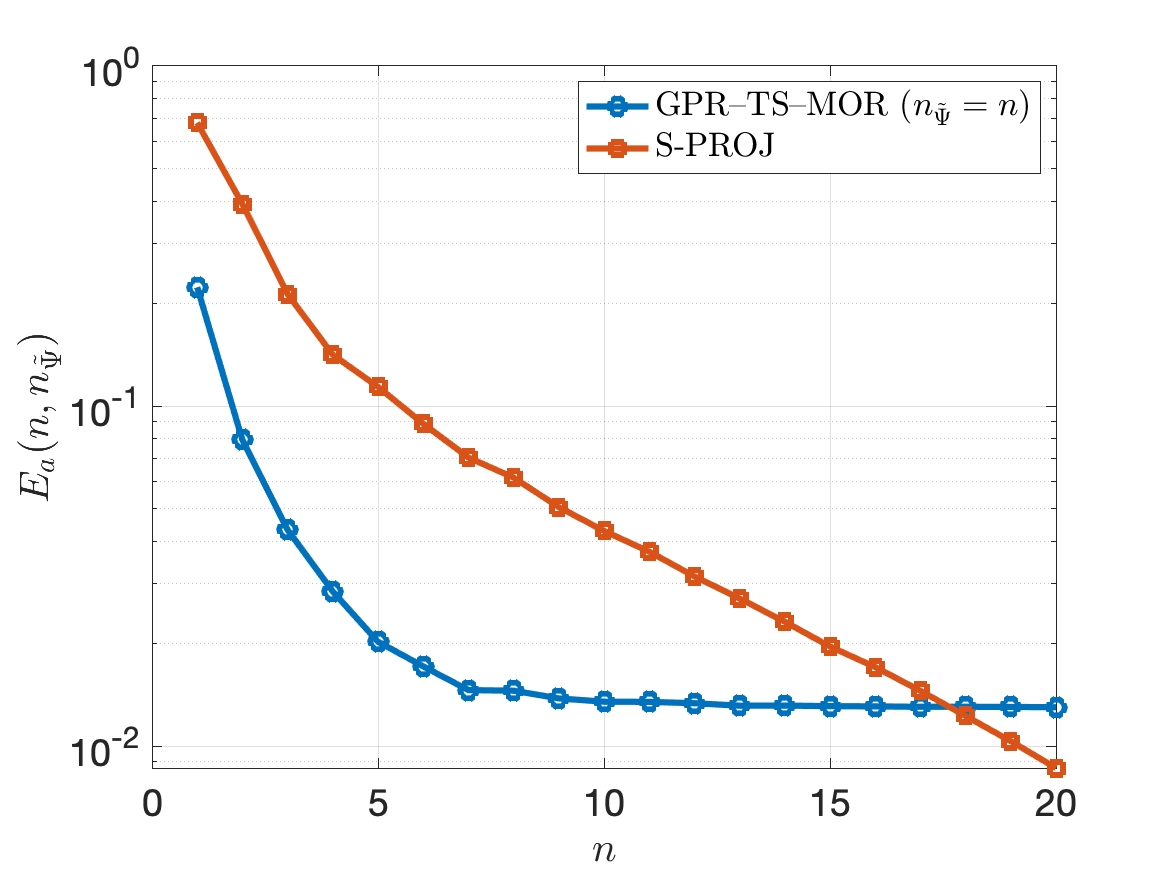}\label{fig: test-1 Ecomp}} 
\hfill
\subfloat[]{\includegraphics[width=2.5in]{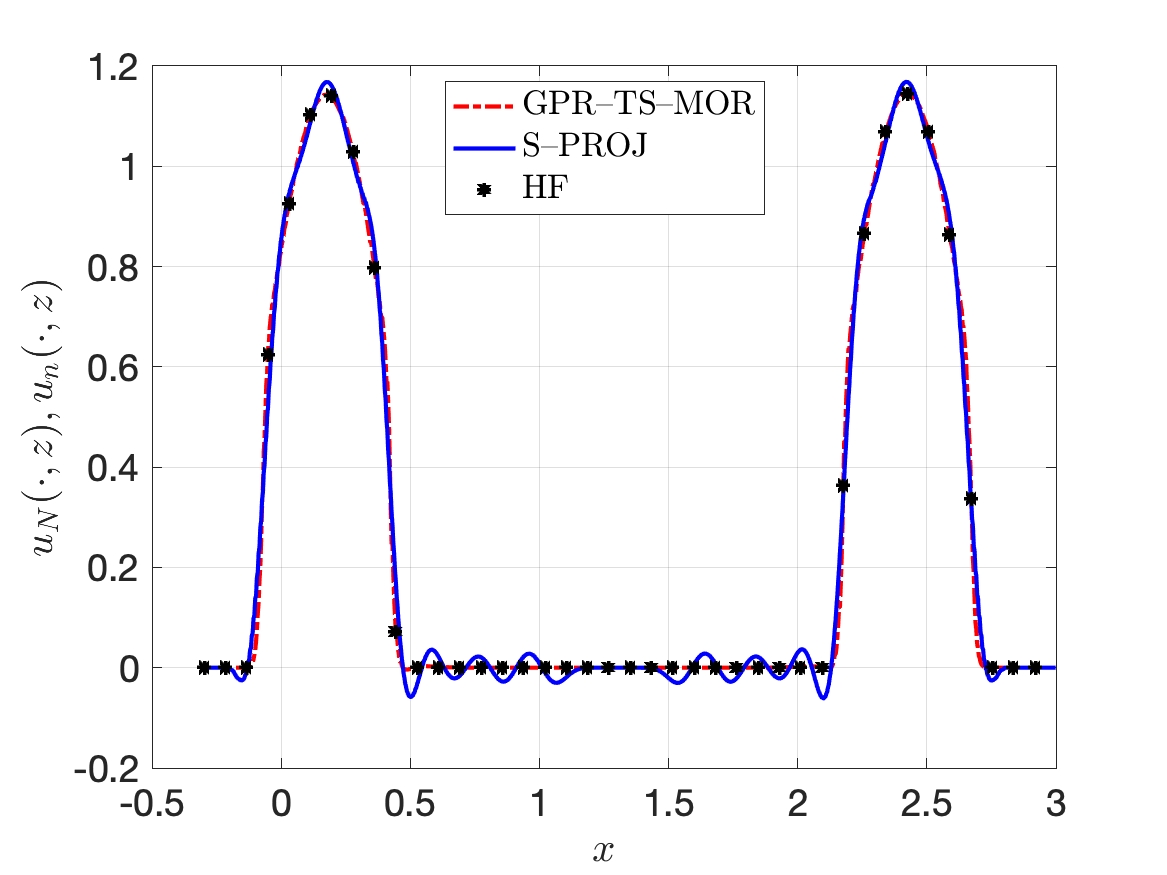}\label{fig: test-1 solcomp}} 
\hfill
\subfloat[]{\includegraphics[width=2.5in]{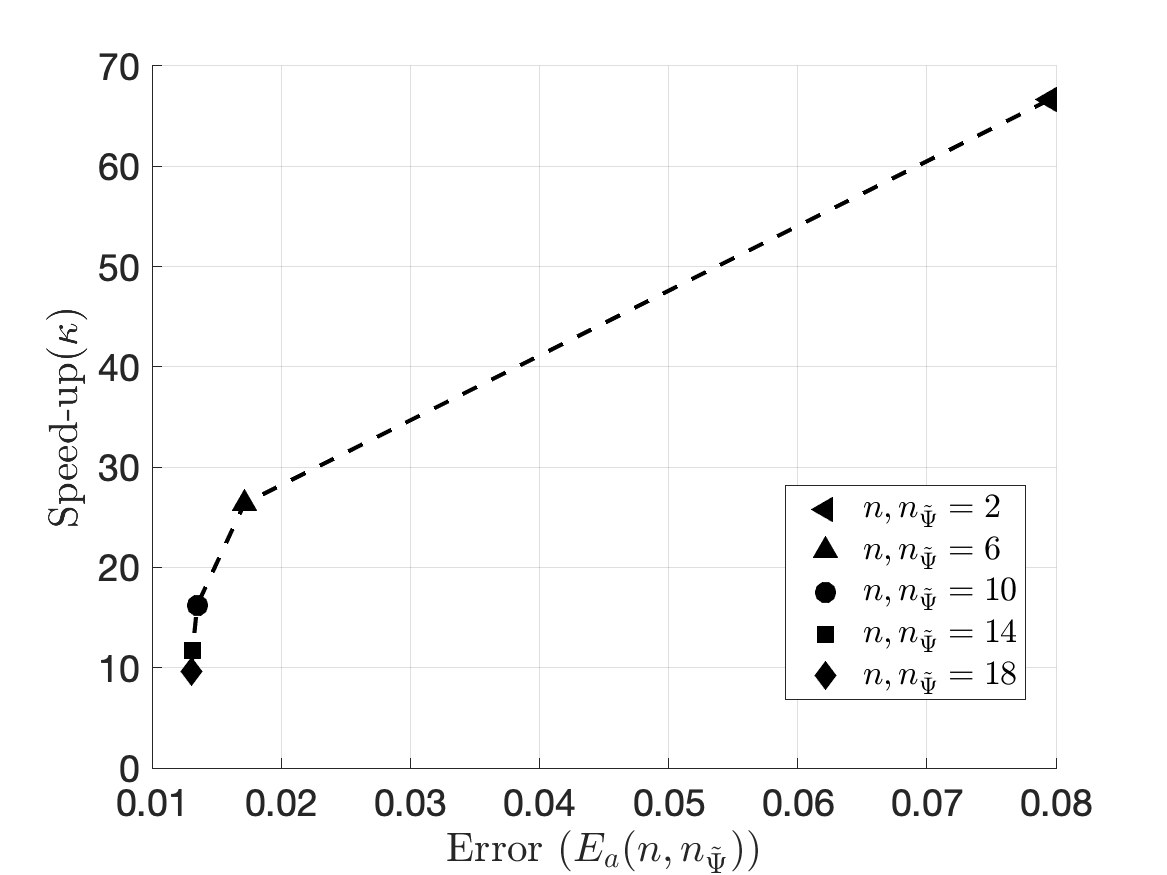}\label{fig: test-1 speedup}} 
\caption{\textit{Results for test-1. (a) Convergence of the average-error with $n$ and $\nPsi$. (b) Error comparison between S-PROJ and GPR-TS-MOR. (c) Solution comparison for $z = (0.8,1)$. (d) Error vs. speed-up for different $n$ and $\nPsi$.}}	\label{fig: test-1 err TSMOR}
\end{figure}

For GPR-TS-MOR, \Cref{fig: test-1 Eavg} presents the average error $E_a(n,\nPsi)$. We vary both $n$ and $\nPsi$ in the set $\{1,2,\dots,20\}$. For smaller values of $n$, the error of approximating the transformed solution outweighs the error of approximating the inverse of the spatial transform and irrespective of how large we make $\nPsi$, it dominates the total error. Consider $n = 1$, for instance, where increasing $\nPsi$ beyond one offers no error decrement---the error stagnates at $\approx11\%$. Similar observation holds for a fixed $\nPsi$ and a variable $n$. 

As anticipated, increasing both $n$ and $\nPsi$ simultaneously, decreases the error monotonically. However, for larger values of $n$ and $\nPsi$, the error almost stagnates at $1.3\%$---\Cref{table: test-1 err} further highlights this stagnation. Following are two plausible reasons. First, the misalignment of spatial discontinuities referred to earlier, which stagnates the singular value decay. Second, due to a limited size of the training set $\mcal D_z$, the GPR can only provide as much accuracy.  The observation that increasing the size of the training data lowers the error at the stagnation point---see \Cref{table: test-1 err}---corroborates our explanation. We emphasize that for practical purposes, an error of $1.3\%$ is reasonable, especially because, compared to the true solution, the HF solver also results in an error of $\approx 2\%$.  Note that other MOR techniques that also rely on regression (driven by neural networks, for instance) report a similar error stagnation \cite{KevinAuto,JanNN2018}. 

For $\nPsi = n$, \Cref{fig: test-1 Ecomp} compares the average error between GPR-TS-MOR and S-PROJ. Recall that for S-PROJ, the value of $\nPsi$ is irrelevant. For $n\lesssim 17$, GPR-TS-MOR outperforms S-PROJ. Already for $n=1$, GPR-TS-MOR results in a relative error of $\approx 22\%$. In comparison, S-PROJ results in an error of $\approx 70\%$. The difference is the most pronounced for $n=5$, for which GPR-TS-MOR results in an error of $2\%$, which is $\approx 5.6$ times smaller than the error resulting from S-PROJ. Due to the stagnation in the POD projection error reported earlier, S-PROJ outperforms GPR-TS-MOR for $n \gtrsim 17$. 

Note that as $n$ is increased, although the solution from S-PROJ appears to converge in the $L^1$-sense, it results in a highly oscillatory approximation---\Cref{fig: test-1 solcomp} compares the different solutions for $z = (0.8,1)$ and $n,\nPsi = 10$. The discontinuities in the (untransformed) solution trigger oscillations in the POD basis, which results in an oscillatory S-PROJ solution. In contrast, owing to the solution transformation, GPR-TS-MOR exhibits no such oscillations. 

\begin{table}
\centering
 \begin{tabular}{c|c|c|c|c|c} 
 & \multicolumn{5}{c}{$E_a(n,\nPsi)$}\\
 \hline
$\#\mcal D_z$  & $\nPsi,n=1$ & $\nPsi,n=3$ & $\nPsi,n=5$ & $\nPsi,n=7$ & $\nPsi,n=9$\\ 
  \hline
 $40\times 20$& $2.2\times 10^{-1}$ & $4.3\times 10^{-2}$ & $2.0\times 10^{-2}$ & $1.4\times 10^{-2}$ & $1.3\times 10^{-2}$\\ 
 \hline
  $80\times 20$& $2.2\times 10^{-1}$ & $4.3\times 10^{-2}$ & $2.0\times 10^{-2}$ & $1.3\times 10^{-2}$ & $1.1\times 10^{-2}$\\ 
 \end{tabular}
 \caption{\textit{Results for test case-1. Average-error for different values of $n$ and $\nPsi$. The set $\mcal D_z$ contains all the training parameters.} }\label{table: test-1 err}
\end{table}

\subsubsection{Speed-up vs. accuracy} We denote the speed-up by $\kappa$ and define it as 
\begin{gather}
\kappa := \frac{\sum_{z\in\mcal D^{tst}_z}\tau_{HF}(z)}{\sum_{z\in\mcal D^{tst}_z}\tau_{MOR}(z)},\label{def kappa}
\end{gather}
where $\tau_{MOR}(z)$ and $\tau_{HF}(z)$ denote the CPU-time required by the GPR-TS-MOR and the high-fidelity solver, respectively, to compute a solution at $z\in\parSp$. We measure the CPU-time with the MATLAB's built-in function \texttt{tic-toc}.

\Cref{fig: test-1 speedup} plots the speed-up against the average error $E_a(n,\nPsi)$. As expected, increasing both $n$ and $\nPsi$ reduces both the error and the speed-up. The minimum speed-up of $\approx 10$ corresponds to $n,\nPsi = 20$ and is associated with an error of $\approx 1.3\%$. The maximum speed-up of $\approx 70$ corresponds to $n,\nPsi = 1$ and is associated with an error of $\approx 22\%$. Note that for time-dependent problems, our reduced approximation does not require any further iterations or time-stepping. To recover an approximation for any $z\in\parSp$, we directly compute the approximation in \eqref{approx u}. This explains why, on average, we observe a significant speed-up compared to a time-stepping based finite-volume solver.

\subsection{Test-2}
We discretize $\Omega$ with $N=300\times 300$ grid cells, which results in a grid size of $\Delta x = 5.3\times 10^{-3}$. The training data $\mcal D_z$ is a set of $100$ uniformly placed points inside $[0,T]$. We approximate the displacement field $\dev{\zRef}{z}$ in the polynomial space $\mcal P_{M=6}$, where the value of $M$ results from the procedure outlined in \Cref{sec: varphi}. 
\subsubsection{Study of the POD projection error} 
\begin{figure}[ht!]
\centering
\includegraphics[width=2.5in]{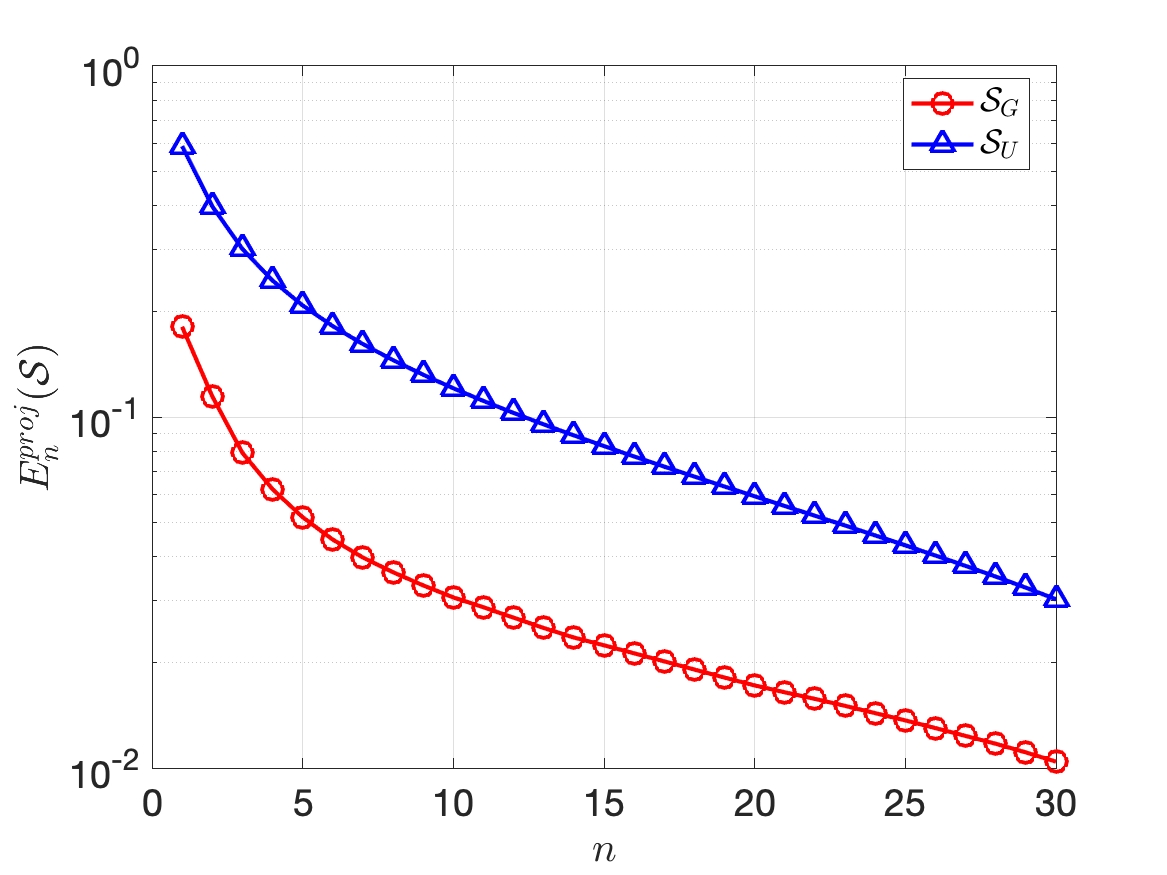} 
\caption{\textit{Results for test-2. Comparison of the POD projection error $E^{proj}_{n}(\mcal S)$ for the transformed $\snapMat{G}$ and the untransformed $\snapMat{U}$ snapshot matrices. }\label{fig: test-2 sigma}}	
\end{figure}

For the two snapshot matrices $\snapMat{G}$ and $\snapMat{U}$ defined in \eqref{def S} and \eqref{def SU}, respectively, \Cref{fig: test-2 sigma} compares the relative POD projection error defined in \eqref{def Eproj}. Similar to the last test case, $\Eproj{\snapMat{G}}$ decays drastically for $n\lesssim 10$, followed by a steady decay for $n\gtrsim 10$. For all $n \leq 30$, it remains at least $2.2$ times smaller than $\Eproj{\snapMat{U}}$. The difference is the most pronounced for $n=13$, for which we find 
\begin{gather}
E^{proj}_{n=13}(\snapMat{G}) \approx 2\times 10^{-2},\hspB E^{proj}_{n=13}(\snapMat{U})\approx 9\times 10^{-2}. 
\end{gather}
Let us recall that $n=13$ is a tiny fraction of $N$-precisely, $1.5\times 10^{-2}\%$ of $N$. For such a small fraction of $N$, a relative error of $2\%$ seems reasonable. 
\subsubsection{Study of the transformed solution}
\begin{figure}[ht!]
\centering
\subfloat[$u_N(\cdot,z=0.2)$]{\includegraphics[width=2.5in]{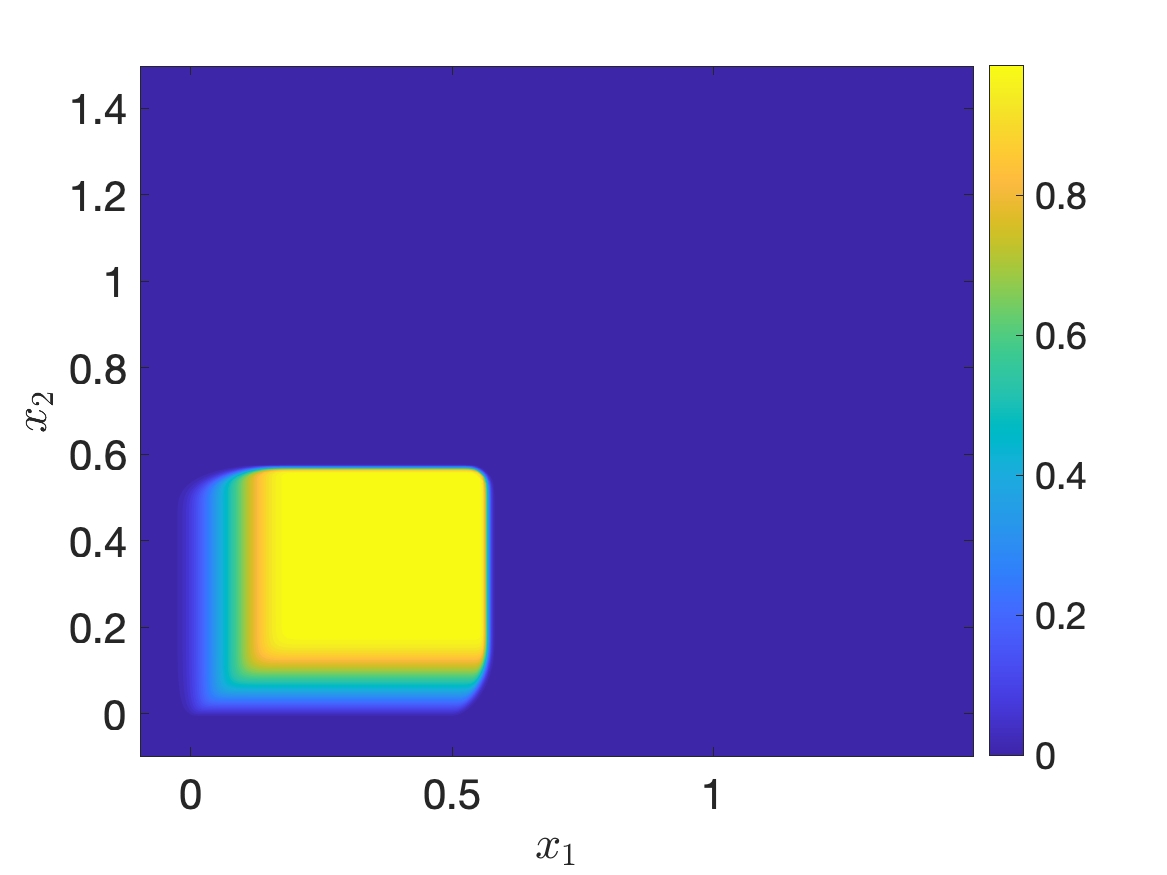}} 
\hfill
\subfloat[$g(\cdot,z=0.2)$]{\includegraphics[width=2.5in]{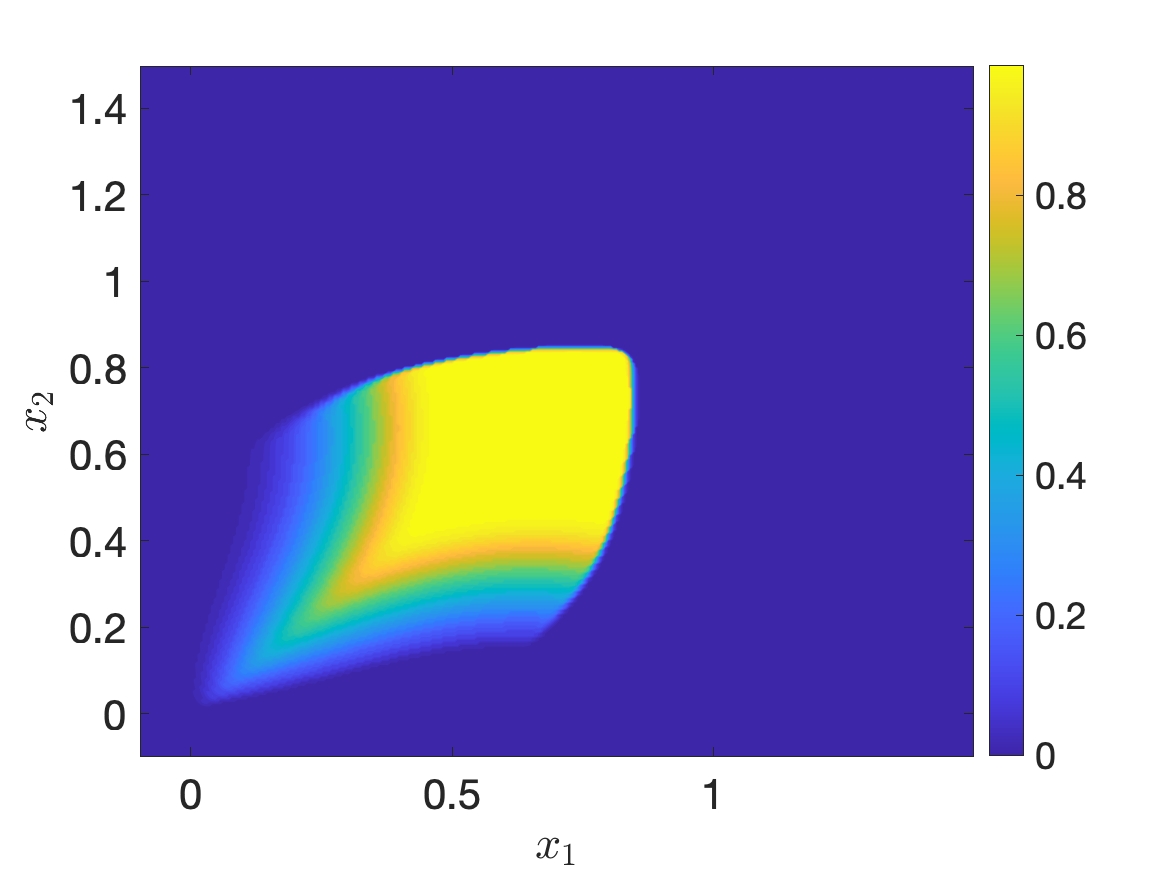}} 
\hfill
\subfloat[$u_N(\cdot,z=1)$]{\includegraphics[width=2.5in]{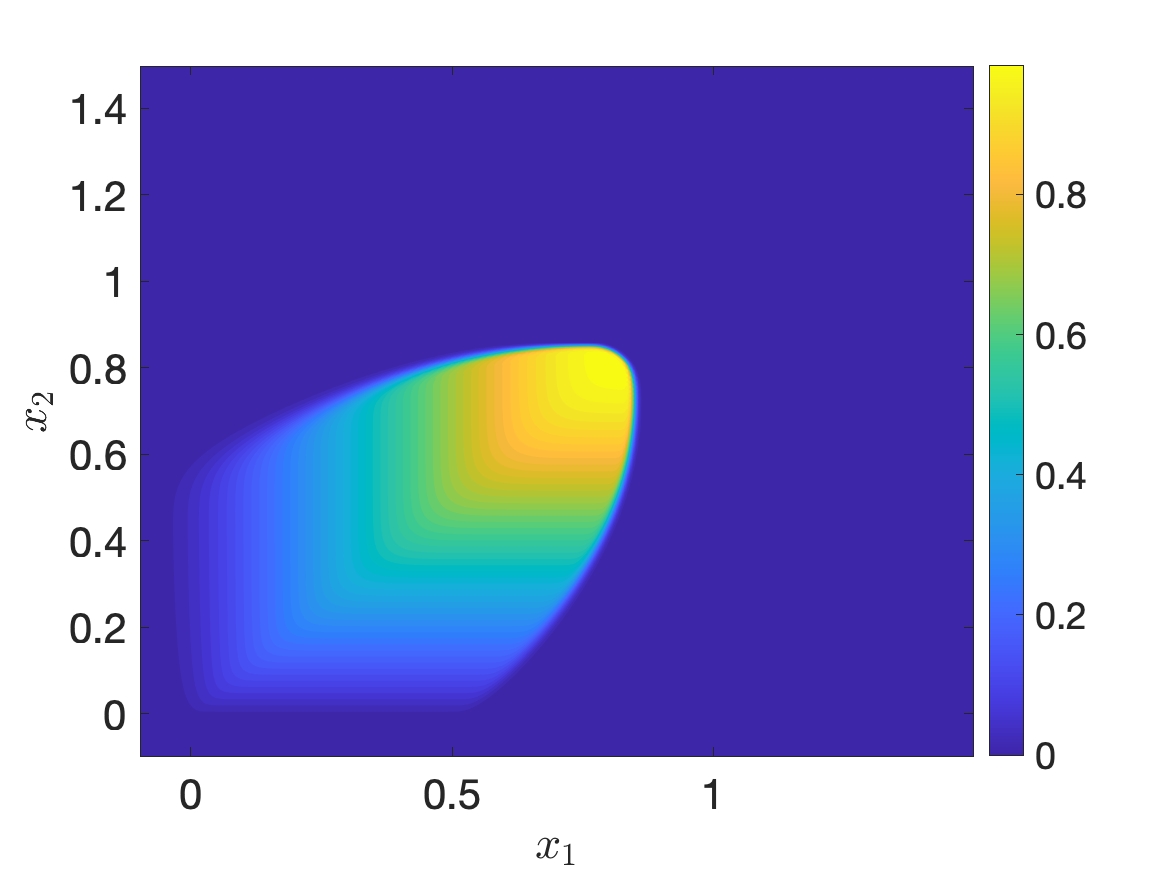}} 
\hfill 
\subfloat[$g(\cdot,z=1)$]{\includegraphics[width=2.5in]{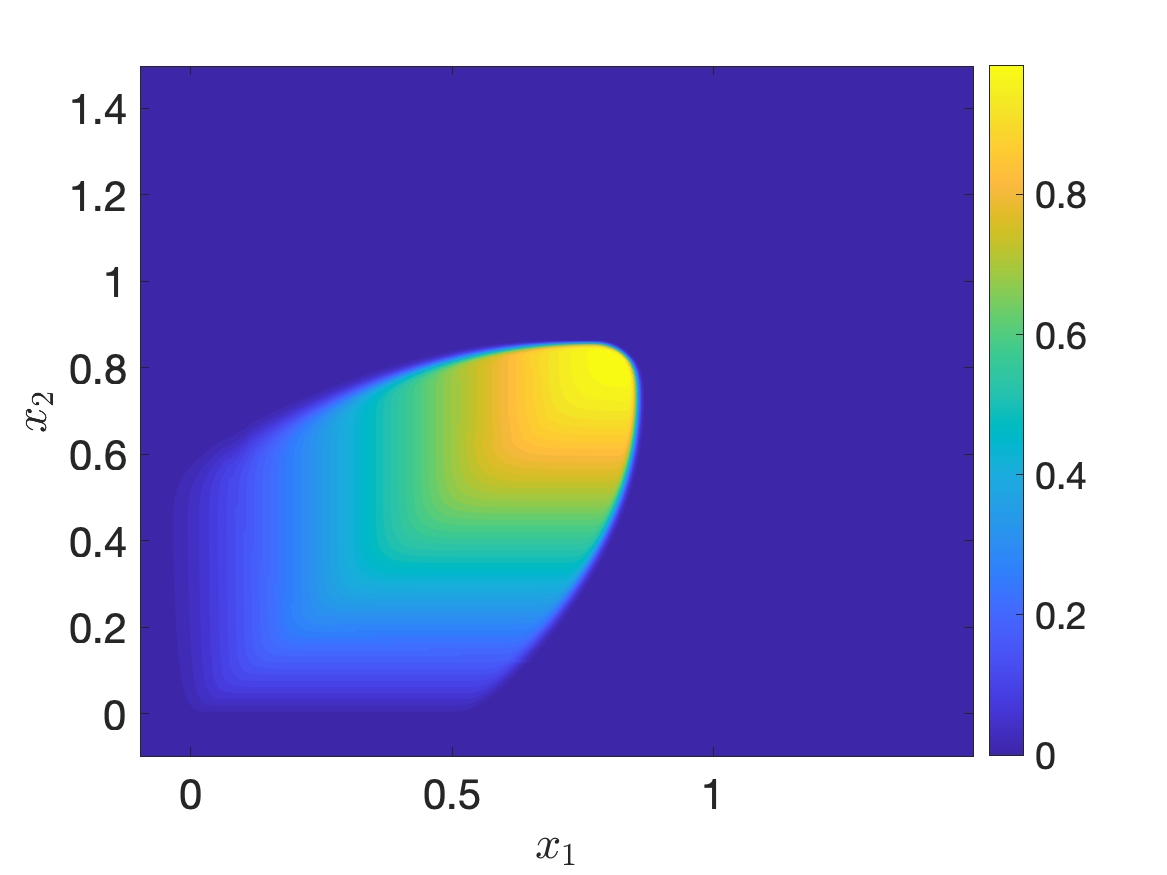}} 
\hfill
\subfloat[$u_N(\cdot,z=1.8)$]{\includegraphics[width=2.5in]{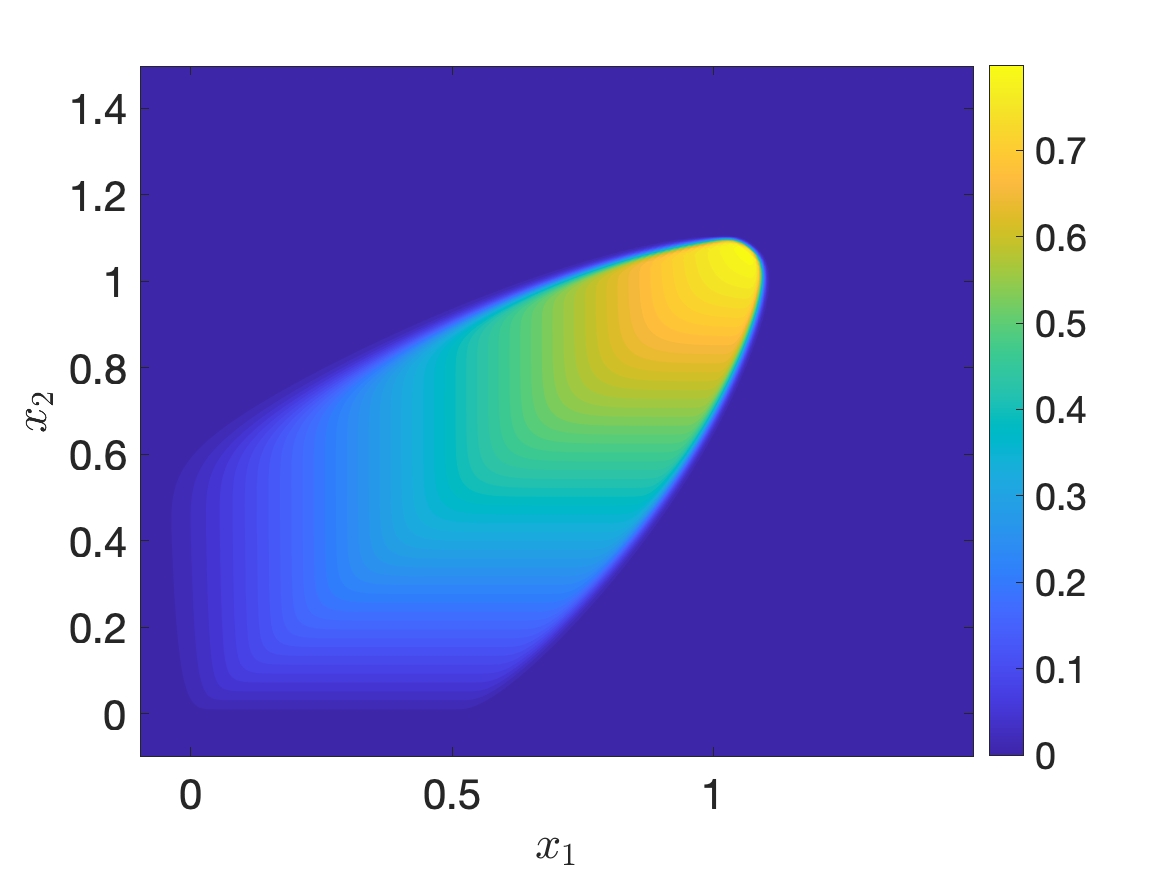}} 
\hfill 
\subfloat[$g(\cdot,z=1.8)$]{\includegraphics[width=2.5in]{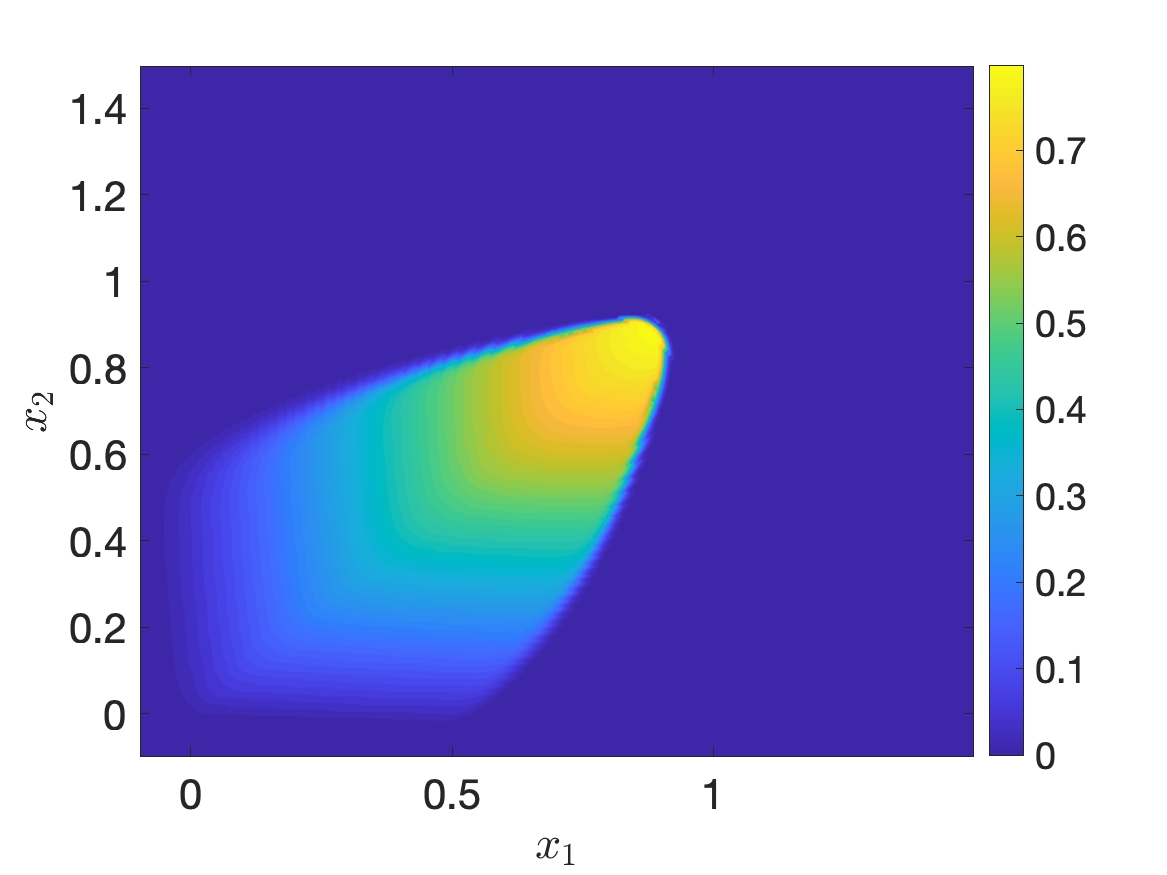}} 
\caption{\textit{Results for test-2. Comparison between the transformed $\solFV{z}$ and the untransformed solution $g(\cdot,z)$. Left panel: untransformed solution. Right panel: transformed solution. }}	\label{fig: test-2 calib sol}
\end{figure}

Let us further elaborate on the reason why the singular values of $\snapMat{G}$ decay faster than those of $\snapMat{U}$.
Consider three different parameter samples: $z_1 = 0.2$, $z_2 = 1$, and $z_3 = 1.8$. For these parameters, \Cref{fig: test-2 calib sol} compares the solution $\solFV{z}$ to the transformed solution $g(\cdot,z)$. We observe that $\solFV{z_2}=g(\cdot,z_2)$. This is because $z_2$ is our reference parameter for which $\spTransf{\zRef}{z_2} = \opn{Id}$. For the untransformed solution $\solFV{z}$, the surface of discontinuity in the spatial domain moves (almost along the diagonal of the spatial domain) as $z$ changes. This movement is the reason why the singular values of $\snapMat{U}$ (reported in \Cref{fig: test-2 sigma}) decay slowly. In comparison, the surface of discontinuity in the transformed solution shows very little movement with $z$, which then induces a fast singular value decay in the snapshot matrix $\snapMat{G}$. 
\subsubsection{Study of $\varphi$ and $\varphi^{-1}$}
Consider the Jacobians 
\begin{gather}
\mcal J(x,z):= \opn{det}(\nabla\varphi(x,\zRef,z)),\hspB \td{\mcal J}(x,z):= \opn{det}(\nabla\td{\varphi}_{\nPsi}(x,\zRef,z)). \label{def jacobians}
\end{gather}
Recall that $\varphi$ results from the optimization problem in \eqref{L2 opt}, and $\td\varphi_{\nPsi}$ is a POD based approximation to $\varphi^{-1}$ and is given in \eqref{approx phiInv}. For the present study, we fix $\nPsi = 20$. Qualitatively, results remain the same for other values of $\nPsi$. 

\begin{figure}[ht!]
\centering
\subfloat[]{\includegraphics[width=2.5in]{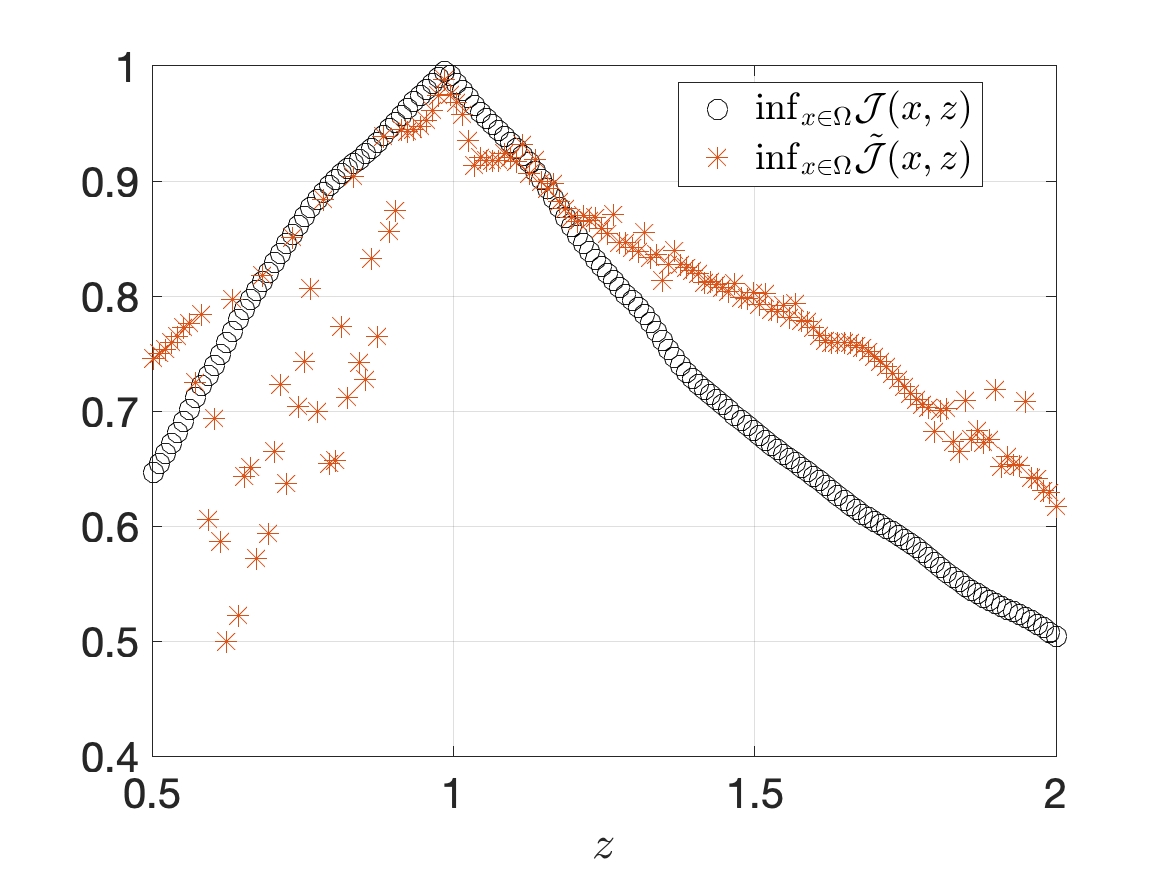}\label{fig: test-2 jacob}} 
\hfill
\subfloat[]{\includegraphics[width=2.5in]{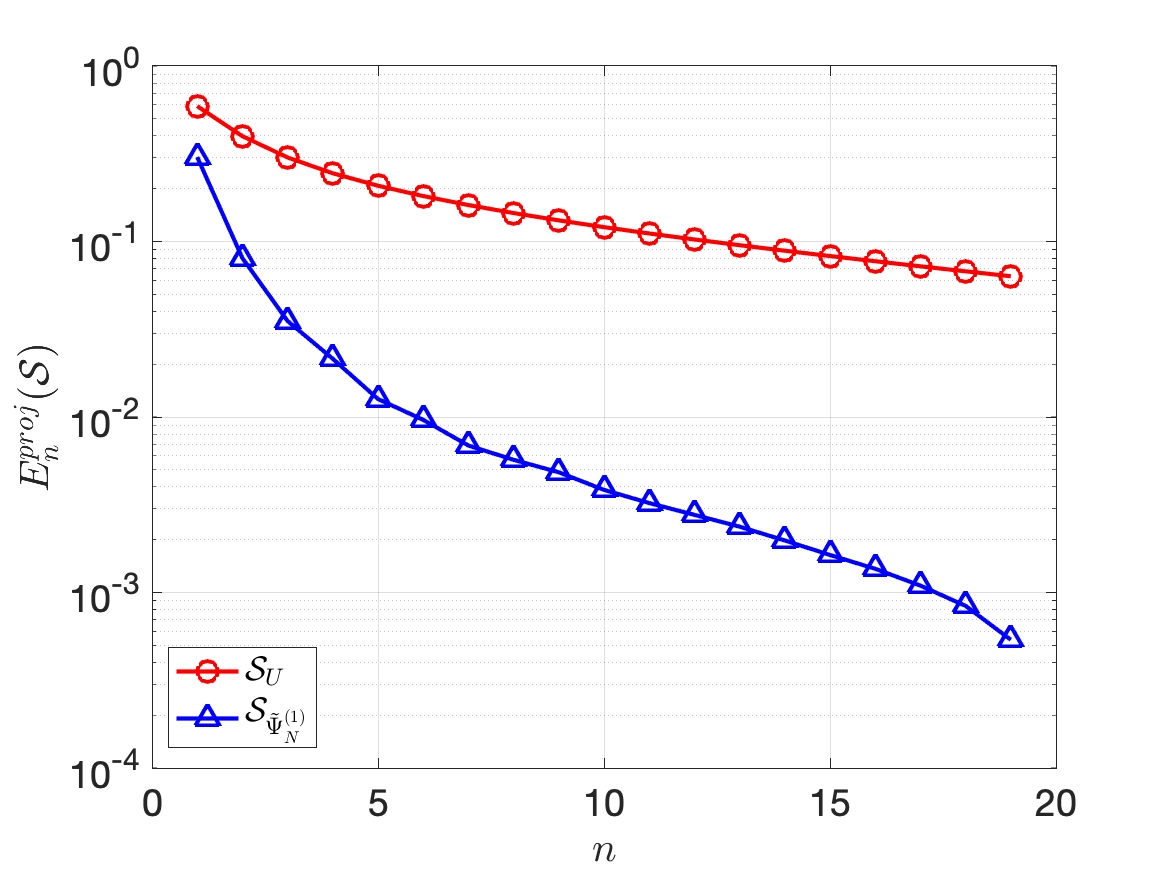}\label{fig: test-2 pod phi}} 
\caption{\textit{Results for test-2. (a) Minimum value of the Jacobian defined in \eqref{def jacobians}. (b) POD projection error for the snapshot matrix $\snapMat{\td\Psi_N^{(1)}}$ and $\snapMat{U}$. }}	
\end{figure}

For different parameter instances, \Cref{fig: test-2 jacob} presents the minimum values of the above two Jacobians. Both the Jacobians stay well above zero, which is desirable. Since $\varphi(\cdot,\zRef,z)\in C^1(\Omega)$, a positive $\mcal J(x,z)$ together with Theorem-2.1 of \cite{RegisterMOR} implies that $\varphi(\cdot,\zRef,z)$ is a diffeomorphism from $\Omega$ to $\Omega$. Equivalently, in the sense of \Cref{remark: diffeo phi}, the present example exhibits a small displacement. Similarly, a positive $\td{\mcal J}(x,z)$ implies that $\td\varphi_\nPsi(\cdot,\zRef,z)$ is a homeomorphism from $\Omega$ to $\Omega$. 

Note that $\zRef = 1$. Therefore, as $z\to 1$, $\spTransf{\zRef}{z}\to \opn{Id}$, which, for all $x\in\Omega$, results in $\mcal J(x,z) \to 1$. However, as $z\to 1$, $\td{\mcal J}(x,z) \not\to 1$. This is because $\td\varphi_\nPsi(\cdot,\zRef,z)$ is only an approximation to $\varphi(\cdot,\zRef,z)^{-1}$. Also note that  $\inf_x \mcal J(x,z)$ is better behaved than $\inf_x \td{\mcal J}(x,z)$ because the former results from directly solving the optimization problem in \eqref{L2 opt}, whereas the latter results from a POD and GPR-based approximation.

Observe that as $|z-\zRef|$ increases, the minimum value $\inf_x \mcal J(x,z)$ (and also, in general, the value $\inf_x \td{\mcal J}(x,z)$) decreases monotonically. This is expected because increasing $|z-\zRef|$ increases the distance between the surface of spatial discontinuities in the solution $\sol{z}$ and the reference $\sol{\zRef}$. Therefore, to move these surfaces closer and align them, the spatial transform needs to displace points by larger distances. This significantly \textit{compresses} some regions of the spatial domain, resulting in smaller Jacobians. We speculate that by increasing $\sup_{z\in\parSp}|z-\zRef|$--i.e., by increasing the size of the parameter domain--one can make $\inf_x \mcal J(x,z)$ negative thus, violating the small displacement assumption referred to in \Cref{remark: diffeo phi}. For such problems, one might have to resort to the large deformation registration considered in \cite{LDDM}. 

In \Cref{sec: approx u}, we assumed that the displacement $\td{\Psi}_N(\cdot,\zRef,z)$ is well approximable in a sufficiently low-dimensional POD space. Here, we justify this assumption empirically. Consider the first component of $\td{\Psi}_N(\cdot,\zRef,z)$ given by $\td{\Psi}^{(1)}_N(\cdot,\zRef,z)$. Let $\snapMat{\td{\Psi}_N^{(1)}}$ be a snapshot matrix for $\td{\Psi}^{(1)}_N(\cdot,\zRef,z)$. For this snapshot matrix, \Cref{fig: test-2 pod phi} presents the POD projection error defined in \eqref{def Eproj}. The projection error decays fast. Already for $\nPsi = 5$, we get a relative error of less than $1\%$. The fast decay becomes obvious when we compare with the projection error of the snapshot matrix $\snapMat{U}$ defined in \eqref{def SU}. This matrix contains untransformed snapshots and thus, has a slow decay in the projection error.

\subsubsection{Study of the average error} 
\begin{figure}[ht!]
\centering
\subfloat[]{\includegraphics[width=2.5in]{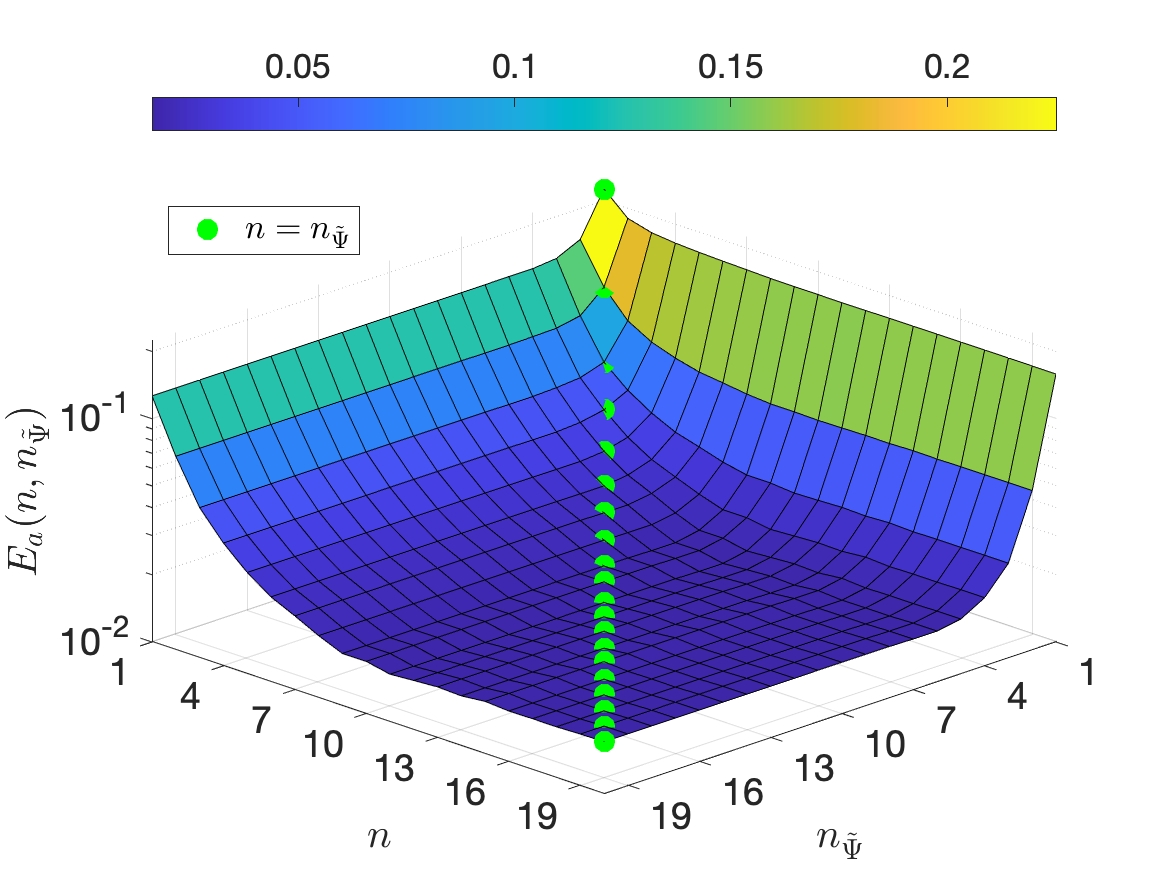}\label{fig: test-2 Eavg}} 
\hfill
\subfloat[]{\includegraphics[width=2.5in]{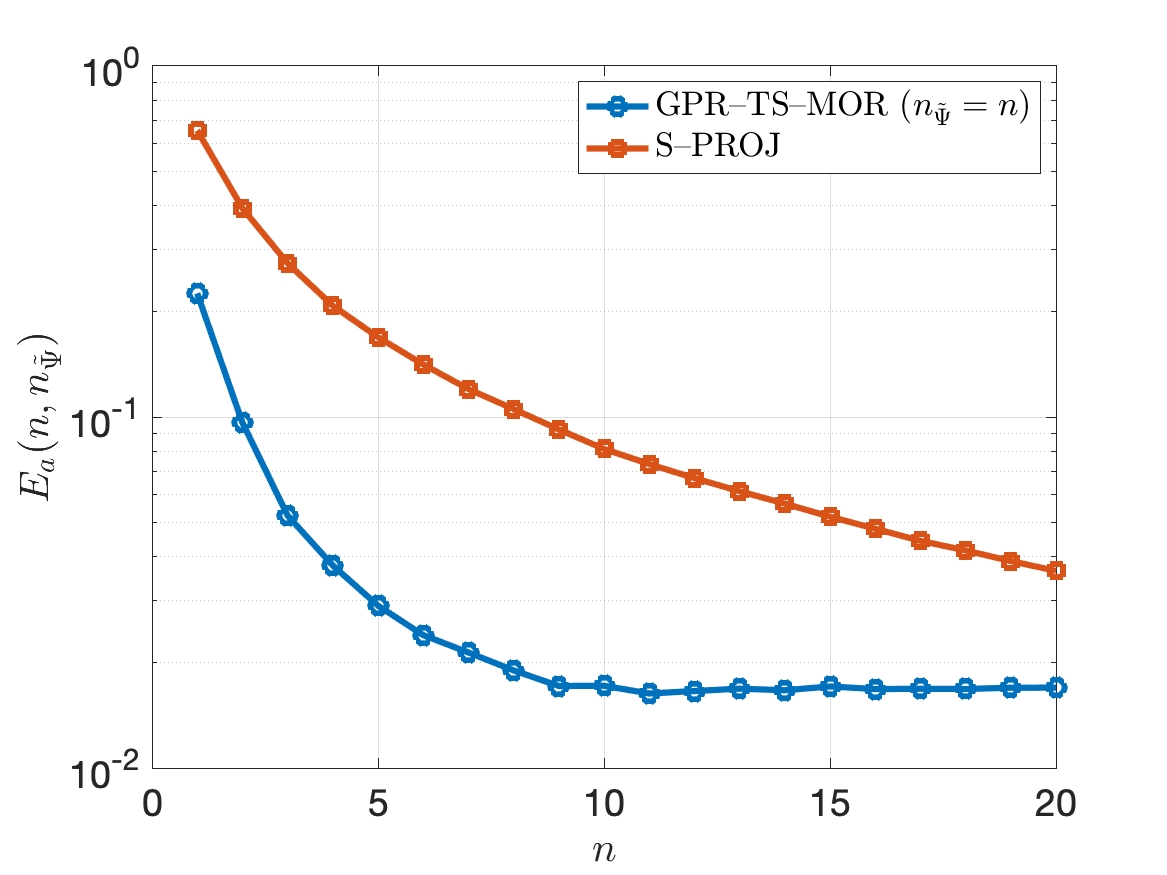}\label{fig: test-2 Ecomp}} 
\caption{\textit{Results for test-2. (a) Convergence of the average error for GPR-TS-MOR with $n$ and $\nPsi$. (b) Error comparison between S-PROJ and GPR-TS-MOR.}}	
\end{figure}

For GPR-TS-MOR, \Cref{fig: test-2 Eavg} presents the average error $E_a(n,\nPsi)$. We vary both $n$ and $\nPsi$ in the set $\{1,2,\dots,20\}$. Results are similar to that of the previous test case. Keeping $n$ fixed and increasing $\nPsi$ beyond a certain point (or vice-versa) offers no benefit. Nevertheless, increasing both $n$ and $\nPsi$ simultaneously decreases the error monotonically. However, due to the reasons explained earlier, eventually, the error starts to stagnates. The lowest relative error we attain is of $1.6\%$. 

\Cref{fig: test-2 Ecomp} compares the average error between GPR-TS-MOR and S-PROJ. For all values of $n\leq 20$, GPR-TS-MOR outperforms S-PROJ. It results in an error that is at least two times smaller than the error resulting from S-PROJ. The difference is the most pronounced for $n=6$; GPR-TS-MOR results in an error of $\approx 2.3\%$, which is $\approx 6$ times smaller than the error resulting from S-PROJ.  Let us recall that $n=6$ is just $5.6\times 10^{-3}\%$ of $N$. For such small value of $n$, an error of $2.3\%$ can be considered reasonable. 

The error from GPR-TS-MOR stagnates at a value of $1.6\%$. The minimum value of $n$ that provides this error is $n=11$. Although not shown in the plot, to achieve the same error, S-PROJ requires $50$ POD modes. This is four times the number of modes required by GPR-TS-MOR. The takeaway being that despite the stagnation, GPR-TS-MOR significantly outperforms S-PROJ.

\subsubsection{Solution comparison}
\begin{figure}[ht!]
\centering
\subfloat[GPR-TS-MOR ($z=1.5$)]{\includegraphics[width=2.5in]{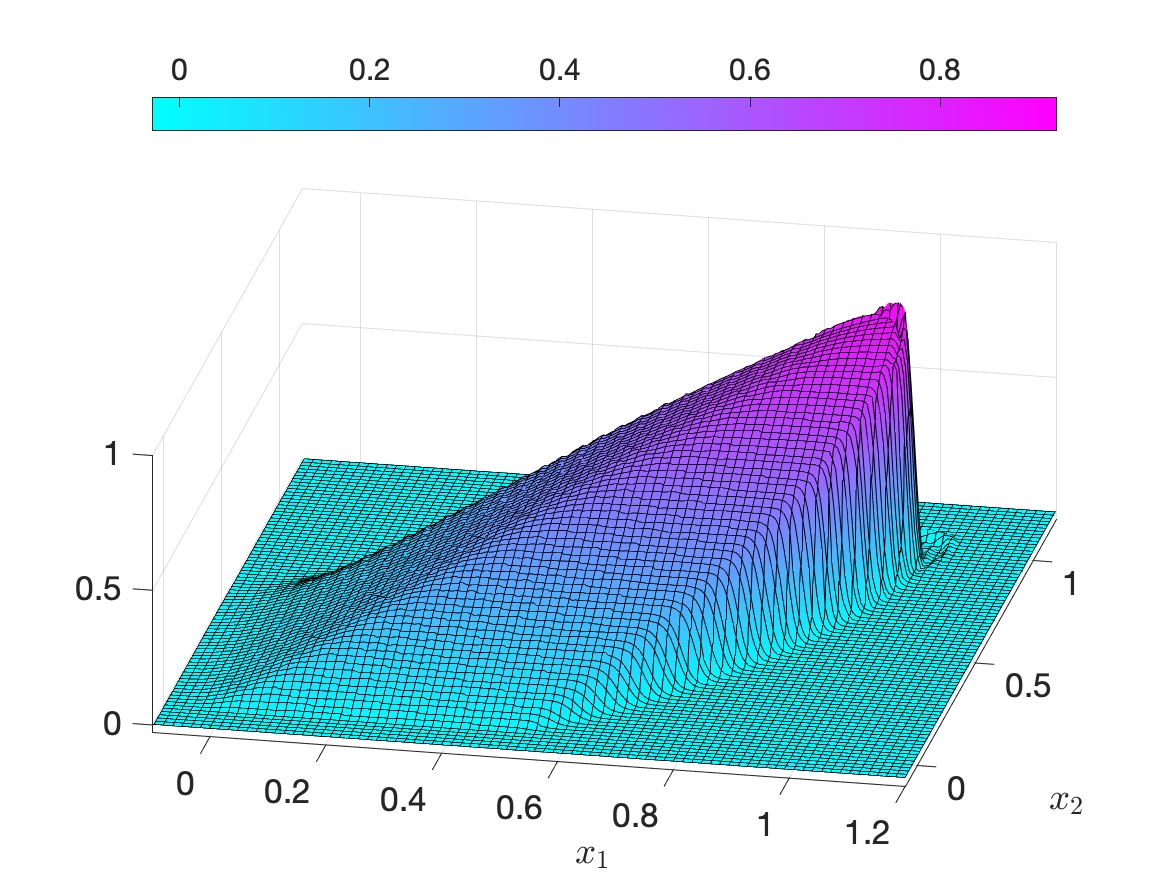}} 
\hfill
\subfloat[S-PROJ ($z=1.5$)]{\includegraphics[width=2.5in]{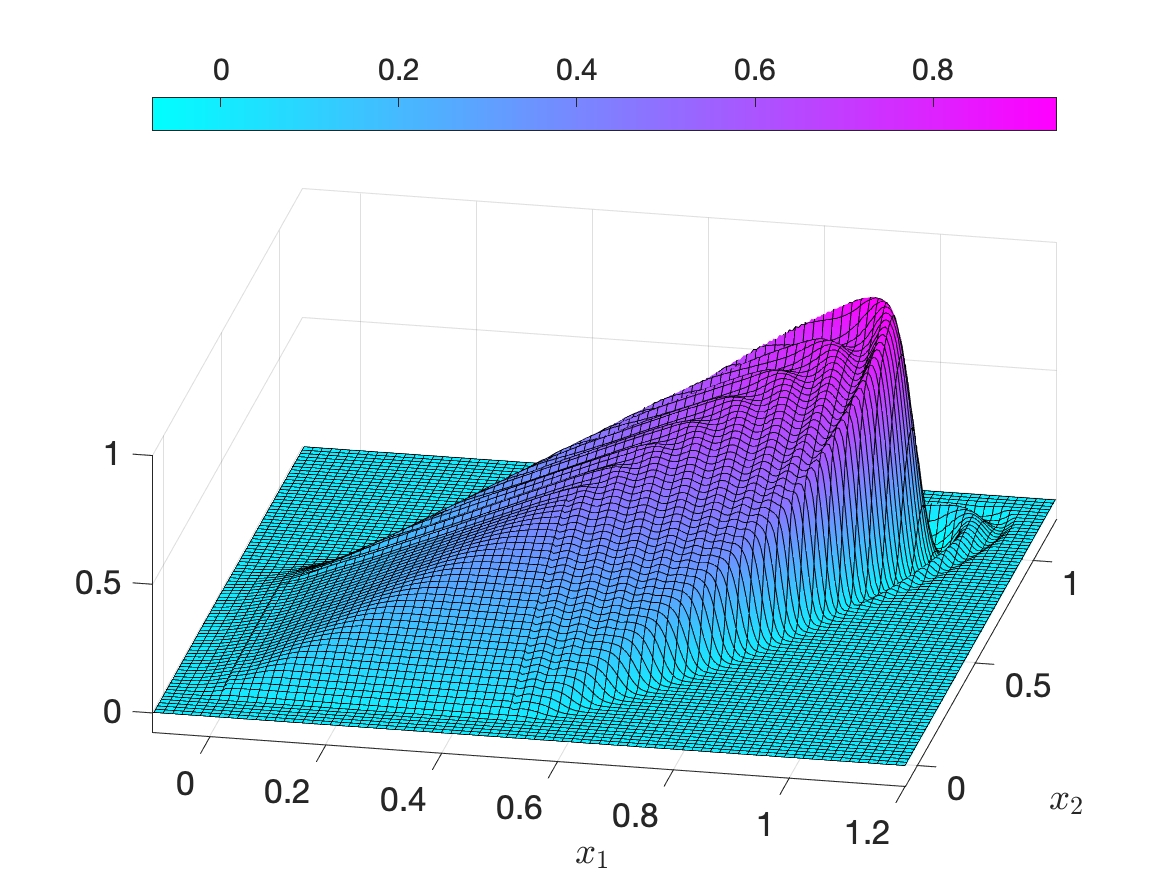}} 
\hfill
\subfloat[HF ($z=1.5$)]{\includegraphics[width=2.5in]{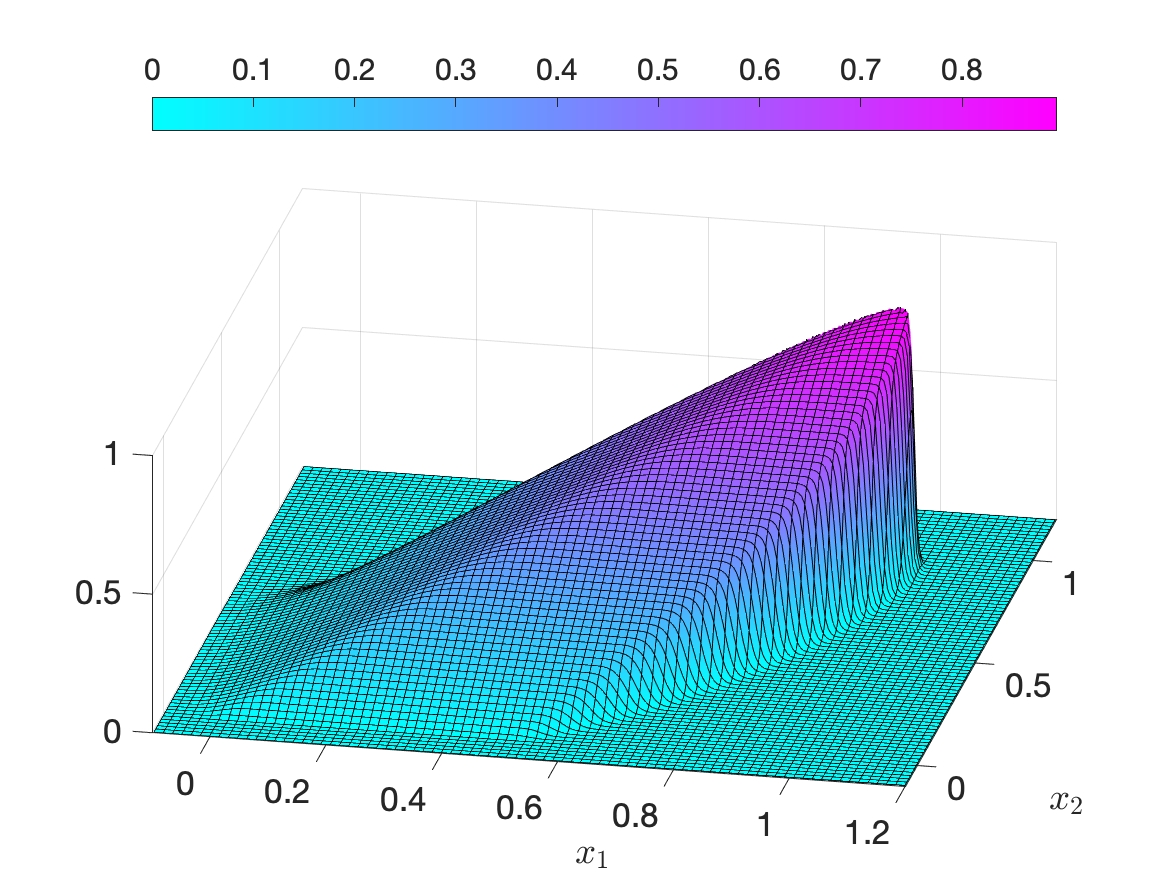}} 
\hfill 
\subfloat[Solution along the line $x_2 = x_1$]{\includegraphics[width=2.5in]{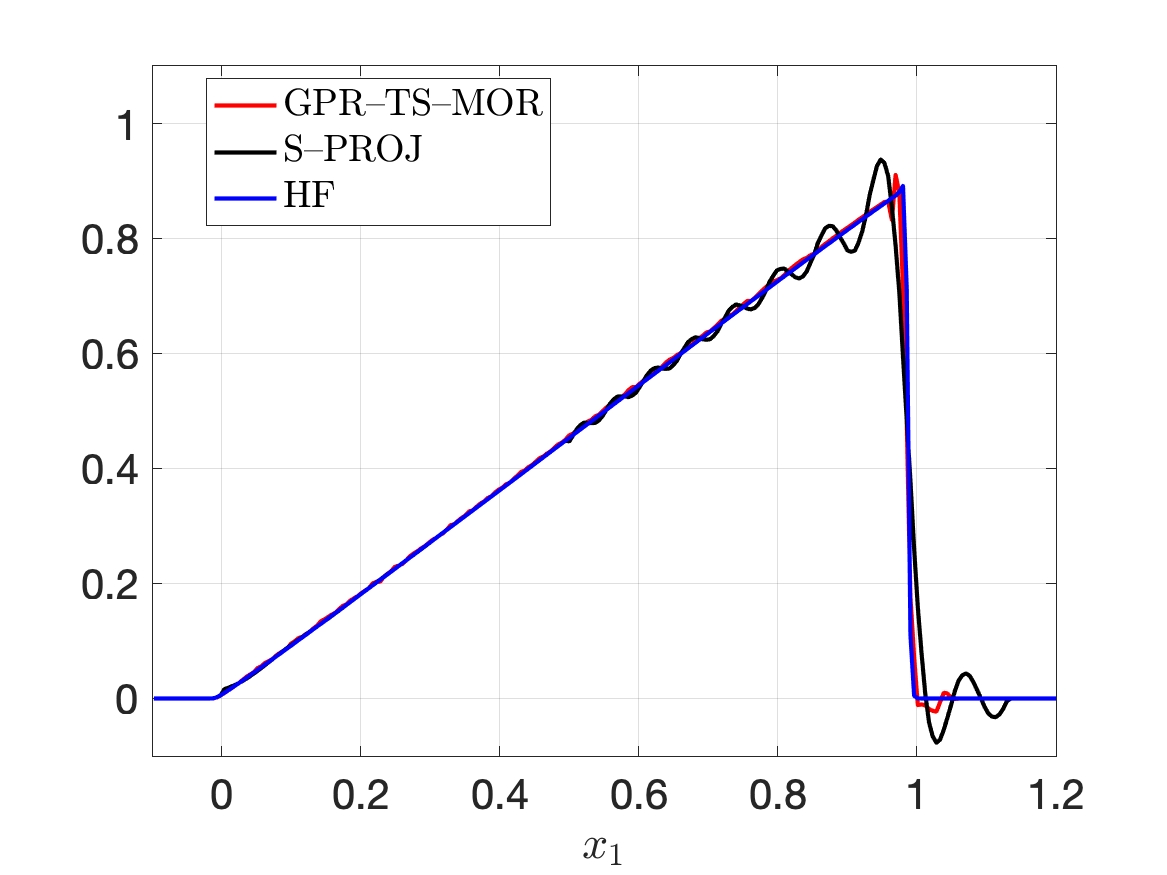}} 
\hfill 
\caption{\textit{Results for test-2. Comparison between the different solutions at $z = 1.5$.}}	\label{fig: test-2 sol comp}
\end{figure}

In \Cref{fig: test-2 sol comp}, we visually compare the solution resulting from GPR-TS-MOR and S-PROJ. We set $z= 1.2$, and $n,\nPsi = 10$.
As anticipated, due to the moving discontinuities in the solution, S-PROJ results in an oscillatory solution. These oscillations are spread-out over the entire spatial domain and appear to originate close to the discontinuity. In contrast, GPR-TS-MOR exhibits no such oscillations. The reason being that it approximates the transformed solution in the POD basis. As studied earlier, discontinuities in the transformed solution do not \textit{move} much, resulting in non-oscillatory POD modes. Observe that GPR-TS-MOR exhibits some minor over and under shoots near the front-end of the surface of discontinuity. The reason being the minor misalignment of discontinuities reported earlier. 

\subsubsection{Study of the error surrogate}
\begin{figure}[ht!]
\centering
\subfloat[]{\includegraphics[width=2.5in]{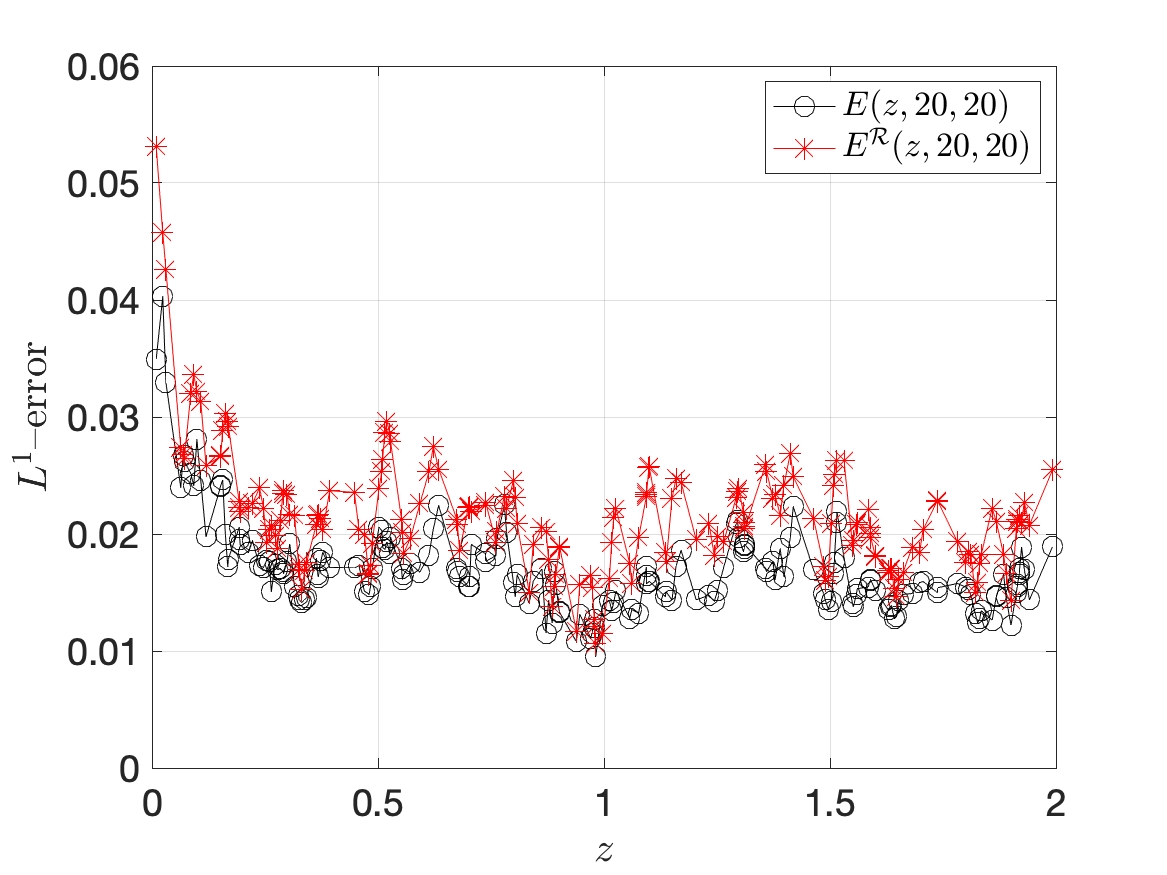}\label{fig: test-2 err comp}} 
\hfill 
\subfloat[]{\includegraphics[width=2.5in]{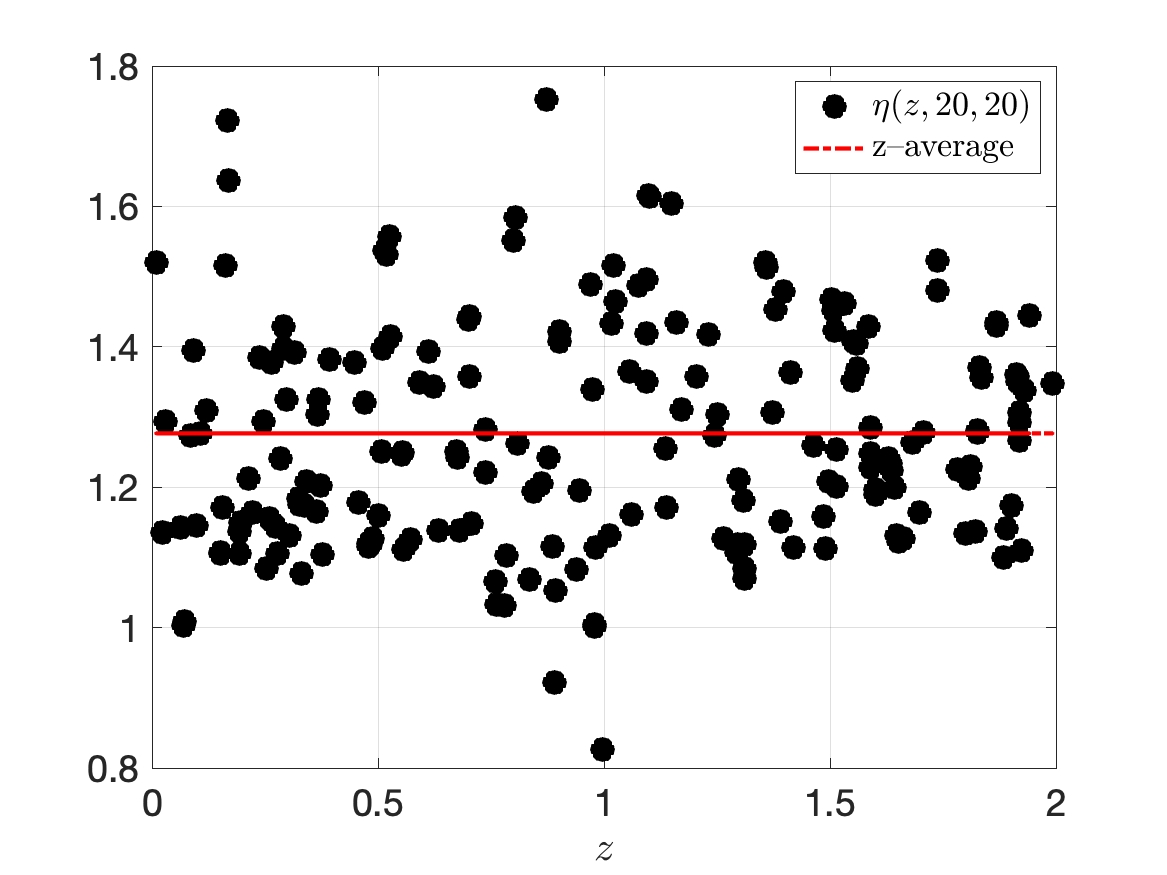}\label{fig: test-2 eff}} 
\caption{\textit{Results for test-2. (a) Compares the true average-error to its surrogate. (b) Presents the efficiency index defined in \eqref{eff index}.  }}	
\end{figure}

For $n,\nPsi = 20$, \Cref{fig: test-2 err comp} compares the true error $E(z,n,\nPsi)$ to its surrogate $E^{\mcal R}(z,n,\nPsi)$ defined in \eqref{def ER}. For most parts of the parameter domain, the surrogate well-approximates the error. This results in an efficiency index (see \Cref{fig: test-2 eff}) that stays close to one. The efficiency index fluctuates between $0.8$ and $1.7$. The average value of the efficiency index is $1.2$ i.e., on average, we over-estimate the error by $20\%$.

Observe that the error is particularly large close to $z=0$. This is because the solution at $z=0$ has a different discontinuity-topology than all the other solution instances---see \Cref{sec: snap transformed} for the relevance of discontinuity-topology. The initial data has two surfaces of discontinuities: (i) the lower and the left edge of the square over which the characteristic function in \eqref{ic burgers} is defined, and (ii) the top and the right edge of the same square. For $z\neq 0$, the first surface manifests into a rarefaction fan, which is continuous. The second surface, however, results in a moving shock. This sudden breakdown of the discontinuity-topology at $z=0$ results in a slightly inaccurate solution transformation close to $z=0$, which, then, results in comparatively larger error values. 
\subsubsection{Speed-up vs. accuracy}
\begin{figure}[ht!]
\centering
\includegraphics[width=2.5in]{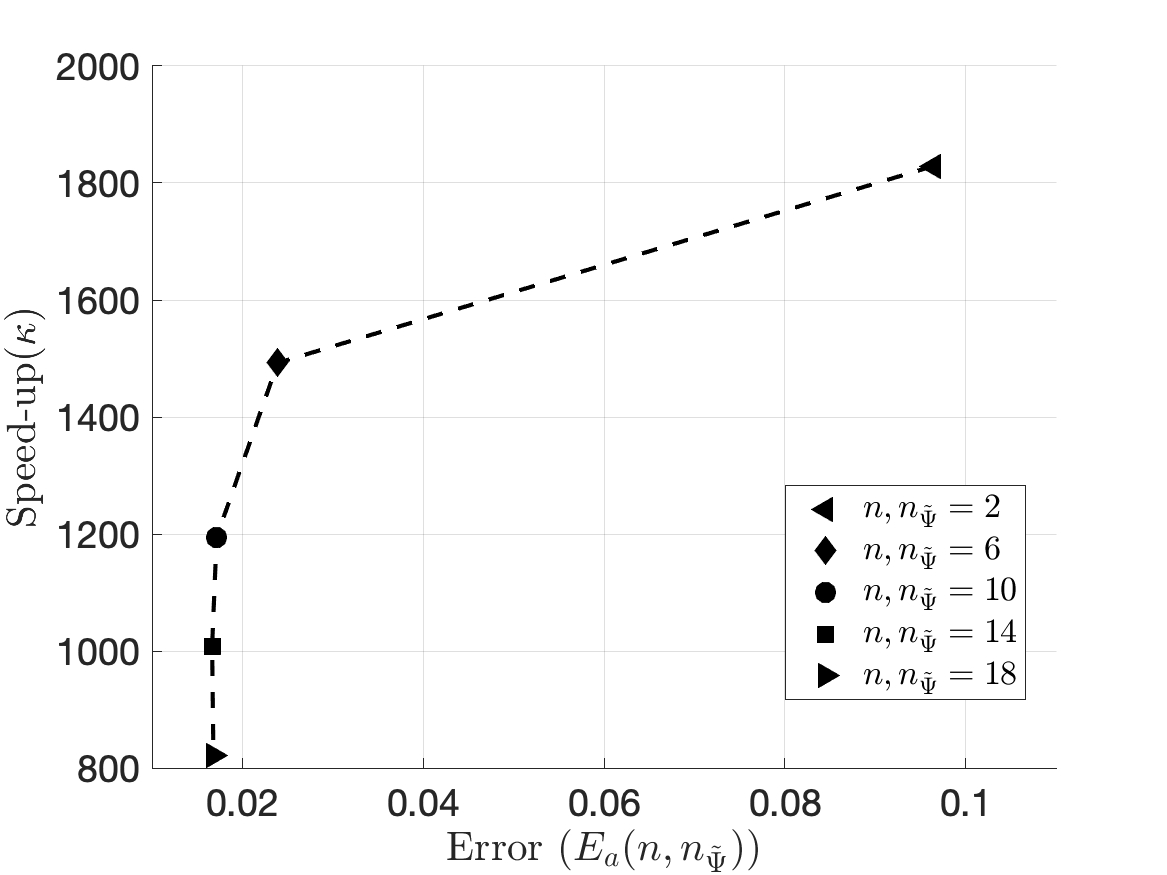}
\caption{\textit{Results for test-2. Speed-up vs. average error plot.}\label{fig: test-2 speedup}}	
\end{figure}

For the GPR-TS-MOR, \Cref{fig: test-2 speedup} plots the speed-up defined in \eqref{def kappa} against the average error. As anticipated, both the speed-up and the error decrease upon simultaneously increasing $n$ and $\nPsi$. The minimum speed-up of $800$ and an error of $1.6\%$ corresponds to $n,\nPsi = 20$. The maximum speed-up of $2000$ and an error of $12\%$ corresponds to $n,\nPsi=1$. Note that the speed-up is at least two orders-of-magnitude larger than in the previous test case. Since the current problem is two-dimensional, the HF space has a large dimension of $300^2$, which makes a HF solver much more expensive than in the previous test case. 

\subsection{Test-3} We discretize $\Omega$ with $300\times 300$ grid cells. The training data $\mcal D_z$ is a set of uniformly placed $40$ points inside $\parSp$. We approximate the displacement field $\dev{\zRef}{z}$ in the polynomial space $\mcal P_{M=6}$, where the value of $M$ results from the procedure outlined in \Cref{sec: varphi}. 

\subsubsection{Solution comparison, average error and speed-up}
\begin{figure}[ht!]
\centering
\subfloat[Error comparison]{\includegraphics[width=2.5in]{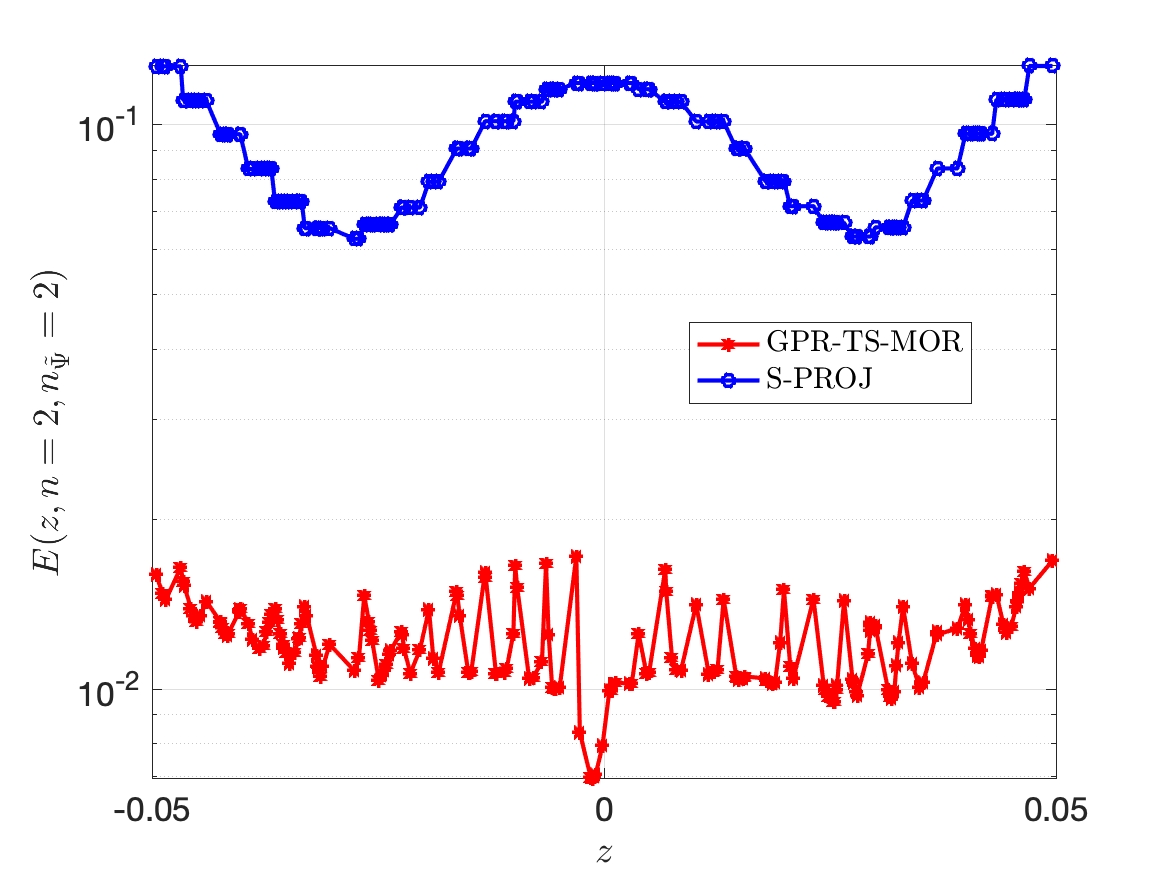}\label{fig: test-3 error}}
\hfill
\subfloat[True vs. predicted error]{\includegraphics[width=2.5in]{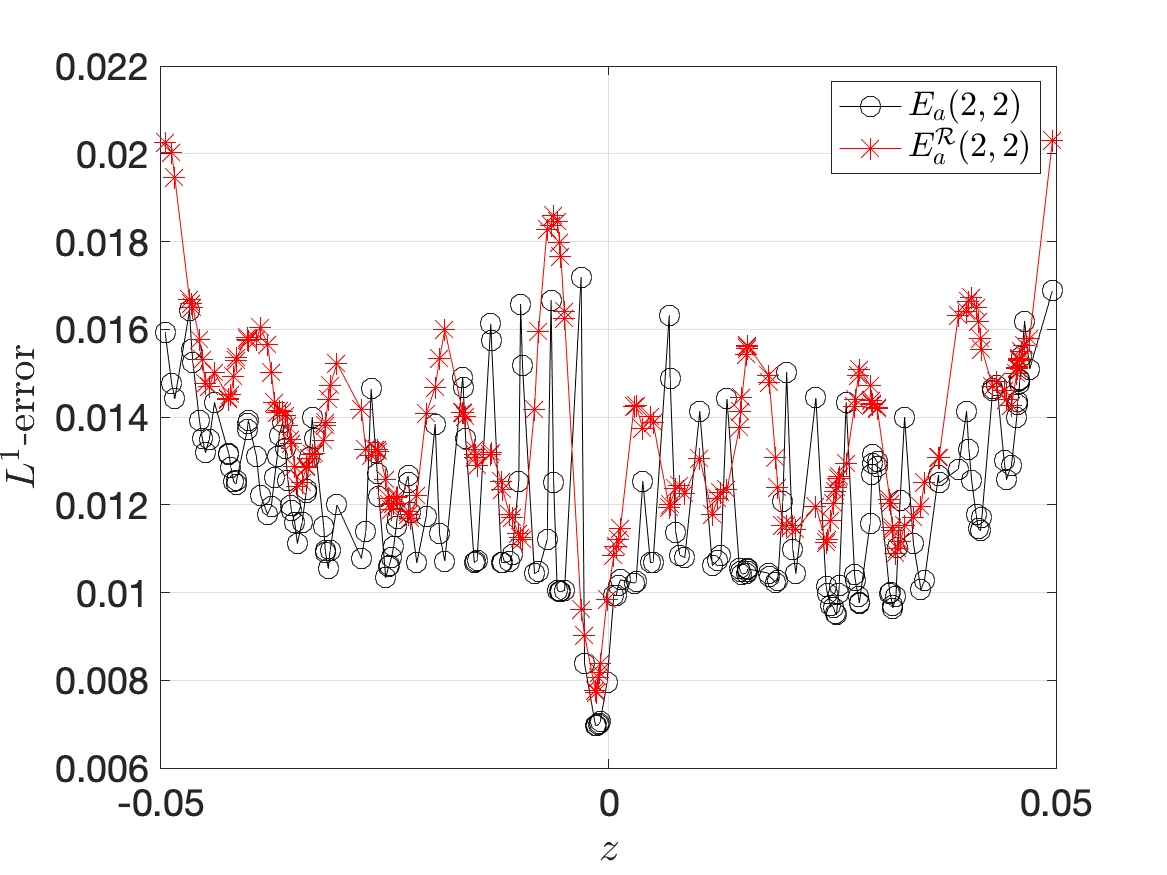}\label{fig: test-3 err comp}} 
\hfill 
\subfloat[Efficiency index for error prediction]{\includegraphics[width=2.5in]{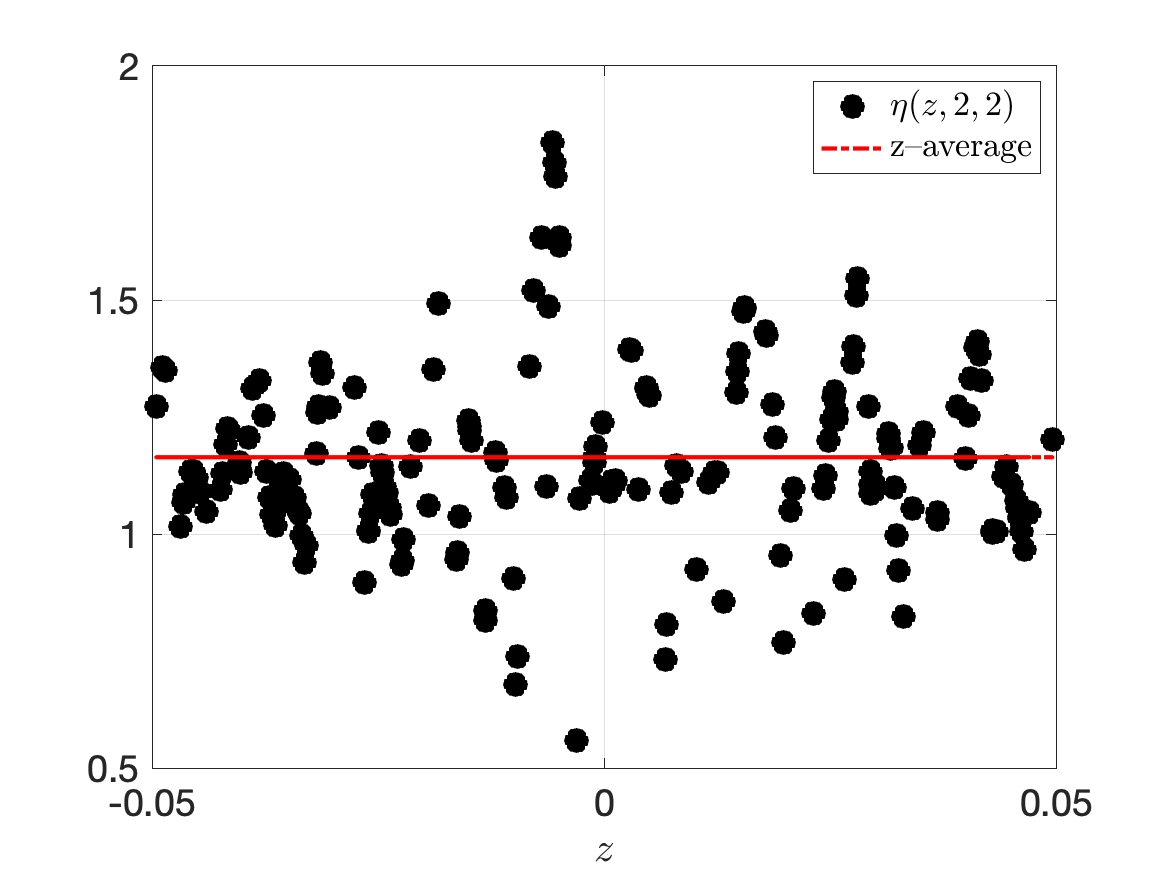}\label{fig: test-3 eff}} 
\caption{\textit{Results for test-3.}}	
\end{figure}

Results for the convergence study remain similar to the previous test case and we do not repeat them here for brevity. We set $n,\nPsi = 2$, study the resulting average error and perform solution comparison. As for the average error, we find
\begin{equation}
\begin{aligned}
\text{GPR-TS-MOR:}\hspB E_a(n=2,\nPsi = 2) = &1.2\times 10^{-2},\\
  \text{S-PROJ:}\hspB E_a(n=2,\nPsi = 2) =&9\times 10^{-2}.
\end{aligned}
\end{equation}
Clearly, GPR-TS-MOR outperforms S-PROJ. It results in an average error that is almost $5.5$ times smaller than that resulting from S-PROJ. \Cref{fig: test-3 error} depicts the $L^1$-error for several different parameter instances. For each of the tested parameter instances, the error from GPR-TS-MOR is five to ten times smaller than that from S-PROJ. \Cref{fig: test-3 err comp} compares the error surrogate to the true error and \Cref{fig: test-3 eff} presents the corresponding efficiency index. As before, on average, our surrogate accurately approximates the true error, with an efficiency index that oscillates between $0.6$ and $1.8$. 

\begin{figure}[ht!]
\centering
\subfloat[GPR-TS-MOR ($z=2.5\times 10^{-2}$)]{\includegraphics[width=2.5in]{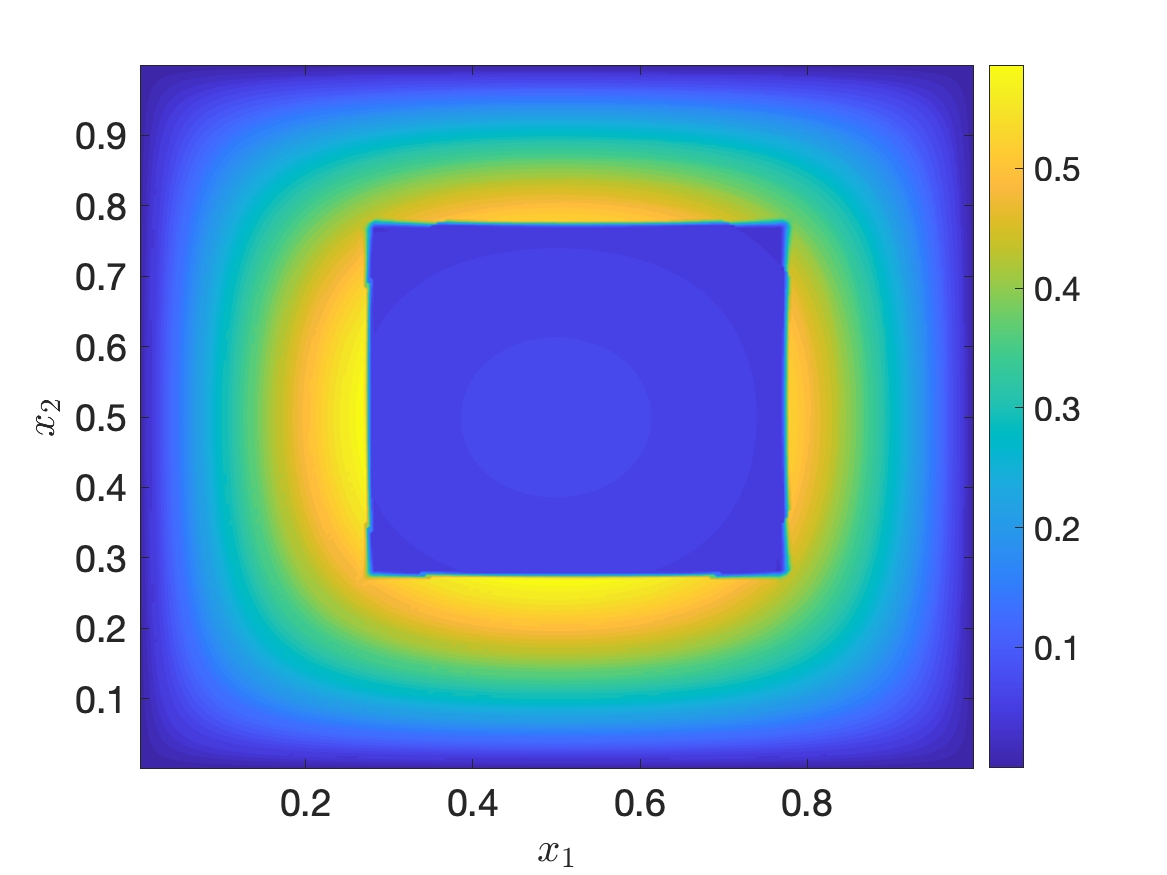}} 
\hfill
\subfloat[S-PROJ ($z=2.5\times 10^{-2}$)]{\includegraphics[width=2.5in]{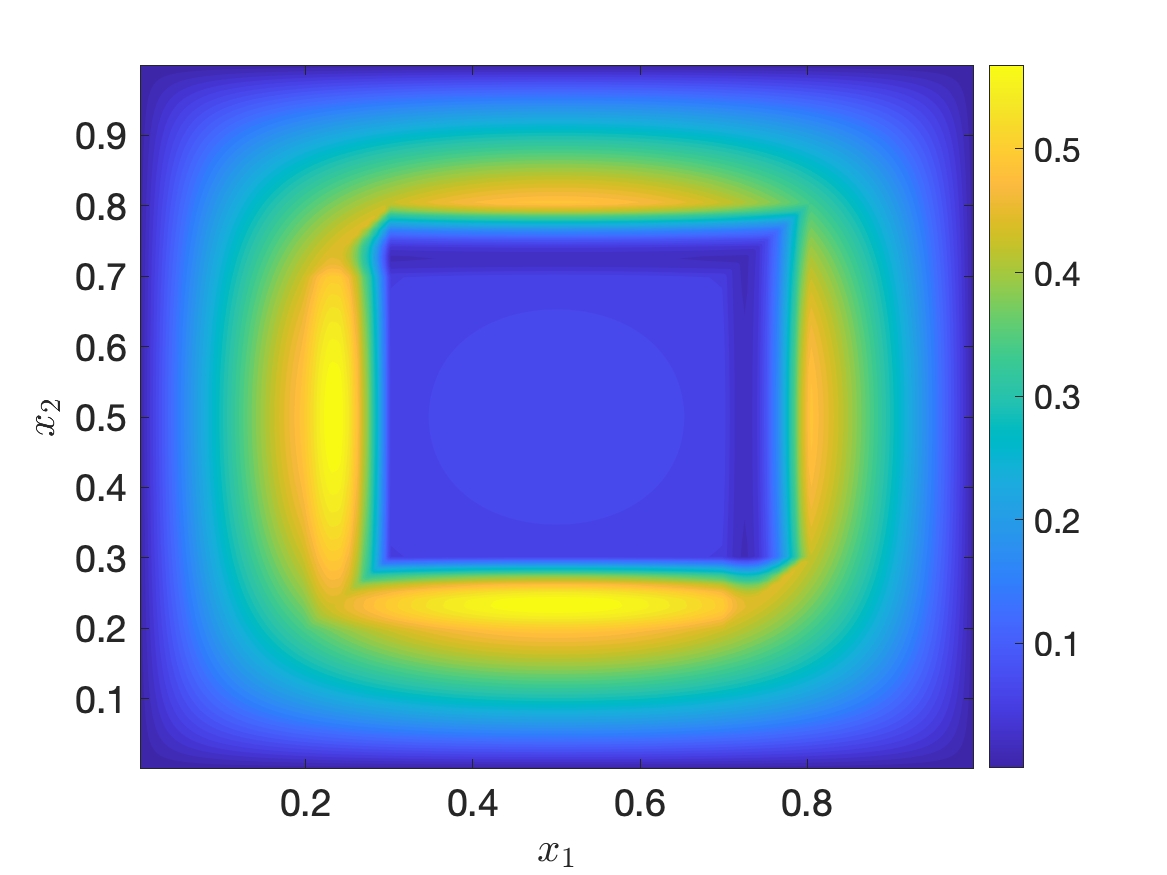}} 
\hfill
\subfloat[HF ($z=2.5\times 10^{-2}$)]{\includegraphics[width=2.5in]{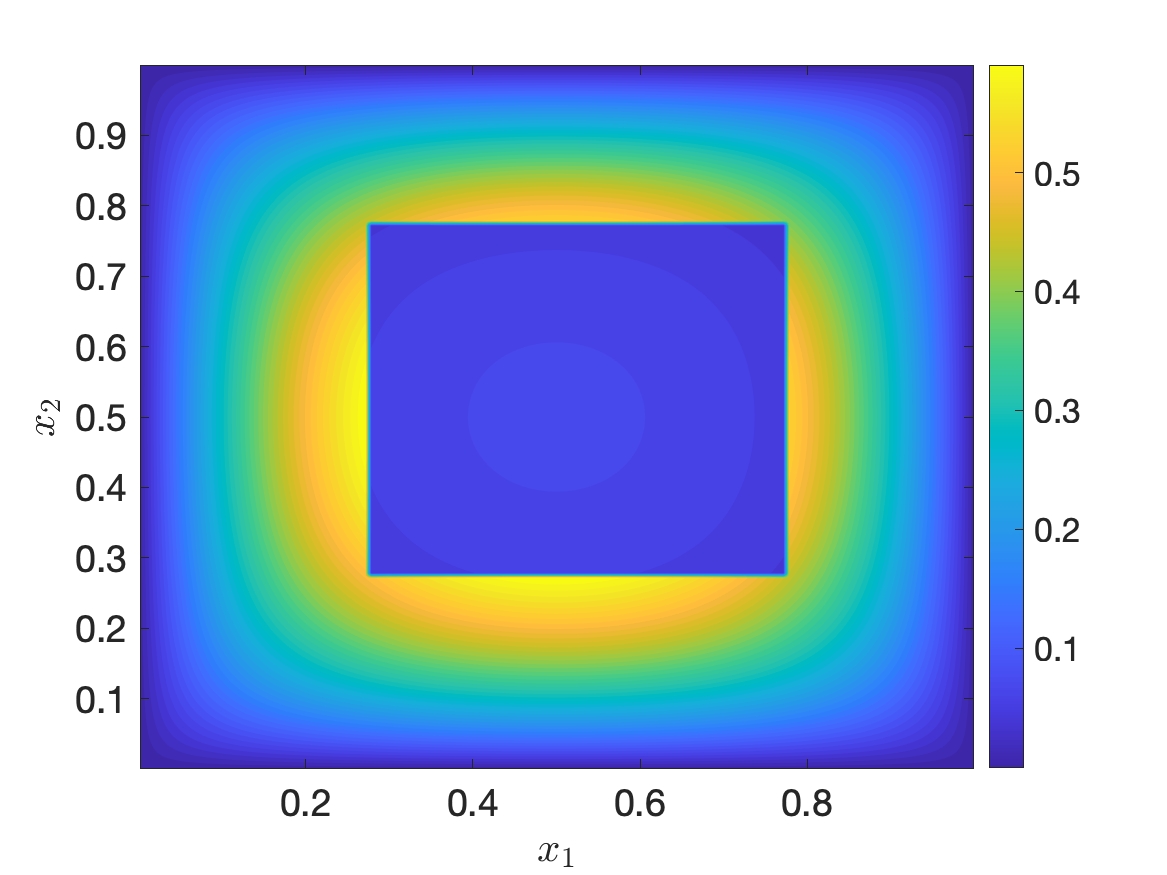}} 
\hfill 
\caption{\textit{Results for test-3. Comparison between the different solutions at $z=2.5\times 10^{-2}$.}}	\label{fig: test-3 sol comp}
\end{figure}

 For $z = 2.5\times 10^{-2}$, \Cref{fig: test-3 sol comp} compares the different solutions. With just two modes, GPR-TS-MOR provides an accurate approximation of the solution. Note that the solution from S-PROJ does not exhibit spurious oscillations reported in the previous test case because the value of $n$ is smaller---oscillations are only present in the higher order POD modes. However, close to the boundaries of the inner-box $\Omega_z$, it does exhibit a staircase type effect. This staircase effect is typical for linear reduced approximations of such problems---see \cite{Nair,Welper2017}, for further examples.
 
Consider the speed-up $\kappa$ defined in \eqref{def kappa}. For the current test case, we observed a speed of $51.5$ with an average error of $1.2\%$. Note that the speed-up is lower than in the previous test case because the HF solver of the current problem requires no time iterations.
\section{Conclusions and discussion} \label{sec: conclusion}
We have proposed a data-driven MOR technique to approximate parameterized partial differential equations that exhibit parameter-dependent jump-discontinuities. Our technique hinges on a two step procedure: transformation followed by de-transformation. The transformation step (ideally) removes the parametric discontinuities by composing the solution with a spatial transform. This results in a transformed solution that is well-approximable in a low-dimensional reduced space. After we approximate this transformed solution, we de-transform the approximation by composing it with an inverse of the spatial transform and recover an approximation to the solution of the differential equation. An offline-online paradigm based procedure guarantees an efficient transformation and de-transformation step. 

Two data-driven methodologies are the building-blocks of our MOR technique: (i) Gaussian process regression, and (ii) optimization-based image registration. With GPR, we approximate the map between the parameter domain and the expansion coefficients of the reduced basis. With image registration on the other hand, we compute the spatial transform that allows for the solution transformation. Owing to the GPR, our MOR technique is purely data-driven i.e., it does not even require the knowledge of the structural non-linearities in the differential equation. This way it treats linear, non-linear, affine-in-parameter and non-affine-in-parameter problems alike and doesn't rely on any hyper-reduction technique.


We performed numerical experiments involving hyperbolic and parabolic differential equations. We compared our technique to a standard MOR technique that does not perform any solution transformation. Main takeaways from our experiments are as follows. Firstly, our technique results in an almost oscillation free solution. We attribute this to the solution transformation that halts the movement of the discontinuities in the parameter domain, resulting in a well-behaved set of POD modes. In contrast, due to moving discontinuities, the standard technique results in a highly oscillatory solution. Secondly, for a given number of POD modes, we outperform the standard technique in terms of accuracy, with the error being two to ten times smaller than that resulting from the standard technique. Lastly, for moderate parameter dimensions, at least for the test cases we considered, we observed speed-ups of one to upto three orders-of-magnitude. 

\section*{Acknowledgements}
\addcontentsline{toc}{section}{Acknowledgments}
N.S  and P.B are supported by the German Federal Ministry for Economic Affairs and Energy (BMWi) in the joint project "MathEnergy - Mathematical Key Technologies for Evolving Energy Grids", sub-project: Model Order Reduction (Grant number: 0324019B).
 \bibliographystyle{abbrv}

\bibliography{papers}

\begin{thebibliography}{10}

\bibitem{NonlinearDMD}
A.~Alla and J.~N. Kutz.
\newblock Nonlinear model order reduction via dynamic mode decomposition.
\newblock {\em SIAM Journal on Scientific Computing}, 39(5):B778--B796, 2017.

\bibitem{GPRVectorValued}
M.~A. {\'A}lvarez, L.~Rosasco, and N.~D. Lawrence.
\newblock Kernels for vector-valued functions: A review.
\newblock {\em Foundations and Trends in Machine Learning}, 4(3):195--266,
  2012.

\bibitem{Amsallem2015}
D.~Amsallem, M.~Zahr, Y.~Choi, and C.~Farhat.
\newblock Design optimization using hyper-reduced-order models.
\newblock {\em Structural and Multidisciplinary Optimization}, 51(4):919--940,
  2015.

\bibitem{NairNonIntrusive2013}
C.~Audouze, F.~De~Vuyst, and P.~B. Nair.
\newblock Nonintrusive reduced-order modeling of parametrized time-dependent
  partial differential equations.
\newblock {\em Numerical Methods for Partial Differential Equations},
  29(5):1587--1628, 2013.

\bibitem{LDDM}
M.~F. Beg, M.~I. Miller, A.~Trouv{\'e}, and L.~Younes.
\newblock Computing large deformation metric mappings via geodesic flows of
  diffeomorphisms.
\newblock {\em International Journal of Computer Vision}, 61(2):139--157, 2005.

\bibitem{Kramer2020}
P.~Benner, P.~Goyal, B.~Kramer, B.~Peherstorfer, and K.~Willcox.
\newblock Operator inference for non-intrusive model reduction of systems with
  non-polynomial nonlinear terms.
\newblock {\em Computer Methods in Applied Mechanics and Engineering},
  372:113433, 2020.

\bibitem{PeterReview}
P.~Benner, S.~Gugercin, and K.~Willcox.
\newblock A survey of projection-based model reduction methods for parametric
  dynamical systems.
\newblock {\em SIAM Review}, 57(4):483--531, 2015.

\bibitem{PeterBook}
P.~Benner, M.~Ohlberger, A.~Cohen, and K.~Willcox.
\newblock {\em Model Reduction and Approximation}.
\newblock Society for Industrial and Applied Mathematics, Philadelphia, PA,
  2017.

\bibitem{Benner2014}
P.~Benner, E.~Sachs, and S.~Volkwein.
\newblock Model order reduction for {PDE} constrained optimization.
\newblock In G.~Leugering, P.~Benner, S.~Engell, A.~Griewank, H.~Harbrecht,
  M.~Hinze, R.~Rannacher, and S.~Ulbrich, editors, {\em Trends in PDE
  Constrained Optimization}, pages 303--326. Springer, 2014.

\bibitem{Benner2015uncertainty}
P.~Benner and J.~Schneider.
\newblock Uncertainty quantification for {M}axwell's equations using stochastic
  collocation and model order reduction.
\newblock {\em International Journal for Uncertainty Quantification},
  5(3):195--208, 2015.

\bibitem{Cagniart2019}
N.~Cagniart, Y.~Maday, and B.~Stamm.
\newblock Model order reduction for problems with large convection effects.
\newblock In B.~N. Chetverushkin, W.~Fitzgibbon, Y.~Kuznetsov,
  P.~Neittaanm{\"a}ki, J.~Periaux, and O.~Pironneau, editors, {\em
  Contributions to Partial Differential Equations and Applications}, pages
  131--150. Springer International Publishing, Cham, 2019.

\bibitem{Kevin2015adaptive}
K.~Carlberg.
\newblock Adaptive h-refinement for reduced-order models.
\newblock {\em International Journal for Numerical Methods in Engineering},
  102(5):1192--1210, 2015.

\bibitem{GreedyJan2018}
W.~Chen, J.~S. Hesthaven, B.~Junqiang, Y.~Qiu, Z.~Yang, and Y.~Tihao.
\newblock Greedy nonintrusive reduced order model for fluid dynamics.
\newblock {\em AIAA Journal}, 56(12):4927--4943, 2018.

\bibitem{Crisovan2019model}
R.~Crisovan, D.~Torlo, R.~Abgrall, and S.~Tokareva.
\newblock Model order reduction for parametrized nonlinear hyperbolic problems
  as an application to uncertainty quantification.
\newblock {\em Journal of Computational and Applied Mathematics}, 348:466--489,
  2019.

\bibitem{ROMES}
M.~Drohmann and K.~Carlberg.
\newblock The {ROMES} method for statistical modeling of reduced-order-model
  error.
\newblock {\em SIAM/ASA Journal on Uncertainty Quantification}, 3(1):116--145,
  2015.

\bibitem{OptimalSVD}
C.~Eckart and G.~Young.
\newblock The approximation of one matrix by another of lower rank.
\newblock {\em Psychometrika,}, 1:211–218, 1936.

\bibitem{Metric2019}
V.~Ehrlacher, D.~Lombardi, O.~Mula, and F.-X. Vialard.
\newblock Nonlinear model reduction on metric spaces. {A}pplication to
  one-dimensional conservative {PDE}s in {W}asserstein spaces.
\newblock {\em ESAIM: Mathematical Modelling and Numerical Analysis}, 2019.

\bibitem{WaveKolmogorov}
C.~Greif and K.~Urban.
\newblock Decay of the {K}olmogorov {N}-width for wave problems.
\newblock {\em Applied Mathematics Letters}, 96:216 -- 222, 2019.

\bibitem{Guo2019}
M.~Guo and J.~S. Hesthaven.
\newblock Data-driven reduced order modeling for time-dependent problems.
\newblock {\em Computer Methods in Applied Mechanics and Engineering},
  345:75--99, 2019.

\bibitem{Han2005}
J.~S. Han, E.~B. Rudnyi, and J.~G. Korvink.
\newblock Efficient optimization of transient dynamic problems in {MEMS}
  devices using model order reduction.
\newblock {\em Journal of Micromechanics and Microengineering}, 15(4):822--832,
  2005.

\bibitem{JanStamm}
J.~S. Hesthaven, G.~Rozza, and B.~Stamm.
\newblock {\em Certified Reduced Basis Methods for Parametrized Partial
  Differential Equations}.
\newblock Springer, Cham, 2016.

\bibitem{JanNN2018}
J.~S. Hesthaven and S.~Ubbiali.
\newblock Non-intrusive reduced order modeling of nonlinear problems using
  neural networks.
\newblock {\em Journal of Computational Physics}, 363:55--78, 2018.

\bibitem{GPRReview}
M.~Kanagawa, P.~Hennig, D.~Sejdinovic, and B.~K. Sriperumbudur.
\newblock Gaussian processes and kernel methods: A review on connections and
  equivalences.
\newblock {\em arXiv:1807.02582}, 2018.

\bibitem{Klein2007}
S.~{Klein}, M.~{Staring}, and J.~P.~W. {Pluim}.
\newblock Evaluation of optimization methods for nonrigid medical image
  registration using mutual information and {B}-{S}plines.
\newblock {\em IEEE Transactions on Image Processing}, 16(12):2879--2890, 2007.

\bibitem{KevinAuto}
K.~Lee and K.~T. Carlberg.
\newblock Model reduction of dynamical systems on nonlinear manifolds using
  deep convolutional autoencoders.
\newblock {\em Journal of Computational Physics}, 404:108973, 2020.

\bibitem{MovingBC2021}
Z.~Ma and W.~Pan.
\newblock Data-driven nonintrusive reduced order modeling for dynamical systems
  with moving boundaries using {G}aussian process regression.
\newblock {\em Computer Methods in Applied Mechanics and Engineering},
  373:113495, 2021.

\bibitem{SamplingLattice}
M.~D. McKay, R.~J. Beckman, and W.~J. Conover.
\newblock A comparison of three methods for selecting values of input variables
  in the analysis of output from a computer code.
\newblock {\em Technometrics}, 21(2):239--245, 1979.

\bibitem{MojganiLagrangian}
R.~Mojgani and M.~Balajewicz.
\newblock Lagrangian basis method for dimensionality reduction of convection
  dominated nonlinear flows.
\newblock {\em arXiv:1701.04343}, 2017.

\bibitem{Nair}
N.~J. Nair and M.~Balajewicz.
\newblock Transported snapshot model order reduction approach for parametric,
  steady-state fluid flows containing parameter-dependent shocks.
\newblock {\em International Journal for Numerical Methods in Engineering},
  117(12):1234--1262, 2019.

\bibitem{Noblet2008}
V.~Noblet, C.~Heinrich, F.~Heitz, and J.-P. Armspach.
\newblock Accurate inversion of 3-{D} transformation fields.
\newblock {\em IEEE transactions on image processing}, 17(10):1963--1968, 2008.

\bibitem{Benjamin2016}
B.~Peherstorfer and K.~Willcox.
\newblock Data-driven operator inference for nonintrusive projection-based
  model reduction.
\newblock {\em Computer Methods in Applied Mechanics and Engineering},
  306:196--215, 2016.

\bibitem{Qian2020}
E.~Qian, B.~Kramer, B.~Peherstorfer, and K.~Willcox.
\newblock Lift {\&} learn: Physics-informed machine learning for large-scale
  nonlinear dynamical systems.
\newblock {\em Physica D: Nonlinear Phenomena}, 406:132401, 2020.

\bibitem{RBBook}
A.~Quarteroni, A.~Manzoni, and F.~Negri.
\newblock {\em Reduced Basis Methods for Partial Differential Equations: An
  Introduction}.
\newblock Springer International Publishing, 2016.

\bibitem{GPR}
C.~E. Rasmussen and C.~K. Williams.
\newblock {\em Gaussian Processes in Machine Learning}.
\newblock MIT press Cambridge, 2006.

\bibitem{MATS}
D.~Rim, B.~Peherstorfer, and K.~T. Mandli.
\newblock Manifold approximations via transported subspaces: Model reduction
  for transport-dominated problems.
\newblock {\em arXiv:1912.13024}, 2019.

\bibitem{SarnaCalib2020}
N.~Sarna, J.~Giesselmann, and P.~Benner.
\newblock Data-driven snapshot calibration via monotonic feature matching.
\newblock {\em arXiv:2009.08414}, 2020.

\bibitem{ReviewDeformReg}
A.~{Sotiras}, C.~{Davatzikos}, and N.~{Paragios}.
\newblock Deformable medical image registration: A survey.
\newblock {\em IEEE Transactions on Medical Imaging}, 32(7):1153--1190, 2013.

\bibitem{Limiter1984}
P.~K. Sweby.
\newblock High resolution schemes using flux limiters for hyperbolic
  conservation laws.
\newblock {\em SIAM Journal on Numerical Analysis}, 21(5):995--1011, 1984.

\bibitem{RegisterMOR}
T.~Taddei.
\newblock A registration method for model order reduction: Data compression and
  geometry reduction.
\newblock {\em SIAM Journal on Scientific Computing}, 42(2):A997--A1027, 2020.

\bibitem{Taddei2020ST}
T.~Taddei and L.~Zhang.
\newblock Space-time registration-based model reduction of parameterized
  one-dimensional hyperbolic {PDE}s.
\newblock {\em arXiv:2004.06693}, 2020.

\bibitem{CurvedDomains}
T.~Taddei and L.~Zhang.
\newblock Registration-based model reduction in complex two-dimensional
  geometries.
\newblock {\em arXiv:2101.10259}, 2021.

\bibitem{MRAdetect2014}
M.~J. Vuik and J.~K. Ryan.
\newblock Multiwavelet troubled-cell indicator for discontinuity detection of
  discontinuous {G}alerkin schemes.
\newblock {\em Journal of Computational Physics}, 270:138--160, 2014.

\bibitem{WelperAdaptive}
G.~Welper.
\newblock $h$ and $hp$-adaptive interpolation by transformed snapshots for
  parametric and stochastic hyperbolic {PDE}s.
\newblock {\em arXiv:1710.11481}, 2017.

\bibitem{Welper2017}
G.~Welper.
\newblock Interpolation of functions with parameter dependent jumps by
  transformed snapshots.
\newblock {\em SIAM Journal on Scientific Computing}, 39(4):A1225--A1250, 2017.

\bibitem{Xiao2015}
D.~Xiao, F.~Fang, A.~Buchan, C.~Pain, I.~Navon, and A.~Muggeridge.
\newblock Non-intrusive reduced order modelling of the {N}avier–{S}tokes
  equations.
\newblock {\em Computer Methods in Applied Mechanics and Engineering},
  293:522--541, 2015.

\bibitem{RBFXiao2015}
D.~Xiao, F.~Fang, C.~Pain, and G.~Hu.
\newblock Non-intrusive reduced-order modelling of the {N}avier–{S}tokes
  equations based on {RBF} interpolation.
\newblock {\em International Journal for Numerical Methods in Fluids},
  79(11):580--595, 2015.

\bibitem{Xiao2016}
D.~Xiao, P.~Yang, F.~Fang, J.~Xiang, C.~Pain, and I.~Navon.
\newblock Non-intrusive reduced order modelling of fluid–structure
  interactions.
\newblock {\em Computer Methods in Applied Mechanics and Engineering},
  303:35--54, 2016.

\bibitem{Yue2013}
Y.~Yue and K.~Meerbergen.
\newblock Accelerating optimization of parametric linear systems by model order
  reduction.
\newblock {\em SIAM Journal on Optimization}, 23(2):1344--1370, 2013.

\end{thebibliography}

\end{document}